\documentclass[11pt]{article}%
\usepackage{amssymb}
\usepackage{amsfonts}
\usepackage{amsmath}
\usepackage{graphicx}%
\providecommand{\U}[1]{\protect\rule{.1in}{.1in}}
\newtheorem{theorem}{Theorem}

\newcounter{article}

\begin{document}

\title{The Mittag-Leffler function}
\author{P. Van Mieghem\thanks{ Delft University of Technology, Faculty of Electrical
Engineering, Mathematics and Computer Science, P.O Box 5031, 2600 GA Delft,
The Netherlands; \emph{email}: P.F.A.VanMieghem@tudelft.nl }}
\date{in honor of G\"{o}sta Mittag-Leffler\\
for Huijuan\\
25 September 2021}
\maketitle

\begin{abstract}
We review the function theoretical properties of the Mittag-Leffler function
$E_{a,b}\left(  z\right)  $ in a self-contained manner, but also add new
results; more than half is new!

\end{abstract}

\section{Introduction}

We investigate the Mittag-Leffler function,%
\begin{equation}
E_{a,b}\left(  z\right)  =\sum_{k=0}^{\infty}\frac{z^{k}}{\Gamma\left(
b+ak\right)  } \label{Mittag_Leffler_function_Eab}%
\end{equation}
introduced by G\"{o}sta Mittag-Leffler
\cite{Mittag_Leffler_1903_March,Mittag_Leffler_1903} in 1903 with $b=1$, which
he denoted as $E_{a}\left(  z\right)  =\sum_{k=0}^{\infty}\frac{z^{k}}%
{\Gamma\left(  1+ak\right)  }\overset{\Delta}{=}E_{a,1}\left(  z\right)  $.

We consider the broader definition $E_{a,b}\left(  z\right)  $ and not
$E_{a}\left(  z\right)  $, because the functional relations for $E_{a,b}%
\left(  z\right)  $ are closed and expressed in terms of $E_{a,b}\left(
z\right)  $, whereas confinement to $E_{a}\left(  z\right)  $ only, deprives
the analysis from a complete and more elegant picture. There exist
generalizations\footnote{The generalized Mittag-Leffler function is defined as
$E_{a,b,c}\left(  z\right)  =\sum_{k=0}^{\infty}\frac{\Gamma\left(
c+k\right)  }{\Gamma\left(  c\right)  }\frac{z^{k}}{\Gamma\left(  b+ak\right)
}$, where $E_{a,b,0}(z)=E_{a,b}\left(  z\right)  $.} of the Mittag-Leffler
function $E_{a,b}\left(  z\right)  $, which are beyond the present scope, but
discussed by Haubold \emph{et al}. \cite{Haubold_2011} and also covered in the
recent book by Gorenflo\emph{ et al.} \cite{Gorenflo_2020} on Mittag-Leffler
functions and their applications. The Mittag-Leffler function $E_{a}\left(
z\right)  $ is treated by Erdelyi \emph{et al.} \cite[Sec. 18.1 on p.
206-211]{Erdelyi_v3} and by Sansone and Gerretsen \cite[Sec. 6.13 on p.
345-348]{Sansone}. Our aim\footnote{During my sabattical at Stanford in 2015,
I encountered the rather exotic Mittag-Leffler distribution, which I intended
to include in my Performance Analysis book \cite{PVM_PAComplexNetsCUP} as
another example of a power law-like distribution. However, the functional
properties of the Mittag-Leffler function $E_{a,b}\left(  z\right)  $ require
attention first, before one can turn to probability theory.} here is to deduce
the most relevant functional properties of the Mittag-Leffler function
$E_{a,b}\left(  z\right)  $ defined in (\ref{Mittag_Leffler_function_Eab}).
Since about half of the results have been established before, the manuscript
in the form of articles (\textbf{art.}) as in our book \cite{PVM_graphspectra}%
, is more a review, without detailed historical citations as in
\cite{Gorenflo_2020}, but enriched with new results: \textbf{art.}
\ref{art_differentiation}, \ref{art_Hadamard_series}
\ref{art_Taylor_expansion_rond_z0}, \ref{art_Taylor_series_logE_a,b(z)},
\ref{art_Taylor_series_to_Ea,b}, \ref{art_Mobius_inversion},
\ref{art_Mertens_function_and_MittagLeffler},
\ref{art_integral_duplication_formula},
\ref{art_integral_multiplication_formula},
\ref{art_integral_multiplication_formula_special},
\ref{art_Beta_function_like_asymptotic_result_Gauss},
\ref{art_Apelblat_series}, \ref{art_Berberan_Santos_eigen_bewijs},
\ref{art_Euler-Maclaurin_summation}, \ref{art_Euler-Maclaurin_remainder},
\ref{art_1_op_sinpiw_integral_for_0<a<1}, \ref{art_integral_Q_a,b(u)},
\ref{art_bounds_integrals_0<a<1}, \ref{art_Cauchy_type_integral},
\ref{art_Mellin_transform_product_Gamma} and part of \textbf{art.
}\ref{art_Laplace_transform}, \ref{art_evaluation_basic_complex_integral},
\ref{art_complex_integral_cosec_type}, \ref{art_alternative_integrals_0<a<1},
\ref{art_monotonicity_Mittag_Leffler_function},
\ref{art_E_a(-x)_probability_theory}, \ref{art_prob_prop_Mittag-Leffler_rv}.
For self-consistency, we have included an appendix
\ref{sec_theory_Gamma_function} on the Gamma function $\Gamma\left(  z\right)
$, whose properties are essential for the Mittag-Leffler function
$E_{a,b}\left(  z\right)  $.

The Mittag-Leffler function $E_{a,b}\left(  z\right)  $ naturally arises in
fractional calculus and in the solution of \emph{fractional} order integral or
differential equations, as illustrated in \cite{Haubold_2011},
\cite{Gorenflo_2020} and \cite{Sandev_Tomovski}; for example, in the
fractional generalization of the heat equation, random walks, L\'{e}vy
flights, superdiffusive transport and viscoelasticity, and fractional Ohm's
Law. Abel's integral equation, whose solution involves $E_{a,b}\left(
z\right)  $, is treated in \cite[Chapter 7]{Gorenflo_2014}. A
\textquotedblleft fractional\textquotedblright\ generalization of the Poisson
renewal process, discussed in \cite[Sec. 9.4]{Gorenflo_2014}, consists of
replacing the exponential interarrival time between events by a Mittag-Leffler
distribution $E_{a,1}\left(  -t\right)  $ with real $t\geq0$ and $0<a\leq1$.
Section \ref{sec_Mittag_Leffler_probability_theory} presents in \textbf{art}.
\ref{art_monotonicity_Mittag_Leffler_function} a different proof of the
monotonicity of the Mittag-Leffler function $E_{a,b}\left(  z\right)  $, while
\textbf{art}. \ref{art_E_a(-x)_probability_theory} and \textbf{art}.
\ref{art_prob_prop_Mittag-Leffler_rv} focus on the Mittag-Leffler random
variable. The last section \ref{sec_integral_Iab(z)}, also new, explores the
related integral%
\begin{equation}
I_{a,b}\left(  z\right)  =\int_{0}^{\infty}\frac{z^{u}}{\Gamma\left(
b+au\right)  }\,du\label{def_integral_Ia,b(z)}%
\end{equation}
that naturally appears in the Euler-Maclaurin summation for $E_{a,b}\left(
z\right)  $ in \textbf{art}. \ref{art_Euler-Maclaurin_summation}.

\section{Complex function theory}

\label{sec_complex_function_theory}Since the Mittag-Leffler function
$E_{a,b}\left(  z\right)  $ is defined in (\ref{Mittag_Leffler_function_Eab})
as a power series in $z$, a first concern is its validity range in $z$. The
radius $R$ of convergence of the power series $f\left(  z\right)  =\sum
_{k=0}^{\infty}f_{k}\,z^{k}$ of a function $f$ satisfies $\frac{1}{R}=\lim
\sup_{k\rightarrow\infty}|f_{k}|^{1/k}$ or $\frac{1}{R}=\lim_{k\rightarrow
\infty}\left\vert \frac{f_{k+1}}{f_{k}}\right\vert $ when the latter exists
\cite{Titchmarshfunctions}. Using $\frac{\Gamma\left(  z+a\right)  }%
{\Gamma\left(  z+b\right)  }\sim z^{a-b}$ in \cite[6.1.47]{Abramowitz} when
$z\rightarrow\infty$, the radius of convergence of the power series in
(\ref{Mittag_Leffler_function_Eab}) is
\[
\frac{1}{R}=\lim_{k\rightarrow\infty}\left\vert \frac{\Gamma\left(
b+ak\right)  }{\Gamma\left(  b+a+ak\right)  }\right\vert =\lim_{k\rightarrow
\infty}|ak|^{-a}=0
\]
for $\operatorname{Re}\left(  a\right)  >0$ and all complex $b$. An entire
function has a power series with infinite radius of convergence and is, thus,
analytic in the entire complex plane with an essential singularity at
infinity. Associated to entire complex functions is the concept of the order
$\rho$, which is defined \cite[p. 248]{Titchmarshfunctions} for any
$\varepsilon>0$ as $f\left(  z\right)  =O\left(  e^{\left\vert z\right\vert
^{\rho+\varepsilon}}\right)  $ when $\left\vert z\right\vert \rightarrow
\infty$. A necessary and sufficient condition \cite[p. 253]%
{Titchmarshfunctions} that $f\left(  z\right)  =\sum_{k=0}^{\infty}%
f_{k}\,z^{k}$ should be an entire function of order $\rho$ is that%
\[
\frac{1}{\rho}=\lim_{k\rightarrow\infty}\frac{-\log\left\vert f_{k}\right\vert
}{k\log k}%
\]
Applied to the power series in (\ref{Mittag_Leffler_function_Eab}), after
invoking Stirling's asymptotic formula \cite[6.1.39]{Abramowitz} that follows
from (\ref{LogGamma(z)_bound_Re(z)>0}) in Appendix
\ref{sec_theory_Gamma_function},%
\[
\Gamma\left(  ak+b\right)  \sim\sqrt{2\pi}e^{-ak}\left(  ak\right)
^{ak+b-\frac{1}{2}}%
\]
yields%
\[
\frac{1}{\rho}=\lim_{k\rightarrow\infty}\frac{\log\left\vert \Gamma\left(
b+ak\right)  \right\vert }{k\log k}=\lim_{k\rightarrow\infty}\frac{\log
\sqrt{2\pi}-ak+\left(  ak+b-\frac{1}{2}\right)  \log\left(  ak\right)  }{k\log
k}=a
\]
Hence, we find one of the important properties that the Mittag-Leffler
function $E_{a,b}\left(  z\right)  $ is an entire complex function in $z$ of
order $\rho=\frac{1}{a}$ for $\operatorname{Re}\left(  a\right)  >0$ and any
$b$. For real $a$ and $b$, the Mittag-Leffler function $E_{a,b}\left(
z\right)  $ is real on the real axis. Moreover, for real positive $a$ and $b$,
the Mittag-Leffler function $E_{a,b}\left(  z\right)  $ attains its maximum on
the real positive axis, because $\left\vert E_{a,b}\left(  re^{i\theta
}\right)  \right\vert \leq\sum_{k=0}^{\infty}\frac{r^{k}\left\vert
e^{ik\theta}\right\vert }{\left\vert \Gamma\left(  b+ak\right)  \right\vert
}=\sum_{k=0}^{\infty}\frac{r^{k}}{\Gamma\left(  b+ak\right)  }=E_{a,b}\left(
r\right)  $.

The number $n\left(  r\right)  $ of zeros $z\,_{1},z_{2},\ldots$ of an entire
function \thinspace$f\left(  z\right)  $ of order $\rho$, for which
$\left\vert z_{n}\right\vert \leq r$ is non-decreasing in $r$, is $n\left(
r\right)  =O\left(  r^{\rho+\varepsilon}\right)  $. Roughly stated \cite[p.
249]{Titchmarshfunctions}, the higher the order $\rho$ of an entire function,
the more zeros it may have in a given region of the complex plane. Moreover,
if the modulus of a zero $z_{n}$ is $r_{n}=\left\vert z_{n}\right\vert $,
then
\[
\sum_{n=1}^{\infty}\frac{1}{\left(  r_{n}\right)  ^{\alpha}}\text{ converges
if }\alpha>\rho
\]
and the lower bound of $\alpha$ is the exponent $\rho_{1}$ of convergence;
thus $\rho_{1}\leq\rho$. If $\sum_{n=1}^{\infty}\left(  \frac{r}{r_{n}%
}\right)  ^{p+1}$ converges for an integer $p$, then $p+1>\rho_{1}$ and the
smallest integer $p$ is called the genus of $f\left(  z\right)  $. In any
case, $p\leq\rho_{1}\leq\rho$. We may have $\rho_{1}<\rho$, for example for
$f\left(  z\right)  =e^{z}$, whose order is $\rho=1$, but the exponent of
convergence $\rho_{1}=0$, because $e^{z}$ does not have zeros. Applied to the
Mittag-Leffler function $E_{a,b}\left(  z\right)  $ of order $\rho=\frac{1}%
{a}$, the theory indicates that more zeros are expected for small $a$ than for
large $a$, which seems contradictory to the monotonicity of $E_{a,b}\left(
z\right)  $ for $0<a<1$ on the negative real axis in \textbf{art}.
\ref{art_monotonicity_Mittag_Leffler_function}. The determination of the zeros
of $E_{a,b}\left(  z\right)  $ is generally difficult \cite[Sec.
4.6]{Gorenflo_2014},\cite{Popov_2013} and omitted here.

For some special values of the parameter $a$, the Taylor series
(\ref{Mittag_Leffler_function_Eab}) reduces to known functions, such as
$E_{0,b}\left(  z\right)  =\frac{1}{\Gamma\left(  b\right)  }\frac{1}{1-z}$
for $\left\vert z\right\vert <1$, $E_{1,1}\left(  z\right)  =e^{z}$ and
$E_{2,1}\left(  z\right)  =\cosh\left(  \sqrt{z}\right)  $. From the
incomplete Gamma function \cite[6.5.29]{Abramowitz}, we have
\begin{align}
E_{1,b}\left(  z\right)   &  =e^{z}\gamma^{\ast}\left(  b-1,z\right)
=\frac{e^{z}}{\Gamma\left(  b-1\right)  }\sum_{k=0}^{\infty}\frac{\left(
-z\right)  ^{k}}{k!\left(  b-1+k\right)  }\nonumber\\
&  =z^{1-b}e^{z}\left(  1-\frac{1}{\Gamma\left(  b-1\right)  }\int_{z}%
^{\infty}e^{-t}t^{b-2}dt\right)  \label{E_1,b}%
\end{align}
which can also be written in terms of Kummer's confluent hypergeometric
function \cite[13.1.2]{Abramowitz}
\[
M\left(  a,b;z\right)  =\frac{\Gamma\left(  b\right)  }{\Gamma\left(
a\right)  }\sum_{k=0}^{\infty}\frac{\Gamma\left(  a+k\right)  }{\Gamma\left(
b+k\right)  }\frac{z^{k}}{k!}%
\]
as%
\[
E_{1,b}\left(  z\right)  =\sum_{k=0}^{\infty}\frac{z^{k}}{\Gamma\left(
b+k\right)  }=\frac{M\left(  1,b,z\right)  }{\Gamma\left(  b\right)  }%
\]
Generalizations of (\ref{E_1,b}) to $E_{\frac{1}{n},b}\left(  z\right)  $ for
fractional $a=\frac{1}{n}$, where $n$ is an integer, are derived in
(\ref{Mittag_Leffler_1/n,b_Wiman}) in \textbf{art}. \ref{art_fractional_a}
below. Many analytic functions can be expressed in terms of the hypergeometric
function, defined by Gauss's series \cite[15.1.1]{Abramowitz}%
\begin{equation}
F\left(  a,b;c;z\right)  =\frac{\Gamma(c)}{\Gamma(a)\Gamma(b)}\sum
_{k=0}^{\infty}\frac{\Gamma(a+k)\Gamma(b+k)}{\Gamma(c+k)k!}z^{k}\hspace
{1cm}\text{convergent for }\left\vert z\right\vert <1
\label{Gauss_hypergeometric_series}%
\end{equation}
where the argument of the Gamma functions is of the form $\alpha k+\beta$ with
$\alpha=1$, in contrast to the Mittag-Leffler function in
(\ref{Mittag_Leffler_function_Eab}) where $\alpha=a$ is real positive. Just
the fact that $a$ is real and not an integer colors the theory of the
Mittag-Leffler function $E_{a,b}\left(  z\right)  $ and causes its main challenges.

\section{Deductions from the definition (\ref{Mittag_Leffler_function_Eab}) of
$E_{a,b}\left(  z\right)  $}

\medskip\medskip\refstepcounter{article}{\noindent\textbf{\thearticle. }%
}\ignorespaces\label{art_special_values_a_and_b} \emph{Special values of }%
$a$\emph{ and }$b$. If $b=0$, then, since $\lim_{z\rightarrow0}\frac{1}%
{\Gamma\left(  z\right)  }=0$, we have%
\[
E_{a,0}\left(  z\right)  =\sum_{k=1}^{\infty}\frac{z^{k}}{\Gamma\left(
ak\right)  }=\sum_{k=0}^{\infty}\frac{z^{k+1}}{\Gamma\left(  ak+a\right)  }%
\]
and, hence,%
\begin{equation}
E_{a,0}\left(  z\right)  =zE_{a,a}\left(  z\right)  \label{E_a,0}%
\end{equation}
If $a=0$, then%
\[
E_{0,b}\left(  z\right)  =\frac{1}{\Gamma\left(  b\right)  }\frac{1}%
{1-z}\hspace{1cm}\text{for }\left\vert z\right\vert <1
\]
From $E_{a,b}\left(  z\right)  =\frac{1}{\Gamma\left(  b\right)  }+\sum
_{k=1}^{\infty}\frac{z^{k}}{\Gamma\left(  b+ak\right)  }$, we observe that
$\lim_{a\rightarrow\infty}E_{a,b}\left(  z\right)  =\frac{1}{\Gamma\left(
b\right)  }$.

\medskip\refstepcounter{article}{\noindent\textbf{\thearticle. }%
}\ignorespaces\label{art_even_and_odd} After splitting odd and even indices in
the $k$-sum of (\ref{Mittag_Leffler_function_Eab}), we obtain%
\[
E_{a,b}\left(  -z\right)  =\sum_{k=0}^{\infty}\frac{\left(  -1\right)
^{k}z^{k}}{\Gamma\left(  b+ak\right)  }=\sum_{k=0}^{\infty}\frac{z^{2k}%
}{\Gamma\left(  b+2ak\right)  }-\sum_{k=0}^{\infty}\frac{z^{2k+1}}%
{\Gamma\left(  b+a+2ak\right)  }%
\]
and%
\begin{equation}
E_{a,b}\left(  -z\right)  =E_{2a,b}\left(  z^{2}\right)  -zE_{2a,b+a}\left(
z^{2}\right)  \label{Mittag_Leffler_arg_to_sqr(arg)}%
\end{equation}
Property (\ref{Mittag_Leffler_arg_to_sqr(arg)}) cannot be expressed for
$E_{a}\left(  -z\right)  $ in terms of itself and motivates our viewpoint that
the complex function theory of the Mittag-Leffler function should focus on
$E_{a,b}\left(  z\right)  $, rather than on $E_{a}\left(  z\right)  $. The
differentiation rule in \textbf{art}. \ref{art_Differentiation_recursion}
below is the more fundamental motivation.

Adding $E_{a,b}\left(  z\right)  =E_{2a,b}\left(  z^{2}\right)  +zE_{2a,b+a}%
\left(  z^{2}\right)  $ to (\ref{Mittag_Leffler_arg_to_sqr(arg)}) leads to%
\begin{equation}
E_{2a,b}\left(  z^{2}\right)  =\frac{E_{a,b}\left(  z\right)  +E_{a,b}\left(
-z\right)  }{2} \label{Sum_pos_neg_argument_E_a,b}%
\end{equation}
and, similarly,%
\begin{equation}
E_{2a,b+a}\left(  z^{2}\right)  =\frac{E_{a,b}\left(  z\right)  -E_{a,b}%
\left(  -z\right)  }{2z} \label{Verschil_pos_neg_argument_E_a,b}%
\end{equation}

\textbf{Examples} From $E_{1,1}\left(  z\right)  =E_{1}\left(  z\right)
=e^{z}$, the relation (\ref{Sum_pos_neg_argument_E_a,b}) indicates that
$E_{2}\left(  z\right)  =\cosh\left(  \sqrt{z}\right)  $ and next
\[
E_{4}\left(  z\right)  =\frac{1}{2}\left\{  \cosh\left(  z^{\frac{1}{4}%
}\right)  +\cos\left(  z^{\frac{1}{4}}\right)  \right\}
\]
The odd variant (\ref{Verschil_pos_neg_argument_E_a,b}) gives $E_{2,2}\left(
z^{2}\right)  =\frac{\sinh z}{z}$.

\medskip\refstepcounter{article}{\noindent\textbf{\thearticle. }%
}\ignorespaces\label{art_cyclotomic_property} \emph{Cyclotomic property}. When
introducing the identity $\sum_{r=0}^{m-1}e^{i\frac{2\pi kr}{m}}%
=\frac{1-e^{2\pi ki}}{1-e^{i\frac{2\pi k}{m}}}=m1_{m|k}$ into
(\ref{Mittag_Leffler_function_Eab}),%
\[
\sum_{r=0}^{m-1}E_{a,b}\left(  e^{i\frac{2\pi r}{m}}z\right)  =\sum
_{k=0}^{\infty}\frac{\sum_{r=0}^{m-1}e^{i\frac{2\pi kr}{m}}z^{k}}%
{\Gamma\left(  b+ak\right)  }=m\sum_{k=0}^{\infty}\frac{1_{m|k}z^{k}}%
{\Gamma\left(  b+ak\right)  }=m\sum_{l=0}^{\infty}\frac{z^{ml}}{\Gamma\left(
b+aml\right)  }%
\]
we obtain%
\begin{equation}
E_{am,b}\left(  z^{m}\right)  =\frac{1}{m}\sum_{r=0}^{m-1}E_{a,b}\left(
ze^{i\frac{2\pi r}{m}}\right)  \label{recursion_Mittag_Leffler_function_E_a,b}%
\end{equation}
For $m=2$ in (\ref{recursion_Mittag_Leffler_function_E_a,b}),
\[
2E_{2a,b}\left(  z^{2}\right)  =E_{a,b}\left(  z\right)  +E_{a,b}\left(
-z\right)
\]
we retrieve the relation equivalent to (\ref{Mittag_Leffler_arg_to_sqr(arg)}),
because (\ref{recursion_Mittag_Leffler_function_E_a,b}) essentially follows by
multisectioning of a power series \cite[Section 4.3]{Riordan} of which
splitting in odd and even terms is the obvious case in $m=2$ sections.

\textbf{Example} The case $a=b=1$ in
(\ref{recursion_Mittag_Leffler_function_E_a,b}),%
\begin{equation}
E_{m,1}\left(  z\right)  =E_{m}\left(  z\right)  =\frac{1}{m}\sum_{r=0}%
^{m-1}e^{z^{\frac{1}{m}}e^{i\frac{2\pi r}{m}}} \label{E_m,1=E_m}%
\end{equation}
can be extended to certain integer values of $b$. Indeed, using
\[
E_{1,n}\left(  z\right)  =\sum_{k=0}^{\infty}\frac{z^{k}}{\left(
k+n-1\right)  !}=\sum_{k=n-1}^{\infty}\frac{z^{k-\left(  n-1\right)  }}%
{k!}=z^{-\left(  n-1\right)  }\left(  e^{z}-\sum_{j=0}^{n-2}\frac{z^{j}}%
{j!}\right)
\]
in (\ref{recursion_Mittag_Leffler_function_E_a,b}) yields%
\[
E_{m,n}\left(  z\right)  =\frac{1}{m}\sum_{r=0}^{m-1}E_{1,n}\left(
z^{\frac{1}{m}}e^{i\frac{2\pi r}{m}}\right)  =\frac{z^{-\frac{n-1}{m}}}{m}%
\sum_{r=0}^{m-1}e^{-i\frac{2\pi r}{m}\left(  n-1\right)  }\left(
e^{z^{\frac{1}{m}}e^{i\frac{2\pi r}{m}}}-\sum_{j=0}^{n-2}\frac{z^{\frac{j}{m}%
}e^{i\frac{2\pi r}{m}j}}{j!}\right)
\]
The last double sum%
\[
\sum_{r=0}^{m-1}\sum_{j=0}^{n-2}e^{-i\frac{2\pi r}{m}\left(  n-1-j\right)
}\frac{z^{\frac{j}{m}}}{j!}=\sum_{q=1}^{n-1}\frac{z^{\frac{n-1-q}{m}}}{\left(
n-1-q\right)  !}\sum_{r=0}^{m-1}e^{-i\frac{2\pi qr}{m}}=\sum_{q=1}^{n-1}%
\frac{z^{\frac{n-1-q}{m}}m1_{m|q}}{\left(  n-1-q\right)  !}%
\]
vanishes for all integers $n\leq m$. Thus, for $0\leq n\leq m$, we arrive at%
\begin{equation}
E_{m,n}\left(  z\right)  =\frac{z^{-\frac{n-1}{m}}}{m}\sum_{r=0}%
^{m-1}e^{-i\frac{2\pi r}{m}\left(  n-1\right)  }e^{z^{\frac{1}{m}}%
e^{i\frac{2\pi r}{m}}} \label{Mittag_Leffler_integer_a_integer_b<a}%
\end{equation}

\medskip\refstepcounter{article}{\noindent\textbf{\thearticle. }%
}\ignorespaces\label{art_Mittag-Leffler_function_b=am} \emph{Mittag-Leffler
function with }$b=\beta+am$\emph{ where }$m\in\mathbb{Z}$. Another rewriting
of the definition (\ref{Mittag_Leffler_function_Eab}) of $E_{a,b}\left(
z\right)  $,%
\[
E_{a,b}\left(  z\right)  =\sum_{k=0}^{\infty}\frac{z^{k}}{\Gamma\left(
b+ak\right)  }=\sum_{k=1}^{\infty}\frac{z^{k-1}}{\Gamma\left(  b-a+ak\right)
}=\frac{1}{z}\left(  \sum_{k=0}^{\infty}\frac{z^{k}}{\Gamma\left(
b-a+ak\right)  }-\frac{1}{\Gamma\left(  b-a\right)  }\right)
\]
leads to \textquotedblleft the shift down of $b$ by $a$\textquotedblright%
\ formula%
\begin{equation}
E_{a,b}\left(  z\right)  =\frac{1}{z}\left(  E_{a,b-a}\left(  z\right)
-\frac{1}{\Gamma\left(  b-a\right)  }\right)  \label{E_a,b_in_E_a,b-a}%
\end{equation}
or, similarly after $b\rightarrow a+b$,%
\[
E_{a,b}\left(  z\right)  =\frac{1}{\Gamma\left(  b\right)  }+zE_{a,b+a}\left(
z\right)
\]
from which, for $b=0$, we find again (\ref{E_a,0}). If $b=\beta+ma$ in
(\ref{E_a,b_in_E_a,b-a}), then we obtain a recursion in $m$%
\begin{equation}
E_{a,\beta+ma}\left(  z\right)  =\frac{1}{z}\left(  E_{a,\beta+\left(
m-1\right)  a}\left(  z\right)  -\frac{1}{\Gamma\left(  \beta+\left(
m-1\right)  a\right)  }\right)  \label{recursion_E_a,b+ma}%
\end{equation}
After iteration of (\ref{recursion_E_a,b+ma}), we find%
\begin{equation}
z^{m}E_{a,\beta+ma}\left(  z\right)  =E_{a,\beta}\left(  z\right)  -\sum
_{l=0}^{m-1}\frac{z^{l}}{\Gamma\left(  \beta+la\right)  } \label{E_a,ma}%
\end{equation}
The case for $m=1$ in (\ref{E_a,ma}) is again (\ref{E_a,b_in_E_a,b-a}).
Relation (\ref{E_a,ma}) is directly retrieved from the definition
(\ref{Mittag_Leffler_function_Eab}) of $E_{a,b}\left(  z\right)  $,%
\[
E_{a,\beta+ma}\left(  z\right)  =\sum_{k=0}^{\infty}\frac{z^{k}}{\Gamma\left(
\beta+a\left(  m+k\right)  \right)  }=\sum_{k=m}^{\infty}\frac{z^{k-m}}%
{\Gamma\left(  \beta+ak\right)  }=z^{-m}\left(  \sum_{k=0}^{\infty}\frac
{z^{k}}{\Gamma\left(  \beta+ak\right)  }-\sum_{k=0}^{m-1}\frac{z^{k}}%
{\Gamma\left(  \beta+ak\right)  }\right)
\]
Similarly,%
\[
E_{a,\beta-ma}\left(  z\right)  =\sum_{k=0}^{\infty}\frac{z^{k}}{\Gamma\left(
\beta+a\left(  k-m\right)  \right)  }=\sum_{k=-m}^{\infty}\frac{z^{k+m}%
}{\Gamma\left(  \beta+ak\right)  }=z^{m}\left(  \sum_{k=0}^{\infty}\frac
{z^{k}}{\Gamma\left(  \beta+ak\right)  }+\sum_{k=1}^{m}\frac{z^{-k}}%
{\Gamma\left(  \beta-ak\right)  }\right)
\]
and, hence,%
\begin{equation}
z^{-m}E_{a,\beta-ma}\left(  z\right)  =E_{a,\beta}\left(  z\right)
+\sum_{k=1}^{m}\frac{\left(  \frac{1}{z}\right)  ^{k}}{\Gamma\left(
\beta-ak\right)  } \label{E_a,minus(ma)}%
\end{equation}

Subtracting (\ref{E_a,minus(ma)}) from (\ref{E_a,ma}) yields\footnote{which we
can write as%
\begin{align*}
-\sum_{l=-m}^{m-1}\frac{z^{l}}{\Gamma\left(  \beta+la\right)  }  &  =\int
_{-m}^{m}\frac{d}{du}\left\{  z^{u}E_{a,\beta+ua}\left(  z\right)  \right\}
du=\int_{-m}^{m}\sum_{k=0}^{\infty}\frac{d}{du}\frac{z^{k+u}}{\Gamma\left(
\beta+a\left(  u+k\right)  \right)  }du\\
&  =\log z\int_{-m}^{m}z^{u}E_{a,\beta+ua}\left(  z\right)  du-a\int_{-m}%
^{m}z^{u}\sum_{k=0}^{\infty}\frac{z^{k}\psi\left(  \beta+a\left(  u+k\right)
\right)  }{\Gamma\left(  \beta+a\left(  u+k\right)  \right)  }du
\end{align*}
where $\psi\left(  z\right)  =\frac{d\log\Gamma\left(  z\right)  }{dz}$ is the
digamma function \cite[6.3.1]{Abramowitz}.}%
\[
z^{m}E_{a,\beta+ma}\left(  z\right)  -z^{-m}E_{a,\beta-ma}\left(  z\right)
=-\sum_{l=-m}^{m-1}\frac{z^{l}}{\Gamma\left(  \beta+la\right)  }%
\]
illustrating for real $a>0$, $x\geq0$ and positive $\beta>ma$ that
$x^{2m}E_{a,\beta+ma}\left(  x\right)  <E_{a,\beta-ma}\left(  x\right)  $.

\medskip\refstepcounter{article}{\noindent\textbf{\thearticle. }%
}\ignorespaces\label{art_differentiation} \emph{Differentiation with respect
to }$z$. The derivative of the definition (\ref{Mittag_Leffler_function_Eab})
with respect to $z$ is
\begin{align*}
\frac{d}{dz}E_{a,b}\left(  z\right)   &  =\sum_{k=1}^{\infty}\frac{kz^{k-1}%
}{\Gamma\left(  b+ak\right)  }=\sum_{k=0}^{\infty}\frac{\left(  k+1\right)
z^{k}}{\Gamma\left(  b+a+ak\right)  }\\
&  =\frac{1}{a}\sum_{k=0}^{\infty}\frac{\left(  ak+b+a-1\right)  z^{k}}%
{\Gamma\left(  b+a+ak\right)  }-\frac{b-1}{a}\sum_{k=0}^{\infty}\frac{z^{k}%
}{\Gamma\left(  b+a+ak\right)  }\\
&  =\frac{1}{a}\sum_{k=0}^{\infty}\frac{z^{k}}{\Gamma\left(  b+a-1+ak\right)
}-\frac{b-1}{a}\sum_{k=0}^{\infty}\frac{z^{k}}{\Gamma\left(  b+a+ak\right)  }%
\end{align*}
and, with the definition (\ref{Mittag_Leffler_function_Eab}),%
\[
\frac{d}{dz}E_{a,b}\left(  z\right)  =\frac{1}{a}E_{a,b+a-1}\left(  z\right)
-\frac{b-1}{a}E_{a,b+a}\left(  z\right)
\]
Using (\ref{E_a,b_in_E_a,b-a}), $E_{a,b+a}\left(  z\right)  =\frac{1}%
{z}\left(  E_{a,b}\left(  z\right)  -\frac{1}{\Gamma\left(  b\right)
}\right)  $, simplifies to%
\begin{equation}
az\frac{d}{dz}E_{a,b}\left(  z\right)  =E_{a,b-1}\left(  z\right)  -\left(
b-1\right)  E_{a,b}\left(  z\right)  \label{D_E_a,b(z)}%
\end{equation}

The $m$-derivative%
\[
a\frac{d^{m}}{dz^{m}}E_{a,b}\left(  z\right)  =\frac{d^{m-1}}{dz^{m-1}%
}E_{a,b+a-1}\left(  z\right)  -\left(  b-1\right)  \frac{d^{m-1}}{dz^{m-1}%
}E_{a,b+a}\left(  z\right)
\]
has a recursive structure, when denoting $h_{m}\left(  b\right)  =\frac{d^{m}%
}{dz^{m}}E_{a,b}\left(  z\right)  $ and $h_{0}\left(  b\right)  =E_{a,b}%
\left(  z\right)  $,
\[
ah_{m}\left(  b\right)  =h_{m-1}\left(  b+a-1\right)  -\left(  b-1\right)
h_{m-1}\left(  b+a\right)
\]
which can be iterated resulting in%
\[
a^{m}\frac{d^{m}}{dz^{m}}E_{a,b}\left(  z\right)  =\sum_{j=0}^{m}q_{j}\left(
a,b,m\right)  E_{a,b+ma-j}\left(  z\right)
\]
where the coefficients $q_{m}\left(  a,b,m\right)  =1$, $q_{m-1}\left(
a,b,m\right)  =-\left(  mb+\frac{m\left(  m-1\right)  }{2}a-\frac{m\left(
m+1\right)  }{2}\right)  $ and $q_{0}\left(  a,b,m\right)  =\prod_{k=0}%
^{m-1}\left(  b-1+ka\right)  $. In general, $q_{j}\left(  a,b,m\right)  $ are
polynomials in $a$ and $b$ of order $m-j$. Unfortunately\footnote{For example,%
\begin{align*}
a^{2}\frac{d^{2}E_{a,b}\left(  z\right)  }{dz^{2}}  &  =E_{a,b+2a-2}\left(
z\right)  -\left(  2b+a-3\right)  E_{a,b+2a-1}\left(  z\right)  +\left(
b-1\right)  \left(  b+a-1\right)  E_{a,b+2a}\left(  z\right) \\
a^{3}\frac{d^{3}E_{a,b}\left(  z\right)  }{dz^{3}}  &  =E_{a,b+3a-3}\left(
z\right)  -\left(  3b+3a-6\right)  E_{a,b+3a-2}\left(  z\right) \\
&  +\left\{  \left(  2b+a-3\right)  \left(  b+2a-2\right)  +\left(
b-1\right)  \left(  b+a-1\right)  \right\}  E_{a,b+3a-1}\left(  z\right)
-\left(  b-1\right)  \left(  b-1+a\right)  \left(  b-1+2a\right)
E_{a,b+3a}\left(  z\right) \\
a^{4}\frac{d^{4}E_{a,b}\left(  z\right)  }{dz^{4}}  &  =E_{a,b+4a-4}\left(
z\right)  -\left(  4b+6a-10\right)  E_{a,b+4a-3}\left(  z\right) \\
&  +\left\{  \left(  2b+a-3\right)  \left(  b+2a-2\right)  +\left(
b-1\right)  \left(  b+a-1\right)  +3\left(  b+a-2\right)  \left(
b+3a-3\right)  \right\}  E_{a,b+4a-2}\left(  z\right) \\
&  -\left\{  \left(  b+3a-2\right)  \left(  2b+a-3\right)  \left(
b+2a-2\right)  +\left(  b-1\right)  \left(  b+a-1\right)  \left(
b+3a-2\right)  +\left(  b-1\right)  \left(  b-1+a\right)  \left(
b-1+2a\right)  \right\}  E_{a,b+4a-1}\left(  z\right) \\
&  +\left(  b-1\right)  \left(  b-1+a\right)  \left(  b-1+2a\right)  \left(
b-1+3a\right)  E_{a,b+4a}\left(  z\right)
\end{align*}
}, it is not easy to write all coefficients $q_{j}\left(  a,b,m\right)  $ in
closed form. With (\ref{E_a,ma}), we have%
\begin{equation}
\left(  az\right)  ^{m}\frac{d^{m}}{dz^{m}}E_{a,b}\left(  z\right)
=\sum_{j=0}^{m}q_{j}\left(  a,b,m\right)  \left(  E_{a,b-j}\left(  z\right)
-\sum_{l=0}^{m-1}\frac{z^{l}}{\Gamma\left(  b-j+la\right)  }\right)
\label{m_th_derivative_E_a,b(z)}%
\end{equation}
For $z=0$, we obtain from (\ref{m_th_derivative_E_a,b(z)}) with $\left.
\frac{d^{m}}{dz^{m}}E_{a,b}\left(  z\right)  \right\vert _{z=0}=\frac
{1}{\Gamma\left(  b+am\right)  }$,%
\[
\frac{1}{\Gamma\left(  b+am\right)  }=\frac{1}{a^{m}}\sum_{j=0}^{m}\frac
{q_{j}\left(  a,b,m\right)  }{\Gamma\left(  b-j+ma\right)  }%
\]
and, with $\frac{\Gamma\left(  b+am\right)  }{\Gamma\left(  b-j+ma\right)
}=\prod_{q=1}^{j}\left(  b+am-q\right)  $, the polynomial nature of
$q_{j}\left(  a,b,m\right)  $ is apparent:
\[
\sum_{j=1}^{m}\prod_{q=1}^{j}\left(  b+am-q\right)  q_{j}\left(  a,b,m\right)
=a^{m}-\prod_{k=0}^{m-1}\left(  b-1+ka\right)
\]
\textbf{Art}. \ref{art_Taylor_expansion_rond_z0} below will present a closed
form for $\frac{d^{m}}{dz^{m}}E_{a,b}\left(  z\right)  $.

\medskip\refstepcounter{article}{\noindent\textbf{\thearticle. }%
}\ignorespaces\label{art_Differentiation_recursion} \emph{Differentiation
recursion}. Using the functional equation of the Gamma function,
$\Gamma\left(  z+1\right)  =z\Gamma\left(  z\right)  $, we write%
\[
E_{a,b}\left(  xz\right)  =\sum_{k=0}^{\infty}\frac{x^{k}z^{k}}{\Gamma\left(
b+ak\right)  }=\sum_{k=0}^{\infty}\frac{x^{k}z^{k}}{\left(  b-1+ak\right)
\Gamma\left(  b-1+ak\right)  }%
\]
Thus,%
\[
z^{b-1}E_{a,b}\left(  xz^{a}\right)  =\sum_{k=0}^{\infty}\frac{x^{k}%
z^{b-1+ak}}{\left(  b-1+ak\right)  \Gamma\left(  b-1+ak\right)  }%
\]
Differentiation with respect to $z$ gives us%
\[
\frac{d}{dz}\left\{  z^{b-1}E_{a,b}\left(  xz^{a}\right)  \right\}
=\sum_{k=0}^{\infty}\frac{x^{k}z^{b-2+ak}}{\Gamma\left(  b-1+ak\right)
}=z^{b-2}\sum_{k=0}^{\infty}\frac{\left(  xz^{a}\right)  ^{k}}{\Gamma\left(
b-1+ak\right)  }%
\]
With the definition (\ref{Mittag_Leffler_function_Eab}), we arrive for any $x$
at the differentiation recursion in $b$,%
\begin{equation}
\frac{d}{dz}\left\{  z^{b-1}E_{a,b}\left(  xz^{a}\right)  \right\}
=z^{b-2}E_{a,b-1}\left(  xz^{a}\right)  \label{D_(z^(b-1)E_a,b(z^a))}%
\end{equation}
Differentiating (\ref{D_(z^(b-1)E_a,b(z^a))}) again $m$-times and using the
recursion (\ref{D_(z^(b-1)E_a,b(z^a))}) yields
\begin{equation}
\frac{d^{m}}{dz^{m}}\left\{  z^{b-1}E_{a,b}\left(  xz^{a}\right)  \right\}
=z^{b-1-m}E_{a,b-m}\left(  xz^{a}\right)  \label{Dm_(z^(b-1)E_a,b(z^a))}%
\end{equation}

\medskip\refstepcounter{article}{\noindent\textbf{\thearticle. }%
}\ignorespaces\label{art_fractional_a} \emph{Fractional values of }$a$. Let
$a=\frac{m}{n}$ where $m\leq n$ are integers, then the differentiation formula
(\ref{Dm_(z^(b-1)E_a,b(z^a))}) becomes, for $x=1$,%
\begin{align*}
\frac{d^{m}}{dz^{m}}\left\{  z^{b-1}E_{\frac{m}{n},b}\left(  z^{\frac{m}{n}%
}\right)  \right\}   &  =z^{b-1}\sum_{k=0}^{\infty}\frac{z^{\frac{m}{n}k-m}%
}{\Gamma\left(  b-m+\frac{m}{n}k\right)  }=z^{b-1}\sum_{k=0}^{\infty}%
\frac{z^{\frac{m}{n}\left(  k-n\right)  }}{\Gamma\left(  b+\frac{m}{n}\left(
k-n\right)  \right)  }\\
&  =z^{b-1}\sum_{k=-n}^{\infty}\frac{z^{\frac{m}{n}k}}{\Gamma\left(
b+\frac{m}{n}k\right)  }=z^{b-1}\sum_{k=1}^{n}\frac{z^{-\frac{m}{n}k}}%
{\Gamma\left(  b-\frac{m}{n}k\right)  }+z^{b-1}\sum_{k=0}^{\infty}%
\frac{z^{\frac{m}{n}k}}{\Gamma\left(  b+\frac{m}{n}k\right)  }%
\end{align*}
and we find that%
\begin{equation}
\frac{d^{m}}{dz^{m}}\left\{  z^{b-1}E_{\frac{m}{n},b}\left(  z^{\frac{m}{n}%
}\right)  \right\}  =z^{b-1}\sum_{k=1}^{n}\frac{z^{-\frac{m}{n}k}}%
{\Gamma\left(  b-\frac{m}{n}k\right)  }+z^{b-1}E_{\frac{m}{n},b}\left(
z^{\frac{m}{n}}\right)  \label{Dm_(z^(b-1)E_a,b(z^a))_fractional_a}%
\end{equation}

For $n=1$, (\ref{Dm_(z^(b-1)E_a,b(z^a))_fractional_a}) is%
\[
\frac{d^{m}}{dz^{m}}\left\{  z^{b-1}E_{m,b}\left(  z^{m}\right)  \right\}
=\frac{z^{b-1-m}}{\Gamma\left(  b-m\right)  }+z^{b-1}E_{m,b}\left(
z^{m}\right)
\]
In particular, when $b=1$, then it holds for $m\geq1$ that $\frac{d^{m}%
}{dz^{m}}E_{m}\left(  z^{m}\right)  =E_{m}\left(  z^{m}\right)  $, which
illustrates that $y=E_{m}\left(  z^{m}\right)  $, explicitly given in
(\ref{E_m,1=E_m}), is a solution of the differential equation $\frac{d^{m}%
y}{dz^{m}}=y$.

For $m=1$, (\ref{Dm_(z^(b-1)E_a,b(z^a))_fractional_a}) reduces to%
\[
\frac{d}{dz}\left\{  z^{b-1}E_{\frac{1}{n},b}\left(  z^{\frac{1}{n}}\right)
\right\}  =\sum_{k=1}^{n}\frac{z^{b-1-\frac{1}{n}k}}{\Gamma\left(  b-\frac
{1}{n}k\right)  }+z^{b-1}E_{\frac{1}{n},b}\left(  z^{\frac{1}{n}}\right)
\]
from which%
\[
\frac{d}{dz}\left\{  e^{-z}z^{b-1}E_{\frac{1}{n},b}\left(  z^{\frac{1}{n}%
}\right)  \right\}  =e^{-z}\sum_{k=1}^{n}\frac{z^{b-1-\frac{1}{n}k}}%
{\Gamma\left(  b-\frac{1}{n}k\right)  }%
\]
Integrating both sides from $0$ to $z$ yields, for $b\geq1$,%
\begin{equation}
E_{\frac{1}{n},b}\left(  z^{\frac{1}{n}}\right)  =z^{1-b}e^{z}\left\{
1_{\left\{  b=1\right\}  }+\sum_{k=1}^{n}\frac{1}{\Gamma\left(  b-\frac{1}%
{n}k\right)  }\int_{0}^{z}t^{b-1-\frac{1}{n}k}e^{-t}dt\right\}
\label{Mittag_Leffler_1/n,b}%
\end{equation}
where the indicator function $1_{x}$ equals one if the condition $x$ is true,
else $1_{x}=0$. After letting $x=z^{\frac{1}{n}}$ in
(\ref{Mittag_Leffler_1/n,b}) and replacing $j=n-k$ in the summation, we
obtain\footnote{Using the recursion $P\left(  a,x\right)  =P\left(
a+1,x\right)  +\frac{x^{a}}{\Gamma\left(  a+1\right)  }e^{-x}$ in
\cite[6.5.21]{Abramowitz} for the incomplete Gamma function $P\left(
a,x\right)  =\frac{1}{\Gamma\left(  a\right)  }\int_{0}^{x}t^{a-1}e^{-t}dt$ in
(\ref{Mittag_Leffler_1/n,b_Wiman}), $E_{\frac{1}{n},b}\left(  x\right)
=x^{\left(  1-b\right)  n}e^{x^{n}}\left\{  1_{\left\{  b=1\right\}  }%
+\sum_{j=0}^{n-1}P\left(  b-1+\frac{j}{n},x^{n}\right)  \right\}  $, shows
that%
\[
E_{\frac{1}{n},b}\left(  x\right)  =\sum_{j=0}^{\infty}\frac{x^{j}}%
{\Gamma\left(  b+\frac{j}{n}\right)  }=x^{\left(  1-b\right)  n}e^{x^{n}%
}\left\{  1_{\left\{  b=1\right\}  }+\sum_{j=0}^{n-1}P\left(  b+\frac{j}%
{n},x^{n}\right)  \right\}  +\sum_{j=0}^{n-1}\frac{x^{j}}{\Gamma\left(
b+\frac{j}{n}\right)  }%
\]
which agrees with the Wiman-recursion (\ref{Mittag_Leffler_1/n,b_Wiman}) when
$b\rightarrow b+1$,
\[
\sum_{j=n}^{\infty}\frac{x^{j}}{\Gamma\left(  b+\frac{j}{n}\right)  }%
=x^{n}E_{\frac{1}{n},b+1}\left(  x\right)  =x^{n}x^{\left(  1-\left(
b+1\right)  \right)  n}e^{x^{n}}\left\{  1_{\left\{  b=0\right\}  }+\sum
_{j=0}^{n-1}P\left(  b+\frac{j}{n},x^{n}\right)  \right\}
\]
}
\begin{equation}
E_{\frac{1}{n},b}\left(  x\right)  =x^{\left(  1-b\right)  n}e^{x^{n}}\left\{
1_{\left\{  b=1\right\}  }+\sum_{j=0}^{n-1}\frac{1}{\Gamma\left(  b-1+\frac
{j}{n}\right)  }\int_{0}^{x^{n}}t^{\left(  b-1+\frac{j}{n}\right)  -1}%
e^{-t}dt\right\}  \label{Mittag_Leffler_1/n,b_Wiman}%
\end{equation}
that reduces, for $b=1$, to Wiman's form in \cite{Wiman_1905}%
\[
E_{\frac{1}{n}}\left(  x\right)  =e^{x^{n}}\left\{  1+\int_{0}^{x^{n}}%
e^{-t}\sum_{j=1}^{n-1}\frac{t^{\frac{j}{n}-1}}{\Gamma\left(  \frac{j}%
{n}\right)  }dt\right\}
\]
For positive real $x$, it holds that $\int_{0}^{x^{n}}t^{\left(  b-1+\frac
{j}{n}\right)  -1}e^{-t}dt<\int_{0}^{\infty}t^{\left(  b-1+\frac{j}{n}\right)
-1}e^{-t}dt=\Gamma\left(  b-1+\frac{j}{n}\right)  $ and
(\ref{Mittag_Leffler_1/n,b_Wiman}) shows that $E_{\frac{1}{n},b}\left(
x\right)  <x^{\left(  1-b\right)  n}e^{x^{n}}\left\{  1_{\left\{  b=1\right\}
}+n-1_{\left\{  b=1\right\}  }\right\}  $ and%
\[
E_{\frac{1}{n},b}\left(  x\right)  <nx^{\left(  1-b\right)  n}e^{x^{n}}%
\]
which illustrates for $a=\frac{1}{n}$ that the entire function $E_{\frac{1}%
{n},b}\left(  z\right)  $ has order $\rho=\frac{1}{a}=n$. This bound also
reappears in \textbf{art}.
\ref{art_deductions_Mittag_Leffler_contour_integral}. Moreover, it is
interesting to compare (\ref{Mittag_Leffler_integer_a_integer_b<a}) formally,
written for an integer $b<m$,%

\[
E_{m,b}\left(  z\right)  =\frac{z^{\frac{1-b}{m}}}{m}\left\{  e^{z^{\frac
{1}{m}}}+\sum_{r=1}^{m-1}e^{-i\frac{2\pi r}{m}\left(  b-1\right)  }%
e^{z^{\frac{1}{m}}e^{i\frac{2\pi r}{m}}}\right\}
\]
with Bieberbach's integral (\ref{E_a,b(z)_Bieberbach_growth_z}) for
$E_{a,b}\left(  z\right)  $ and with the bound, where $\frac{1}{n}$ is
replaced by $m$,%
\[
E_{m,b}\left(  z\right)  <\frac{z^{\frac{1-b}{m}}}{m}e^{z^{\frac{1}{m}}}%
\]
We return to the relation between $E_{a,b}\left(  z\right)  $ and $E_{\frac
{1}{a},b}\left(  z\right)  $ later in \textbf{art}.
\ref{art_deductions_Mittag_Leffler_contour_integral}.

\textbf{Examples} If $n=1$ in (\ref{Mittag_Leffler_1/n,b}), we retrieve
(\ref{E_1,b}) and when $n=2$ in (\ref{Mittag_Leffler_1/n,b}), we find for
$b=1$%
\[
E_{\frac{1}{2}}\left(  z^{\frac{1}{2}}\right)  =e^{z}\left\{  1+\frac
{1}{\Gamma\left(  \frac{1}{2}\right)  }\int_{0}^{z}t^{-\frac{1}{2}}%
e^{-t}dt\right\}  =e^{z}\left\{  1+\frac{2}{\sqrt{\pi}}\int_{0}^{\sqrt{z}%
}e^{-u^{2}}du\right\}
\]
so that, with the definition \cite[7.1.1]{Abramowitz} of the error function
$\operatorname{erf}\left(  x\right)  =\frac{2}{\sqrt{\pi}}\int_{0}%
^{z}e^{-u^{2}}du$,%
\begin{equation}
E_{\frac{1}{2}}\left(  z\right)  =e^{z^{2}}\left\{  1+\frac{2}{\sqrt{\pi}}%
\int_{0}^{z}e^{-u^{2}}du\right\}  =e^{z^{2}}\left\{  1+\operatorname{erf}%
\left(  z\right)  \right\}  \label{E_1/2,1}%
\end{equation}

\medskip\refstepcounter{article}{\noindent\textbf{\thearticle. }%
}\ignorespaces\label{art_Hadamard_series} \emph{Hadamard's series }$\sum
_{k=0}^{\infty}\frac{z^{k}}{\left(  k!\right)  ^{a}}$\emph{ and }$E_{a}\left(
z\right)  $.\emph{ }Sharp bounds of the Gamma function for $ak+b>0$, that
follow from (\ref{LogGamma(z)_bound_Re(z)>0}) in Appendix
\ref{sec_theory_Gamma_function}, are%
\[
\sqrt{2\pi}\left(  ak+b\right)  ^{ak+b-\frac{1}{2}}e^{-\left(  ak+b\right)
}\leq\Gamma\left(  ak+b\right)  \leq\sqrt{2\pi}\left(  ak+b\right)
^{ak+b-\frac{1}{2}}e^{-\left(  ak+b\right)  }e^{\frac{1}{12\left(
ak+b\right)  }}%
\]
Now,%
\begin{align*}
\sqrt{2\pi}\left(  ak+b\right)  ^{ak+b-\frac{1}{2}}e^{-\left(  ak+b\right)  }
&  =\sqrt{2\pi}a^{a\left(  k+\frac{b}{a}\right)  -\frac{1}{2}}\left(
k+\frac{b}{a}\right)  ^{\frac{a-1}{2}}\left(  \left(  k+\frac{b}{a}\right)
^{\left(  k+\frac{b}{a}-\frac{1}{2}\right)  }e^{-\left(  k+\frac{b}{a}\right)
}\right)  ^{a}\\
&  \leq\sqrt{2\pi}a^{a\left(  k+\frac{b}{a}\right)  -\frac{1}{2}}\left(
k+\frac{b}{a}\right)  ^{\frac{a-1}{2}}\left(  \frac{\Gamma\left(  k+\frac
{b}{a}\right)  }{\sqrt{2\pi}}\right)  ^{a}%
\end{align*}
and replacing the inequality, we approximatively have%
\[
\Gamma\left(  ak+b\right)  \approx\frac{\left(  \sqrt{2\pi}\right)
^{1-a}a^{b}}{\sqrt{ak+b}}\left(  a^{k}\sqrt{k+\frac{b}{a}}\Gamma\left(
k+\frac{b}{a}\right)  \right)  ^{a}%
\]
If $b=1$, then for $ak+1>0$, the above reduces to%
\begin{align*}
\Gamma\left(  ak+1\right)   &  \approx\frac{a\left(  \sqrt{2\pi}\right)
^{1-a}}{\sqrt{ak+1}}\left(  a^{k}\sqrt{k+\frac{1}{a}}\Gamma\left(  k+\frac
{1}{a}\right)  \right)  ^{a}\\
&  <\frac{\sqrt{a}\left(  \sqrt{2\pi}\right)  ^{1-a}}{\sqrt{k}}\left(
a^{k}\Gamma\left(  k+\frac{1}{a}+1\right)  \right)  ^{a}<\sqrt{a}\left(
\sqrt{2\pi}\right)  ^{1-a}\left(  a^{k}\Gamma\left(  k+\frac{1}{a}+1\right)
\right)  ^{a}%
\end{align*}
For large $a$, we approximate as%
\begin{equation}
\Gamma\left(  ak+1\right)  \approx\frac{\sqrt{a}\left(  \sqrt{2\pi}\right)
^{1-a}}{\sqrt{k}}\left(  a^{k}\Gamma\left(  k+1\right)  \right)  ^{a}
\label{Gamma(ak+1)_as_power_a_of_Gamma(k+1)_approx}%
\end{equation}
In \textbf{art}. \ref{art_multiplication_formula_Gamma}, Gauss's
multiplication formula is written as $\Gamma\left(  nz+1\right)  =\left(
2\pi\right)  ^{-\frac{n-1}{2}}n^{nz+\frac{1}{2}}\prod_{k=1}^{n}\Gamma\left(
z+\frac{k}{n}\right)  $ and indicates (assuming $n=a$ is an integer) that%
\[
\Gamma\left(  ak+1\right)  =\sqrt{a}\left(  \sqrt{2\pi}\right)  ^{1-a}%
a^{ak}\prod_{j=1}^{a}\Gamma\left(  k+\frac{j}{a}\right)
\]
from which, for $k\geq2$,%
\[
\frac{\sqrt{a}\left(  \sqrt{2\pi}\right)  ^{1-a}}{k^{a}}\left(  a^{k}%
\Gamma\left(  k+1\right)  \right)  ^{a}<\Gamma\left(  ak+1\right)  <\sqrt
{a}\left(  \sqrt{2\pi}\right)  ^{1-a}\left(  a^{k}\Gamma\left(  k+1\right)
\right)  ^{a}%
\]

After introducing (\ref{Gamma(ak+1)_as_power_a_of_Gamma(k+1)_approx}) in the
definition (\ref{Mittag_Leffler_function_Eab}) of the Mittag-Leffler function
for real, positive $x$, we approximately obtain%
\[
E_{a}\left(  x\right)  =\sum_{k=0}^{\infty}\frac{x^{k}}{\Gamma\left(
ak+1\right)  }\approx\frac{\left(  \sqrt{2\pi}\right)  ^{a-1}}{\sqrt{a}}%
\sum_{k=0}^{\infty}\frac{\sqrt{k}\left(  \frac{x}{a}\right)  ^{k}}{\left(
k!\right)  ^{a}}%
\]
but Gauss's multiplication formula produces a lower bound%
\begin{equation}
E_{a}\left(  x\right)  =\sum_{k=0}^{\infty}\frac{x^{k}}{\Gamma\left(
ak+1\right)  }>\frac{\left(  \sqrt{2\pi}\right)  ^{a-1}}{\sqrt{a}}\sum
_{k=0}^{\infty}\frac{\left(  \frac{x}{a}\right)  ^{k}}{\left(  k!\right)
^{a}} \label{lower_bound_E_a(x)}%
\end{equation}

About 10 years before Mittag-Leffler has introduced his function $E_{a}\left(
x\right)  $, Hadamard \cite[p. 180]{Hadamard_Taylor_series} suggests in his
study of entire functions that $E_{a}\left(  x\right)  =\sum_{k=0}^{\infty
}\frac{x^{k}}{\Gamma\left(  ak+1\right)  }\sim\frac{\left(  \sqrt{2\pi
}\right)  ^{a-1}}{\sqrt{a}}\sum_{k=0}^{\infty}\frac{\left(  \frac{x}%
{a}\right)  ^{k}}{\left(  k!\right)  ^{a}}$ for large $x$. Hadamard derives an
exact $a$-fold integral for the last series, from which he deduces that
$\sum_{k=0}^{\infty}\frac{x^{k}}{\left(  k!\right)  ^{a}}<e^{ax^{\frac{1}{a}}%
}$. \textbf{Art}. \ref{art_fractional_a} shows that $E_{a}\left(  x\right)
<\frac{1}{a}e^{x^{\frac{1}{a}}}$. Combined with (\ref{lower_bound_E_a(x)})
leads to%
\[
\sum_{k=0}^{\infty}\frac{x^{k}}{\left(  k!\right)  ^{a}}<\frac{\left(
\sqrt{2\pi}\right)  ^{1-a}}{\sqrt{a}}e^{a^{\frac{1}{a}}x^{\frac{1}{a}}}%
\]
which is considerably sharper for $a>1$ than Hadamard's \cite[p.
180]{Hadamard_Taylor_series} bound.

We also give Hadamard's \cite[p. 180]{Hadamard_Taylor_series} nice argument,
starting from%
\[
e^{x^{\frac{1}{a}}}=\sum_{k=0}^{\infty}\frac{x^{\frac{k}{a}}}{\Gamma\left(
k+1\right)  }%
\]
and letting $m^{\prime}=\frac{k}{a}$, which runs over fractions for $a>1$, so
that $e^{x^{\frac{1}{a}}}=\sum_{m^{\prime}=0}^{\infty}\frac{x^{m^{\prime}}%
}{\Gamma\left(  m^{\prime}a+1\right)  }$. Comparing terms with $E_{a}\left(
x\right)  =\sum_{m=0}^{\infty}\frac{x^{m}}{\Gamma\left(  am+1\right)  }$, then
shows that $E_{a}\left(  x\right)  <e^{x^{\frac{1}{a}}}$. For $a<1$, Hadamard
states that $E_{a}\left(  x\right)  <\left[  \frac{1}{a}\right]  x^{\frac
{1}{a}}e^{x^{\frac{1}{a}}}$. However, the bounds presented here based on the
theory of the Mittag-Leffler function are sharper than Hadamard's estimates.

\medskip\refstepcounter{article}{\noindent\textbf{\thearticle. }%
}\ignorespaces\label{art_logarithmic_derivative} \emph{Logarithmic
derivative.} The logarithmic derivative $\frac{d}{dz}\log E_{a,b}\left(
z\right)  =\frac{\frac{d}{dz}E_{a,b}\left(  z\right)  }{E_{a,b}\left(
z\right)  }$ follows directly from (\ref{D_E_a,b(z)}) as%
\begin{equation}
\frac{d\log E_{a,b}\left(  z\right)  }{dz}=\frac{1}{az}\left(  \frac
{E_{a,b-1}\left(  z\right)  }{E_{a,b}\left(  z\right)  }-\left(  b-1\right)
\right)  \label{logarithmic_derivative_E_a,b(y)}%
\end{equation}

Since $b-1+ak>b-1$ for $k\geq1$ because $a>0$, we have for positive real $z$
and $b>1$,%
\[
E_{a,b}\left(  z\right)  =\sum_{k=0}^{\infty}\frac{z^{k}}{\left(
b-1+ak\right)  \Gamma\left(  b-1+ak\right)  }<\frac{1}{b-1}\sum_{k=0}^{\infty
}\frac{z^{k}}{\Gamma\left(  b-1+ak\right)  }=\frac{E_{a,b-1}\left(  z\right)
}{b-1}%
\]
and that $\left(  b-1\right)  <\frac{E_{a,b-1}\left(  z\right)  }%
{E_{a,b}\left(  z\right)  }$. The logarithmic derivative
(\ref{logarithmic_derivative_E_a,b(y)}) becomes $\frac{d\log E_{a,b}\left(
z\right)  }{dz}>0$, illustrating that $\log E_{a,b}\left(  z\right)  $ is
increasing for real $z\geq0$ and $b>1$. A little more precisely with
$b-1+ak>b-1+a$ for $k\geq1$, we have for real positive $z$ and for $b>1-a$
(because then $\Gamma\left(  b-1+ak\right)  >0$ for $k\geq1$),%
\begin{align*}
E_{a,b}\left(  z\right)   &  =\frac{1}{\Gamma\left(  b\right)  }+\sum
_{k=1}^{\infty}\frac{z^{k}}{\left(  b-1+ak\right)  \Gamma\left(
b-1+ak\right)  }\\
&  <\frac{1}{\Gamma\left(  b\right)  }+\frac{1}{b-1+a}\left(  \sum
_{k=0}^{\infty}\frac{z^{k}}{\Gamma\left(  b-1+ak\right)  }-\frac{1}%
{\Gamma\left(  b-1\right)  }\right)
\end{align*}
and%
\[
E_{a,b}\left(  z\right)  <\frac{1}{\left(  b-1+a\right)  }\left(  \frac
{a}{\Gamma\left(  b\right)  }+E_{a,b-1}\left(  z\right)  \right)
\]
Thus,%
\[
a\left(  1-\frac{1}{\Gamma\left(  b\right)  E_{a,b}\left(  z\right)  }\right)
<\frac{E_{a,b-1}\left(  z\right)  }{E_{a,b}\left(  z\right)  }-\left(
b-1\right)
\]
and (\ref{logarithmic_derivative_E_a,b(y)}) becomes%
\[
\frac{d\log E_{a,b}\left(  z\right)  }{dz}>\frac{1}{z}\left(  1-\frac
{1}{\Gamma\left(  b\right)  E_{a,b}\left(  z\right)  }\right)
\]

\textbf{Example} If $b=1$ and $a=1$, then $E_{1,1}\left(  z\right)  =e^{z}$
and the above inequality results in the well-known bound (see e.g. \cite[p.
103]{PVM_PAComplexNetsCUP}) that $e^{-z}>1-z$ for all real $z$.

\medskip\refstepcounter{article}{\noindent\textbf{\thearticle. }%
}\ignorespaces\label{art_Taylor_expansion_rond_z0} \emph{Taylor expansion
around }$z_{0}$\emph{.} Since the Mittag-Leffler function $E_{a,b}\left(
z\right)  $ is an entire function, any Taylor expansion around an arbitrary
(finite) point $z_{0}$ has infinite radius of convergence,%
\[
E_{a,b}\left(  z\right)  =\sum_{m=0}^{\infty}\frac{1}{m!}\left.  \frac
{d^{m}E_{a,b}\left(  z\right)  }{dx^{m}}\right\vert _{x=z_{0}}\left(
z-z_{0}\right)  ^{m}%
\]
With the $m$-th derivative (\ref{m_th_derivative_E_a,b(z)}), we find%
\[
E_{a,b}\left(  z\right)  =\sum_{m=0}^{\infty}\left\{  \sum_{j=0}^{m}%
q_{j}\left(  a,b,m\right)  \left(  E_{a,b-j}\left(  z_{0}\right)  -\sum
_{l=0}^{m-1}\frac{z_{0}^{l}}{\Gamma\left(  b-j+la\right)  }\right)  \right\}
\frac{\left(  \frac{z-z_{0}}{az_{0}}\right)  ^{m}}{m!}%
\]
but, unfortunately (see \textbf{art}. \ref{art_differentiation}), the
coefficients $q_{j}\left(  a,b,m\right)  $ are not available in closed form.
However, the closed form (\ref{Dm_(z^(b-1)E_a,b(z^a))}), which is a rather
fundamental property of $E_{a,b}\left(  z\right)  $, opens a new avenue. The
Taylor expansion around $z_{0}$ of%
\[
z^{b-1}E_{a,b}\left(  z^{a}\right)  =\sum_{m=0}^{\infty}\frac{1}{m!}\left.
\frac{d^{m}}{dx^{m}}\left\{  x^{b-1}E_{a,b}\left(  x^{a}\right)  \right\}
\right\vert _{x=z_{0}}\left(  z-z_{0}\right)  ^{m}%
\]
becomes with the differentiation recursion (\ref{Dm_(z^(b-1)E_a,b(z^a))})%
\[
z^{b-1}E_{a,b}\left(  z^{a}\right)  =z_{0}^{b-1}\sum_{m=0}^{\infty}%
\frac{E_{a,b-m}\left(  z_{0}^{a}\right)  }{m!}z_{0}^{-m}\left(  z-z_{0}%
\right)  ^{m}%
\]
The function $z^{b-1}E_{a,b}\left(  z^{a}\right)  $ has a branch cut at the
negative real axis for real $a$ and $b$, implying that the radius $R$ of
convergence equals $\left\vert z_{0}\right\vert $. Using its definition in
Section \ref{sec_complex_function_theory},%
\[
\frac{1}{R}=\lim_{m\rightarrow\infty}\left\vert \frac{z_{0}^{-m-1}%
E_{a,b-m-1}\left(  z_{0}^{a}\right)  m!}{z_{0}^{-m}E_{a,b-m}\left(  z_{0}%
^{a}\right)  \left(  m+1\right)  !}\right\vert =z_{0}^{-1}\lim_{m\rightarrow
\infty}\left\vert \frac{E_{a,b-m-1}\left(  z_{0}^{a}\right)  }{E_{a,b-m}%
\left(  z_{0}^{a}\right)  \left(  m+1\right)  }\right\vert
\]
then indicates for any finite $z\neq0$, $a$ and $b$ that
\begin{equation}
\left\vert \frac{E_{a,b-m-1}\left(  z\right)  }{E_{a,b-m}\left(  z\right)
}\right\vert \sim\left(  m+1\right)  \text{ if }m\rightarrow\infty
\label{asymptotic_E_a,b-m(z^a)_large_m}%
\end{equation}
Introducing (\ref{logarithmic_derivative_E_a,b(y)}) for $z\neq0$ that%
\[
1=\lim_{m\rightarrow\infty}\left\vert \frac{E_{a,b-m-1}\left(  z\right)
}{E_{a,b-m}\left(  z\right)  \left(  m+1\right)  }\right\vert =\lim
_{m\rightarrow\infty}\left\vert \frac{az}{m}\frac{d\log E_{a,b-m}\left(
z\right)  }{dz}-\left(  1-\frac{b}{m}\right)  \right\vert
\]
shows that
\[
\lim_{m\rightarrow\infty}\left\vert \frac{1}{m}\frac{d\log E_{a,b-m}\left(
z\right)  }{dz}\right\vert =0
\]
Thus, the logarithmic derivative for $z\neq0$ and finite $a$ and $b$ is of
order $\frac{d\log E_{a,b-m}\left(  z\right)  }{dz}=O\left(  m^{1-\varepsilon
}\right)  $ for any positive $\varepsilon>0$.

We proceed by removing real powers of $z^{\beta}=e^{\beta\log z}$ that destroy
analyticity in the complex plane and introduce formally the Taylor series
$\left(  \frac{z}{z_{0}}\right)  ^{1-b}=\left(  \frac{z}{z_{0}}-1+1\right)
^{1-b}=\sum_{k=0}^{\infty}\binom{1-b}{k}\left(  \frac{z}{z_{0}}-1\right)
^{k}$, valid for any $b$ and $\left\vert z\right\vert <\left\vert
z_{0}\right\vert $, in
\[
E_{a,b}\left(  z^{a}\right)  =\left(  \frac{z}{z_{0}}\right)  ^{1-b}\sum
_{m=0}^{\infty}\frac{E_{a,b-m}\left(  z_{0}^{a}\right)  }{m!}\left(  \frac
{z}{z_{0}}-1\right)  ^{m}%
\]
yielding, after executing the Cauchy product,%
\[
E_{a,b}\left(  z^{a}\right)  =\sum_{m=0}^{\infty}\left\{  \sum_{k=0}^{m}%
\binom{1-b}{m-k}\frac{E_{a,b-k}\left(  z_{0}^{a}\right)  }{k!}\right\}
\left(  \frac{z}{z_{0}}-1\right)  ^{m}\hspace{1cm}\text{for }\left\vert
z\right\vert <\left\vert z_{0}\right\vert
\]
Next, letting $y=z^{a}$ and $y_{0}=z_{0}^{a}$, we first expand in a Taylor
series around $y_{0}$%
\[
\left(  \left(  \frac{y}{y_{0}}\right)  ^{\frac{1}{a}}-1\right)  ^{m}%
=\sum_{k=0}^{\infty}\frac{1}{k!}\left.  \frac{d^{k}}{dz^{k}}\left(  \left(
\frac{y}{y_{0}}\right)  ^{\frac{1}{a}}-1\right)  ^{m}\right\vert _{z=y_{0}%
}\left(  y-y_{0}\right)  ^{k}%
\]
as
\begin{align*}
\left(  \left(  \frac{y}{y_{0}}\right)  ^{\frac{1}{a}}-1\right)  ^{m}  &
=\sum_{q=0}^{m}\binom{m}{q}\left(  \frac{y}{y_{0}}\right)  ^{\frac{q}{a}%
}\left(  -1\right)  ^{m-q}=\sum_{q=0}^{m}\binom{m}{q}\left(  \frac{y-y_{0}%
}{y_{0}}+1\right)  ^{\frac{q}{a}}\left(  -1\right)  ^{m-q}\\
&  =\sum_{j=0}^{\infty}\left\{  \sum_{q=0}^{m}\binom{m}{q}\binom{\frac{q}{a}%
}{j}\left(  -1\right)  ^{m-q}\right\}  \left(  \frac{y-y_{0}}{y_{0}}\right)
^{j}%
\end{align*}
Substitution and reversing the $m$- and $k$-sum yields%

\begin{align*}
E_{a,b}\left(  y\right)   &  =\sum_{j=0}^{\infty}\left[  \sum_{m=0}^{\infty
}\sum_{k=0}^{m}\binom{1-b}{m-k}\frac{E_{a,b-k}\left(  y_{0}\right)  }%
{k!}\left\{  \sum_{q=0}^{m}\binom{m}{q}\binom{\frac{q}{a}}{j}\left(
-1\right)  ^{m-q}\right\}  \right]  \left(  \frac{y-y_{0}}{y_{0}}\right)
^{j}\\
&  =\sum_{j=0}^{\infty}\left[  \sum_{k=0}^{\infty}\frac{E_{a,b-k}\left(
y_{0}\right)  }{k!}\sum_{m=k}^{\infty}\binom{1-b}{m-k}\sum_{q=0}^{m}\binom
{m}{q}\binom{\frac{q}{a}}{j}\left(  -1\right)  ^{m-q}\right]  \left(
\frac{y-y_{0}}{y_{0}}\right)  ^{j}%
\end{align*}
The characteristic coefficients \cite[Appendix]{PVM_ASYM} of a complex
function $f\left(  z\right)  $ with Taylor series $f\left(  z\right)
=\sum_{k=0}^{\infty}f_{k}\left(  z_{0}\right)  \left(  z-z_{0}\right)  ^{k}$,
defined by $\left.  s[k,m]\right\vert _{f}\left(  z_{0}\right)  =\frac{1}%
{m!}\left.  \frac{d^{m}}{dz^{m}}\left(  f\left(  z\right)  -f\left(
z_{0}\right)  ^{k}\right)  \right\vert _{z=z_{0}}$, possesses a general form
\begin{equation}
\left.  s[k,m]\right\vert _{f}\left(  z_{0}\right)  =\sum_{\sum_{i=1}^{k}%
j_{i}=m;j_{i}>0}\prod_{i=1}^{k}f_{j_{i}}\left(  z_{0}\right)  \label{def_s}%
\end{equation}
and obeys $\left.  s[k,m]\right\vert _{f}\left(  z_{0}\right)  =0$ if $k<0$
and $k>m$. Also, $\left.  s[1,m]\right\vert _{f}\left(  z_{0}\right)
=f_{m}\left(  z_{0}\right)  $, while (\ref{def_s}) indicates that $\left.
s[m,m]\right\vert _{f}\left(  z_{0}\right)  =f_{1}^{m}\left(  z_{0}\right)  $.
We can show \cite{PVM_charcoef}\ that%
\begin{equation}
\left.  s[k,m]\right\vert _{(1+z)^{\alpha}}=\sum_{j=1}^{k}(-1)^{j+k}{\binom
{k}{j}}{\binom{\alpha j}{m}=}\frac{k!}{m!}\sum_{j=k}^{m}S_{m}^{(j)}%
\mathcal{S}_{j}^{(k)}\,\alpha^{j} \label{chc_(1+z)^alpha_2}%
\end{equation}
where $S_{m}^{(k)}$ and $\mathcal{S}_{m}^{(k)}$ are the Stirling numbers of
the first and second kind \cite[Sec. 24.1.3 and 24.1.4]{Abramowitz},
respectively. We apply these properties to the Taylor series%
\begin{align*}
E_{a,b}\left(  y\right)   &  =\sum_{j=0}^{\infty}\left[  \sum_{k=0}^{\infty
}\frac{E_{a,b-k}\left(  y_{0}\right)  }{k!}\sum_{m=k}^{\infty}\binom{1-b}%
{m-k}\left.  s[m,j]\right\vert _{(1+z)^{\frac{1}{a}}}\right]  \left(
\frac{y-y_{0}}{y_{0}}\right)  ^{j}\\
&  =\sum_{j=0}^{\infty}\left[  \sum_{k=0}^{j}\frac{E_{a,b-k}\left(
y_{0}\right)  }{k!}\sum_{m=k}^{j}\binom{1-b}{m-k}\left.  s[m,j]\right\vert
_{(1+z)^{\frac{1}{a}}}\right]  \left(  \frac{y-y_{0}}{y_{0}}\right)  ^{j}%
\end{align*}
Finally, we arrive with (\ref{chc_(1+z)^alpha_2}) at the Taylor series of
$E_{a,b}\left(  y\right)  $ around $y_{0}$,
\begin{equation}
E_{a,b}\left(  y\right)  =\sum_{j=0}^{\infty}\left[  \frac{y_{0}^{-j}}{j!}%
\sum_{k=0}^{j}E_{a,b-k}\left(  y_{0}\right)  \sum_{m=k}^{j}\binom{m}{k}%
\frac{\Gamma\left(  2-b\right)  }{\Gamma\left(  2-b-m+k\right)  }\sum
_{q=m}^{j}S_{j}^{(q)}\mathcal{S}_{q}^{(m)}\,\frac{1}{a^{q}}\right]  \left(
y-y_{0}\right)  ^{j} \label{Taylor_series_E_a,b(y)_around_y0}%
\end{equation}
from which the closed form of the $j$-th derivative of the Mittag-Leffer
function, evaluated at $y_{0}=w$, follows as%
\begin{equation}
\left.  \frac{d^{j}E_{a,b}\left(  z\right)  }{dz^{j}}\right\vert _{z=w}%
=w^{-j}\sum_{k=0}^{j}E_{a,b-k}\left(  w\right)  \sum_{m=k}^{j}\binom{m}%
{k}\frac{\Gamma\left(  2-b\right)  }{\Gamma\left(  2-b-m+k\right)  }\sum
_{q=m}^{j}S_{j}^{(q)}\mathcal{S}_{q}^{(m)}\,\frac{1}{a^{q}}
\label{m_th_derivative_E_a,b(z)_closed_form}%
\end{equation}
Since the Mittag-Leffer function $E_{a,b}\left(  z\right)  $ is an entire
function, the $j$-th Taylor coefficient around $w$ decreases faster than any
power of $\left(  z-w\right)  ^{j}$ for large $j$ and any $z$, we deduce that
$\left.  \frac{d^{j}E_{a,b}\left(  z\right)  }{dz^{j}}\right\vert
_{z=w}=o\left(  j!w^{j}\right)  $. Thus, the $j$-th Taylor coefficient
$\frac{1}{j!}\left.  \frac{d^{j}E_{a,b}\left(  z\right)  }{dz^{j}}\right\vert
_{z=w}$ increases at most as a polynomial in $w$ of order $j$. In contrast to
the relatively simple Taylor series of $E_{a,b}\left(  z\right)  $ in
(\ref{Mittag_Leffler_function_Eab}) around the origin, the general form
(\ref{m_th_derivative_E_a,b(z)_closed_form}) emphasizes the complicated nature
of the Mittag-Leffer function $E_{a,b}\left(  z\right)  $ elsewhere in the
complex plane.

As a check of (\ref{m_th_derivative_E_a,b(z)_closed_form}), for $a=b=1$, we
have with the orthogonality condition of the Stirling numbers%
\begin{equation}
\sum_{k=n}^{m}S_{k}^{(n)}\mathcal{S}_{m}^{(k)}=\delta_{nm}
\label{orthogonality_StirlingNumbers_1}%
\end{equation}
that%
\[
E_{a,b}\left(  y\right)  =\sum_{j=0}^{\infty}\left[  \sum_{k=0}^{j}%
E_{1,1-k}\left(  y_{0}\right)  \sum_{m=k}^{j}\binom{m}{k}\frac{\Gamma\left(
2-b\right)  }{\Gamma\left(  2-b-m+k\right)  }\delta_{jm}\,\right]
\frac{\left(  \frac{y-y_{0}}{y_{0}}\right)  ^{j}}{j!}%
\]
Applying (\ref{E_a,minus(ma)}) yields $x^{-m}E_{1,1-m}\left(  x\right)
=E_{1,1}\left(  x\right)  +\sum_{k=1}^{m}\frac{\left(  \frac{1}{x}\right)
^{k}}{\Gamma\left(  1-k\right)  }$ and $E_{1,1-m}\left(  x\right)  =x^{m}%
e^{x}$ so that%
\begin{align*}
E_{a,b}\left(  y\right)   &  =e^{y_{0}}\sum_{j=0}^{\infty}\left[  \sum
_{k=0}^{j}y_{0}^{k-j}\binom{j}{k}\frac{1}{\Gamma\left(  1-j+k\right)
}\right]  \frac{\left(  y-y_{0}\right)  ^{j}}{j!}\\
&  =e^{y_{0}}\sum_{j=0}^{\infty}\frac{\left(  y-y_{0}\right)  ^{j}}%
{j!}=e^{y_{0}}e^{y-y_{0}}=e^{y}%
\end{align*}

\medskip\refstepcounter{article}{\noindent\textbf{\thearticle. }%
}\ignorespaces\label{art_second_order_logarithmic_derivative} \emph{Second
order logarithmic derivative.} Similarly as in \textbf{art}.
\ref{art_logarithmic_derivative}, we can directly differentiate
(\ref{logarithmic_derivative_E_a,b(y)}) again,
\[
\frac{d^{2}\log E_{a,b}\left(  z\right)  }{dz^{2}}=-\frac{1}{az^{2}}\left(
\frac{E_{a,b-1}\left(  z\right)  }{E_{a,b}\left(  z\right)  }-\left(
b-1\right)  \right)  +\frac{1}{az}\frac{d}{dz}\frac{E_{a,b-1}\left(  z\right)
}{E_{a,b}\left(  z\right)  }%
\]
using (\ref{D_E_a,b(z)}), resulting in
\[
\frac{d}{dz}\frac{E_{a,b-1}\left(  z\right)  }{E_{a,b}\left(  z\right)
}=\frac{1}{az}\left\{  \frac{E_{a,b-2}\left(  z\right)  }{E_{a,b}\left(
z\right)  }+\frac{E_{a,b-1}\left(  z\right)  }{E_{a,b}\left(  z\right)
}-\left(  \frac{E_{a,b-1}\left(  z\right)  }{E_{a,b}\left(  z\right)
}\right)  ^{2}\right\}
\]
After tedious manipulations, we arrive at the second order logarithmic
derivative%
\begin{equation}
\frac{d^{2}\log E_{a,b}\left(  z\right)  }{dz^{2}}=\frac{1}{\left(  az\right)
^{2}}\left\{  \frac{E_{a,b-2}\left(  z\right)  }{E_{a,b}\left(  z\right)
}+a\left(  b-1\right)  -\left(  a-1\right)  \frac{E_{a,b-1}\left(  z\right)
}{E_{a,b}\left(  z\right)  }-\left(  \frac{E_{a,b-1}\left(  z\right)
}{E_{a,b}\left(  z\right)  }\right)  ^{2}\right\}
\label{logarithmic_derivative_E_a,b(y)_2nd_order}%
\end{equation}

\medskip\refstepcounter{article}{\noindent\textbf{\thearticle. }%
}\ignorespaces\label{art_Taylor_series_logE_a,b(z)} \emph{The Taylor series of
}$\log E_{a,b}\left(  z\right)  $\emph{ around }$z=0$. From the power series
definition (\ref{Mittag_Leffler_function_Eab}) of $E_{a,b}\left(  z\right)  $,
we have%
\begin{align*}
\lim_{z\rightarrow0}\frac{d\log E_{a,b}\left(  z\right)  }{dz}  &
=\lim_{z\rightarrow0}\frac{1}{\sum_{k=0}^{\infty}\frac{z^{k}}{\Gamma\left(
b+ak\right)  }}\lim_{z\rightarrow0}\frac{1}{az}\left(  \sum_{k=0}^{\infty
}\frac{z^{k}}{\Gamma\left(  b-1+ak\right)  }-\sum_{k=0}^{\infty}\frac{\left(
b-1\right)  z^{k}}{\Gamma\left(  b+ak\right)  }\right) \\
&  =\Gamma\left(  b\right)  \lim_{z\rightarrow0}\sum_{k=1}^{\infty}%
\frac{kz^{k-1}}{\Gamma\left(  b+ak\right)  }%
\end{align*}
and find%
\[
\lim_{z\rightarrow0}\frac{d\log E_{a,b}\left(  z\right)  }{dz}=\frac
{\Gamma\left(  b\right)  }{\Gamma\left(  b+a\right)  }%
\]

Proceeding in this way to higher-order derivatives becomes cumbersome. The
general theory of characteristic coefficients \cite[Appendix]{PVM_ASYM}
provides us with%
\[
\log f(z)=\log f_{0}\left(  z_{0}\right)  +\sum_{m=1}^{\infty}\left(
\sum_{k=1}^{m}\frac{(-1)^{k-1}}{k\;f_{0}^{k}\left(  z_{0}\right)
}\,s[k,m]\left(  z_{0}\right)  \right)  \;\left(  z-z_{0}\right)  ^{m}%
\]
Confining to $z_{0}=0$, the characteristic coefficients (\ref{def_s}) of the
Mittag-Leffler function $E_{a,b}\left(  z\right)  $ are%
\[
\left.  s[k,m]\right\vert _{E_{a,b}\left(  z\right)  }=\sum_{\sum_{i=1}%
^{k}j_{i}=m;j_{i}>0}\prod_{i=1}^{k}\frac{1}{\Gamma\left(  b+j_{i}a\right)  }%
\]
and the Taylor series of $\log E_{a,b}\left(  z\right)  $ around $z_{0}=0$
follows as
\[
\log E_{a,b}\left(  z\right)  =\sum_{m=0}^{\infty}c_{m}\;z^{m}%
\]
where we define $c_{0}=\log\left(  -\Gamma\left(  b\right)  \right)  $ and the
coefficients $c_{m}$ for $m>0$
\begin{align*}
c_{m}  &  =-\sum_{k=1}^{m}\frac{\left(  -\Gamma\left(  b\right)  \right)
^{k}}{k}\,\left.  s[k,m]\right\vert _{E_{a,b}\left(  z\right)  }\\
&  =\frac{\Gamma\left(  b\right)  }{\Gamma\left(  b+ma\right)  }\,+\sum
_{k=2}^{m-1}\frac{\left(  -1\right)  ^{k-1}}{k}\sum_{\sum_{i=1}^{k}%
j_{i}=m;j_{i}>0}\prod_{i=1}^{k}\frac{\Gamma\left(  b\right)  }{\Gamma\left(
b+j_{i}a\right)  }+\frac{\left(  -1\right)  ^{m-1}}{m}\left(  \frac
{\Gamma\left(  b\right)  }{\Gamma\left(  b+a\right)  }\right)  ^{m}%
\end{align*}
The list of the coefficients $c_{m}$ for $m=1$ up to $m=5$ is%
\begin{align*}
c_{1}  &  =\frac{\Gamma\left(  b\right)  }{\Gamma\left(  b+a\right)  }\\
c_{2}  &  =\frac{\Gamma\left(  b\right)  }{\Gamma\left(  b+2a\right)  }%
-\frac{1}{2}\left(  \frac{\Gamma\left(  b\right)  }{\Gamma\left(  b+a\right)
}\right)  ^{2}\\
c_{3}  &  =\frac{\Gamma\left(  b\right)  }{\Gamma\left(  b+3a\right)  }%
+\frac{\Gamma^{2}\left(  b\right)  }{\Gamma\left(  b+2a\right)  \Gamma\left(
b+a\right)  }+\frac{1}{3}\left(  \frac{\Gamma\left(  b\right)  }{\Gamma\left(
b+a\right)  }\right)  ^{3}\\
c_{4}  &  =\frac{\Gamma\left(  b\right)  }{\Gamma\left(  b+4a\right)  }%
-\frac{\Gamma^{2}\left(  b\right)  }{\Gamma\left(  b+3a\right)  \Gamma\left(
b+a\right)  }-\frac{1}{2}\left(  \frac{\Gamma\left(  b\right)  }{\Gamma\left(
b+2a\right)  }\right)  ^{2}+\frac{\Gamma^{3}\left(  b\right)  }{\Gamma\left(
b+2a\right)  \Gamma^{2}\left(  b+a\right)  }-\frac{1}{4}\left(  \frac
{\Gamma\left(  b\right)  }{\Gamma\left(  b+a\right)  }\right)  ^{4}\\
c_{5}  &  =\frac{\Gamma\left(  b\right)  }{\Gamma\left(  b+5a\right)  }%
-\frac{\Gamma^{2}\left(  b\right)  }{\Gamma\left(  b+4a\right)  \Gamma\left(
b+a\right)  }-\frac{\Gamma^{2}\left(  b\right)  }{\Gamma\left(  b+3a\right)
\Gamma\left(  b+2a\right)  }+\frac{\Gamma^{3}\left(  b\right)  }{\Gamma\left(
b+3a\right)  \Gamma^{2}\left(  b+a\right)  }+\frac{\Gamma^{3}\left(  b\right)
}{\Gamma\left(  b+a\right)  \Gamma^{2}\left(  b+2a\right)  }\\
&  -\frac{\Gamma^{4}\left(  b\right)  }{\Gamma^{3}\left(  b+a\right)
\Gamma\left(  b+2a\right)  }+\frac{1}{5}\left(  \frac{\Gamma\left(  b\right)
}{\Gamma\left(  b+a\right)  }\right)  ^{5}%
\end{align*}
Unfortunately, it seems difficult to further sum the terms in $c_{m}$. The
closest zero to $z_{0}=0$ lies at a distance from $z_{0}$ equal to the radius
$R$ of convergence of $\log E_{a,b}\left(  z\right)  $, which is%
\[
\frac{1}{R}=\lim_{m\rightarrow\infty}\left\vert \frac{c_{m+1}}{c_{m}%
}\right\vert =\lim_{m\rightarrow\infty}\frac{\left\vert \sum_{k=1}^{m+1}%
\frac{\left(  -\Gamma\left(  b\right)  \right)  ^{k}}{k}\,\left.
s[k,m+1]\right\vert _{E_{a,b}\left(  z\right)  }\right\vert }{\left\vert
\sum_{k=1}^{m}\frac{\left(  -\Gamma\left(  b\right)  \right)  ^{k}}%
{k}\,\left.  s[k,m]\right\vert _{E_{a,b}\left(  z\right)  }\right\vert }%
\]
The generalization towards the closest zero to $z_{0}$ requires the Taylor
coefficients of (\ref{Taylor_series_E_a,b(y)_around_y0}).

\medskip\refstepcounter{article}{\noindent\textbf{\thearticle. }%
}\ignorespaces\label{art_derivation_vs_a_and_b} \emph{Derivation of }%
$E_{a,b}\left(  z\right)  $\emph{ with respect to the parameters }$a$\emph{
and }$b$. From the definition (\ref{Mittag_Leffler_function_Eab}), partial
differentiating yields%
\begin{align*}
\frac{\partial}{\partial a}E_{a,b}\left(  z\right)   &  =\frac{\partial
}{\partial a}\left(  \frac{1}{\Gamma\left(  b\right)  }+\sum_{k=1}^{\infty
}\frac{z^{k}}{\Gamma\left(  b+ak\right)  }\right) \\
&  =\sum_{k=1}^{\infty}\left.  \frac{d}{dy}\frac{1}{\Gamma\left(  y\right)
}\right\vert _{y=b+ak}\frac{dy}{da}z^{k}=-\sum_{k=1}^{\infty}\left.
\frac{\psi\left(  y\right)  }{\Gamma\left(  y\right)  }\right\vert
_{y=b+ak}kz^{k}%
\end{align*}
while, similarly but containing index $k=0$,%
\[
\frac{\partial}{\partial b}E_{a,b}\left(  z\right)  =-\sum_{k=0}^{\infty
}\left.  \frac{\psi\left(  y\right)  }{\Gamma\left(  y\right)  }\right\vert
_{y=b+ak}z^{k}%
\]
Hence, we observe that%
\[
\frac{\partial}{\partial a}E_{a,b}\left(  z\right)  =z\frac{\partial^{2}%
}{\partial z\partial b}E_{a,b}\left(  z\right)
\]

Partial differentiating $m$-times gives%
\[
\frac{\partial^{m}}{\partial a^{m}}E_{a,b}\left(  z\right)  =\sum
_{k=0}^{\infty}\left.  \frac{d^{m}}{dy^{m}}\frac{1}{\Gamma\left(  y\right)
}\right\vert _{y=b+ak}k^{m}z^{k}%
\]
which suggest to let $z=ye^{w}$ so that%
\[
\frac{\partial^{m}}{\partial a^{m}}E_{a,b}\left(  ye^{w}\right)  =\sum
_{k=0}^{\infty}\left.  \frac{d^{m}}{dy^{m}}\frac{1}{\Gamma\left(  y\right)
}\right\vert _{y=b+ak}k^{m}y^{k}e^{kw}=\frac{\partial^{m}}{\partial w^{m}}%
\sum_{k=0}^{\infty}\left.  \frac{d^{m}}{dy^{m}}\frac{1}{\Gamma\left(
y\right)  }\right\vert _{y=b+ak}y^{k}e^{kw}%
\]
leading to the partial differentiation equation for any integer $m\geq0$ and
any $y$ (independent of $a,b$ and $w$),%
\begin{equation}
\frac{\partial^{m}}{\partial a^{m}}E_{a,b}\left(  ye^{w}\right)
=\frac{\partial^{2m}}{\partial w^{m}\partial b^{m}}E_{a,b}\left(
ye^{w}\right)  \label{partial_dvgl_m_times_E_a,b(e^w)}%
\end{equation}

\medskip\refstepcounter{article}{\noindent\textbf{\thearticle. }%
}\ignorespaces\label{art_Taylor_series_to_Ea,b} \emph{A Taylor series approach
with Fermi-Dirac integrals}. We introduce the Taylor series of the entire
function $\frac{1}{\Gamma\left(  b+ak\right)  }$ around $q$ into the
definition (\ref{Mittag_Leffler_function_Eab}) of $E_{a,b}\left(  z\right)  $,%
\[
E_{a,b}\left(  z\right)  =\frac{1}{\Gamma\left(  b\right)  }+\sum
_{k=1}^{\infty}\frac{z^{k}}{\Gamma\left(  b+ak\right)  }=\frac{1}%
{\Gamma\left(  b\right)  }+\sum_{k=1}^{\infty}z^{k}\sum_{m=0}^{\infty}\frac
{1}{m!}\left.  \frac{d^{m}}{du^{m}}\frac{1}{\Gamma\left(  u\right)
}\right\vert _{u=q}\left(  ak+b-q\right)  ^{m}%
\]
and%
\[
E_{a,b}\left(  z\right)  =\frac{1}{\Gamma\left(  b\right)  }+\frac{1}%
{\Gamma\left(  q\right)  }\frac{z}{1-z}+\sum_{m=1}^{\infty}\frac{1}{m!}\left.
\frac{d^{m}}{du^{m}}\frac{1}{\Gamma\left(  u\right)  }\right\vert _{u=q}%
a^{m}\sum_{k=1}^{\infty}z^{k}\left(  k+\frac{b-q}{a}\right)  ^{m}%
\]
The reversal in the $k$- and $m$-sum leads to a confinement of $\left\vert
z\right\vert <1$. We will now choose $q=b$ and evaluate the series $\sum
_{k=1}^{\infty}k^{m}z^{k}$, that converges for $\left\vert z\right\vert <1$.

The Fermi-Dirac integral of order $p$ is defined as
\begin{equation}
F_{p}(y)=\frac{1}{\Gamma(p+1)}\int_{0}^{\infty}\frac{x^{p}}{1+e^{x-y}}\;dx
\label{def_Fp}%
\end{equation}
The value of the zero argument in $y$ is immediately written in terms of the
Eta function,
\begin{equation}
F_{p}(0)=\eta(p+1) \label{Fp0}%
\end{equation}
where the Eta function $\eta\left(  s\right)  $ is related to the Riemann Zeta
function $\zeta(s)$ as
\begin{equation}
\eta(s)=(1-2^{1-s})\zeta(s) \label{def_Eta_in_terms_Zeta}%
\end{equation}
By expanding $\frac{1}{1+e^{x-y}}=\frac{e^{-x+y}}{1+e^{-x+y}}=\sum
_{k=1}^{\infty}(-1)^{k-1}e^{-k\left(  x-y\right)  }$ for $\operatorname{Re}%
\left(  y\right)  <0$ in (\ref{def_Fp}), the Dirichlet series for all complex
$p$ is readily deduced as
\begin{equation}
F_{p}(y)=\sum_{k=1}^{\infty}(-1)^{k-1}\frac{\left(  e^{y}\right)  ^{k}%
}{k^{p+1}} \label{Fp_dirichlet}%
\end{equation}
In particular, $F_{-1}\left(  y\right)  =\frac{1}{1+e^{-y}}$. Hence, we can
write $\sum_{k=1}^{\infty}k^{m}\left(  -z\right)  ^{k}=-F_{-m-1}(\log z)$ for
$\left\vert z\right\vert <1$ and
\begin{equation}
E_{a,b}\left(  -z\right)  =\frac{1}{\Gamma\left(  b\right)  }-\sum
_{m=0}^{\infty}\frac{1}{m!}\left.  \frac{d^{m}}{du^{m}}\frac{1}{\Gamma\left(
u\right)  }\right\vert _{u=b}F_{-m-1}(\log z)a^{m}
\label{Ea,b(-z)_Taylor_in_a_around_b_in_plusa}%
\end{equation}
For integer negative order and $k>0$, it can be shown \cite{PVM_charcoef}
that
\[
F_{-k}(y)=\frac{d^{k-1}}{dy^{k-1}}\left(  \frac{1}{1+e^{-y}}\right)
=\sum_{m=1}^{k}(m-1)!(-1)^{m-1}\mathcal{S}_{k}^{(m)}\,\left(  \frac
{1}{1+e^{-y}}\right)  ^{m}%
\]
where $\mathcal{S}_{k}^{(m)}$ is the Stirling Number of the Second Kind
\cite[24.1.4]{Abramowitz}. Since $\frac{1}{1+e^{-y}}=1-\frac{1}{1+e^{y}}$,
which is equivalent to $F_{-1}\left(  y\right)  =1-F_{-1}\left(  -y\right)  $,
the $k$-th derivative shows that, for $k>1$,%
\begin{equation}
F_{-k}(y)=\left(  -1\right)  ^{k}F_{-k}(-y) \label{negFD_negargument}%
\end{equation}
and, thus extending the above for $\left\vert z\right\vert <1$ to,%
\begin{equation}
E_{a,b}\left(  -z\right)  =\sum_{m=0}^{\infty}\frac{1}{m!}\left.  \frac{d^{m}%
}{du^{m}}\frac{1}{\Gamma\left(  u\right)  }\right\vert _{u=b}F_{-m-1}(-\log
z)\left(  -a\right)  ^{m} \label{Ea,b(-z)_Taylor_in_a_around_b_in_mina}%
\end{equation}

Stretching the convergence constraint in
(\ref{Ea,b(-z)_Taylor_in_a_around_b_in_mina}) to $z=1$ and using (\ref{Fp0})
results in%
\[
E_{a,b}\left(  -1\right)  =\sum_{m=0}^{\infty}\frac{1}{m!}\left.  \frac{d^{m}%
}{dz^{m}}\frac{1}{\Gamma\left(  z\right)  }\right\vert _{z=b}\eta\left(
-m\right)  \left(  -a\right)  ^{m}%
\]
Further by (\ref{def_Eta_in_terms_Zeta}), it holds that $\eta(-m)=(1-2^{1+m}%
)\zeta(-m)=(1-2^{1+m})\frac{(-1)^{m}}{m+1}B_{m+1}$, because $\zeta
(-n)=\frac{(-1)^{n}}{n+1}B_{n+1}$ and $\zeta(-2n)=0$ for $n>0$. Taking into
account that the odd Bernoulli numbers $B_{2m+1}=0$ for $m>0$, we find%
\[
E_{a,b}\left(  -1\right)  =-\frac{1}{2\Gamma\left(  b\right)  }-\sum
_{k=1}^{\infty}\frac{B_{2k}}{\left(  2k\right)  !}\left.  \frac{d^{2k-1}%
}{dz^{2k-1}}\frac{1}{\Gamma\left(  z\right)  }\right\vert _{z=b}%
(2^{2k}-1)a^{2k-1}%
\]
which converges fast for small $a$. Since $E_{a,b}\left(  -1\right)  $ in
(\ref{Mittag_Leffler_function_Eab}) is an alternating series with decreasing
coefficients for $a>0$ and $b>1.462$, it holds that $\frac{1}{\Gamma\left(
b\right)  }<E_{a,b}\left(  -1\right)  <\frac{1}{\Gamma\left(  a+b\right)  }$.

The major interest of the expansion
(\ref{Ea,b(-z)_Taylor_in_a_around_b_in_mina}) lies in its fast convergence for
small $a$, whereas the definition (\ref{Mittag_Leffler_function_Eab}) is
converging slower for small $a$. Moreover, rewriting
(\ref{recursion_Mittag_Leffler_function_E_a,b}) as%
\[
E_{a,b}\left(  z\right)  =\frac{1}{n}\sum_{r=0}^{n-1}E_{\frac{a}{n},b}\left(
z^{\frac{1}{n}}e^{i\frac{2\pi r}{n}}\right)
\]
illustrates that any real $a$ can be transformed to a value smaller than $1$
by choosing $n=\left[  a\right]  +1$, where $\left[  a\right]  $ is the
largest integer smaller or equal to $a$. Indeed, for $\left\vert z\right\vert
<1$,%
\[
E_{a,b}\left(  z\right)  =\sum_{m=0}^{\infty}\frac{1}{m!}\left.  \frac{d^{m}%
}{du^{m}}\frac{1}{\Gamma\left(  u\right)  }\right\vert _{u=b}\left\{
\frac{\sum_{r=0}^{\left[  a\right]  }F_{-m-1}\left(  -\frac{\log z}{\left[
a\right]  +1}-i\pi\left(  \frac{2r}{\left[  a\right]  +1}+1\right)  \right)
}{\left[  a\right]  +1}\right\}  \left(  -\frac{a}{\left[  a\right]
+1}\right)  ^{m}%
\]

\medskip\refstepcounter{article}{\noindent\textbf{\thearticle. }%
}\ignorespaces\label{art_Mobius_inversion} \emph{M\"{o}bius inversion}. The
first M\"{o}bius inversion pair is%
\begin{equation}
g(x)=\sum_{n=1}^{\infty}f(nx)\hspace{1cm}\Longleftrightarrow\hspace
{1cm}f(x)=\sum_{n=1}^{\infty}\mu(n)g(nx) \label{Mobius_inversion_pair1}%
\end{equation}
where $\mu(n)$ is the M\"{o}bius function. Let $f\left(  x\right)
=\frac{z^{x}}{\Gamma\left(  b+x\right)  }$ in (\ref{Mobius_inversion_pair1}),
then%
\[
g\left(  x\right)  =\sum_{n=1}^{\infty}f(nx)=\sum_{n=1}^{\infty}\frac{z^{nx}%
}{\Gamma\left(  b+nx\right)  }=E_{x,b}\left(  z^{x}\right)  -\frac{1}%
{\Gamma\left(  b\right)  }%
\]
and M\"{o}bius inversion $f(x)=\sum_{n=1}^{\infty}\mu(n)g(nx)$ in
(\ref{Mobius_inversion_pair1}) yields, for $x\neq0$,%
\begin{equation}
\frac{z^{x}}{\Gamma\left(  b+x\right)  }=\sum_{n=1}^{\infty}\mu(n)\left(
E_{nx,b}\left(  z^{nx}\right)  -\frac{1}{\Gamma\left(  b\right)  }\right)
\label{Mittag_Leffler_mobius_inversion}%
\end{equation}
With $y=z^{x}$, (\ref{Mittag_Leffler_mobius_inversion}) simplifies to
$\frac{y}{\Gamma\left(  b+x\right)  }=\sum_{n=1}^{\infty}\mu(n)\left(
E_{nx,b}\left(  y^{n}\right)  -\frac{1}{\Gamma\left(  b\right)  }\right)  $.

Invoking $E_{2a,b}\left(  z^{2}\right)  =\frac{E_{a,b}\left(  z\right)
+E_{a,b}\left(  -z\right)  }{2}$ in \textbf{art}. \ref{art_even_and_odd} in
(\ref{Mittag_Leffler_mobius_inversion})
\begin{align*}
\frac{z^{x}}{\Gamma\left(  b+x\right)  }  &  =\sum_{n=1}^{\infty}\mu(n)\left(
2E_{2nx,b}\left(  z^{2nx}\right)  -E_{nx,b}\left(  -z^{nx}\right)  -\frac
{1}{\Gamma\left(  b\right)  }\right) \\
&  =2\sum_{n=1}^{\infty}\mu(n)\left(  E_{2nx,b}\left(  z^{2nx}\right)
-\frac{1}{\Gamma\left(  b\right)  }\right)  -\sum_{n=1}^{\infty}\mu(n)\left(
E_{nx,b}\left(  -z^{nx}\right)  -\frac{1}{\Gamma\left(  b\right)  }\right)
\end{align*}
and using (\ref{Mittag_Leffler_mobius_inversion}) leads to
\begin{equation}
\sum_{n=1}^{\infty}\mu(n)\left(  E_{nx,b}\left(  -z^{nx}\right)  -\frac
{1}{\Gamma\left(  b\right)  }\right)  =\frac{2z^{2x}}{\Gamma\left(
b+2x\right)  }-\frac{z^{x}}{\Gamma\left(  b+x\right)  }
\label{Mobius_alternating_Mittag_Leffler}%
\end{equation}
which is an instance of the general M\"{o}bius function identity, proved in
\cite{PVM_charcoef},
\begin{equation}
\sum_{j=1}^{n}\mu(j)\sum_{q=1}^{\left[  \frac{n}{j}\right]  }(-1)^{q}f\left(
qj\right)  =2f\left(  2\right)  -f\left(  1\right)
\label{Mobius_identity_alternating}%
\end{equation}
holds for any function $f$ and any $n>1$. With $y=z^{x}$,
(\ref{Mobius_alternating_Mittag_Leffler}) simplifies to $\frac{2y^{2}}%
{\Gamma\left(  b+2x\right)  }-\frac{y}{\Gamma\left(  b+x\right)  }=\sum
_{n=1}^{\infty}\mu(n)\left(  E_{nx,b}\left(  -y^{n}\right)  -\frac{1}%
{\Gamma\left(  b\right)  }\right)  $.

\medskip\refstepcounter{article}{\noindent\textbf{\thearticle. }%
}\ignorespaces\label{art_Mertens_function_and_MittagLeffler} \emph{Mertens
function}. Applying Abel summation using the Mertens\footnote{A sufficient
condition to prove the Riemann Hypothesis is to demonstate that the Mertens
function behaves as $\gamma_{-1}\left(  x\right)  =O\left(  x^{\frac{1}%
{2}+\varepsilon}\right)  $ for large $x$.} function $\gamma_{-1}\left(
k\right)  =\sum_{l=1}^{k}\mu(l)$, we obtain
\[
f(x)=\sum_{k=1}^{\infty}\gamma_{-1}\left(  k\right)  \,(g(kx)-g\left(  \left(
k+1\right)  x)\right)  )+\lim_{N\rightarrow\infty}g(Nx)\gamma_{-1}\left(
N\right)
\]
Hence,%
\[
\frac{z^{x}}{\Gamma\left(  b+x\right)  }=\sum_{k=1}^{\infty}\gamma_{-1}\left(
k\right)  \,\left(  E_{kx,b}\left(  z^{kx}\right)  -E_{kx+x,b}\left(
z^{kx+x}\right)  \right)  +\lim_{N\rightarrow\infty}\left(  E_{Nx,b}\left(
z^{Nx}\right)  -\frac{1}{\Gamma\left(  b\right)  }\right)  \gamma_{-1}\left(
N\right)
\]
and the limit vanishes if $x>0$. With $\int_{kx}^{\left(  k+1\right)  x}%
\frac{d}{da}E_{a,b}\left(  z^{a}\right)  da=E_{kx+x,b}\left(  z^{kx+x}\right)
-E_{kx,b}\left(  z^{kx}\right)  $, the corresponding integral representation
is%
\begin{align*}
\frac{z^{x}}{\Gamma\left(  b+x\right)  }  &  =-\sum_{k=1}^{\infty}\gamma
_{-1}\left(  k\right)  \,\int_{kx}^{\left(  k+1\right)  x}\frac{d}{da}%
E_{a,b}\left(  z^{a}\right)  da\\
&  =-\sum_{k=1}^{\infty}\,\int_{kx}^{\left(  k+1\right)  x}\gamma_{-1}\left(
\left[  \frac{a}{x}\right]  \right)  \frac{d}{da}E_{a,b}\left(  z^{a}\right)
da
\end{align*}
and%
\begin{equation}
\frac{z^{x}}{\Gamma\left(  b+x\right)  }=-\int_{x}^{\infty}\gamma_{-1}\left(
\left[  \frac{a}{x}\right]  \right)  \frac{d}{da}E_{a,b}\left(  z^{a}\right)
da\hspace{1cm}\text{for }x>0
\label{Mobius_inversion_Mertens_function_Mittag_Leffler_integral}%
\end{equation}
Similarly for (\ref{Mobius_alternating_Mittag_Leffler}), it holds that%
\begin{equation}
\frac{2z^{2x}}{\Gamma\left(  b+2x\right)  }-\frac{z^{x}}{\Gamma\left(
b+x\right)  }=-\int_{x}^{\infty}\gamma_{-1}\left(  \left[  \frac{a}{x}\right]
\right)  \frac{d}{da}E_{a,b}\left(  -z^{a}\right)  da\hspace{1cm}\text{for
}x>0 \label{Mobius_alternating_Mittag_Leffler_integral}%
\end{equation}
Although (\ref{Mobius_inversion_Mertens_function_Mittag_Leffler_integral}) and
(\ref{Mobius_alternating_Mittag_Leffler_integral}) are remarkable, further
progress to estimate the order of $\gamma_{-1}\left(  \left[  x\right]
\right)  $ requires a study of $\frac{d}{da}E_{a,b}\left(  yz^{a}\right)  $.

\section{Integrals containing $E_{a,b}\left(  z\right)  $}

\medskip\refstepcounter{article}{\noindent\textbf{\thearticle. }%
}\ignorespaces\label{art_integral_duplication_formula} \emph{Integral
duplication formula for }$E_{a,b}\left(  -z\right)  $. Using the duplication
formula of the Gamma function, $\Gamma\left(  2z\right)  =\frac{1}{\sqrt{\pi}%
}2^{2z-1}\Gamma\left(  z\right)  \Gamma\left(  z+\frac{1}{2}\right)  $, the
definition (\ref{Mittag_Leffler_function_Eab}) is rewritten as%
\[
E_{a,b}\left(  z\right)  =\sum_{k=0}^{\infty}\frac{z^{k}}{\Gamma\left(
b+ak\right)  }=\frac{2^{2b-1}}{\sqrt{\pi}}\sum_{k=0}^{\infty}\frac
{\Gamma\left(  \frac{1}{2}+b+ak\right)  \left(  2^{2a}z\right)  ^{k}}%
{\Gamma\left(  2b+2ak\right)  }%
\]
Invoking the Euler integral $\Gamma\left(  s\right)  =\int_{0}^{\infty}%
t^{s-1}e^{-t}dt$ for $\operatorname{Re}\left(  s\right)  >0$,%
\begin{align*}
E_{a,b}\left(  -z\right)   &  =\frac{2^{2b-1}}{\sqrt{\pi}}\sum_{k=0}^{\infty
}\frac{\left(  -2^{2a}z\right)  ^{k}}{\Gamma\left(  2b+2ak\right)  }\int
_{0}^{\infty}t^{b-\frac{1}{2}}t^{ak}e^{-t}dt\\
&  =\frac{2^{2b-1}}{\sqrt{\pi}}\int_{0}^{\infty}t^{b-\frac{1}{2}}e^{-t}%
\sum_{k=0}^{\infty}\frac{\left(  -z\left(  4t\right)  ^{a}\right)  ^{k}%
}{\Gamma\left(  2b+2ak\right)  }dt
\end{align*}
and the definition (\ref{Mittag_Leffler_function_Eab}), we find an integral
duplication formula for $\operatorname{Re}\left(  z\right)  \geq0$ and
$\operatorname{Re}\left(  b\right)  \geq-\frac{1}{2}$,%
\begin{equation}
E_{a,b}\left(  -z\right)  =\frac{2^{2b-1}}{\sqrt{\pi}}\int_{0}^{\infty
}t^{b-\frac{1}{2}}e^{-t}E_{2a,2b}\left(  -z\left(  4t\right)  ^{a}\right)  dt
\label{integral_recursion_E_a,b(-z)}%
\end{equation}

After substituting the Gamma duplication formula in the slightly rewritten
power series%
\begin{align*}
E_{a,b}\left(  z\right)   &  =\sum_{k=0}^{\infty}\frac{z^{k}}{\Gamma\left(
b+ak\right)  }=\sum_{k=0}^{\infty}\frac{z^{k}}{\left(  b-1+ak\right)
\Gamma\left(  b-1+ak\right)  }\\
&  =\frac{2^{2b-2}}{\sqrt{\pi}}\sum_{k=0}^{\infty}\frac{\left(  4^{a}z\right)
^{k}}{\Gamma\left(  2b-1+2ak\right)  }\Gamma\left(  b+ak-\frac{1}{2}\right)
\end{align*}
we find alternatively,%
\begin{align*}
E_{a,b}\left(  -z\right)   &  =\frac{2^{2b-2}}{\sqrt{\pi}}\sum_{k=0}^{\infty
}\frac{\left(  -4^{a}z\right)  ^{k}}{\Gamma\left(  2b-1+2ak\right)  }\int
_{0}^{\infty}t^{b-\frac{3}{2}}t^{ak}e^{-t}dt\\
&  =\frac{2^{2b-2}}{\sqrt{\pi}}\int_{0}^{\infty}t^{b-\frac{3}{2}}e^{-t}%
\sum_{k=0}^{\infty}\frac{\left(  -z\left(  4t\right)  ^{a}\right)  ^{k}%
}{\Gamma\left(  2b-1+2ak\right)  }dt
\end{align*}
for $\operatorname{Re}\left(  b\right)  \geq\frac{1}{2}$ and
$\operatorname{Re}\left(  z\right)  $,%
\begin{equation}
E_{a,b}\left(  -z\right)  =\frac{2^{2\left(  b-1\right)  }}{\sqrt{\pi}}%
\int_{0}^{\infty}t^{b-\frac{3}{2}}e^{-t}E_{2a,2b-1}\left(  -z\left(
4t\right)  ^{a}\right)  dt \label{integral_recursion_E_a,b(-z)_andere_b}%
\end{equation}
which also applies to $E_{a}\left(  z\right)  =E_{a,1}\left(  z\right)  $
after choosing $b=1$.

Let $u^{a}=z\left(  4t\right)  ^{a}$ or $u=4z^{\frac{1}{a}}t$, then
(\ref{integral_recursion_E_a,b(-z)_andere_b}) is, for real, nonnegative $z$,%
\[
E_{a,b}\left(  -z\right)  =\frac{2^{2b-2}}{\sqrt{\pi}}\left(  \frac{1}%
{4}z^{-\frac{1}{a}}\right)  ^{b-\frac{1}{2}}\int_{0}^{\infty}u^{b-\frac{3}{2}%
}e^{-\left(  \frac{1}{4}z^{-\frac{1}{a}}\right)  u}E_{2a,2b-1}\left(
-u^{a}\right)  du
\]
Let $s=\frac{1}{4}z^{-\frac{1}{a}}$ or $z=\left(  4s\right)  ^{-a}$, then we
arrive at the Laplace transform%
\begin{equation}
\int_{0}^{\infty}e^{-su}u^{b-\frac{3}{2}}E_{2a,2b-1}\left(  -u^{a}\right)
du=\frac{\sqrt{\pi}}{2^{2b-2}s^{b-\frac{1}{2}}}E_{a,b}\left(  -\left(
4s\right)  ^{-a}\right)  \label{Laplace_recursion_E_a,b}%
\end{equation}

\textbf{Example} For $b=1$ and $a=\frac{1}{2}$, the Laplace transform
(\ref{Laplace_recursion_E_a,b}) becomes%
\[
\frac{\sqrt{\pi}}{\sqrt{s}}E_{\frac{1}{2}}\left(  -\frac{1}{2\sqrt{s}}\right)
=\int_{0}^{\infty}u^{-\frac{1}{2}}e^{-su-\sqrt{u}}du=2\int_{0}^{\infty
}e^{-st^{2}-t}dt
\]
With $st^{2}+t=s\left(  t+\frac{1}{2s}\right)  ^{2}-\frac{1}{4s}$, we have%
\[
\frac{\sqrt{\pi}}{\sqrt{s}}E_{\frac{1}{2}}\left(  -\frac{1}{2\sqrt{s}}\right)
=2e^{\frac{1}{4s}}\int_{0}^{\infty}e^{-s\left(  t+\frac{1}{2s}\right)  ^{2}%
}dt=2e^{\frac{1}{4s}}\int_{\frac{1}{2s}}^{\infty}e^{-su^{2}}du=e^{\frac{1}%
{4s}}\frac{2}{\sqrt{s}}\int_{\frac{1}{2\sqrt{s}}}^{\infty}e^{-t^{2}}dt
\]
Simplified with erfc$\left(  x\right)  =\frac{2}{\sqrt{\pi}}\int_{x}^{\infty
}e^{-u^{2}}du=1-\operatorname{erf}\left(  x\right)  $ and $\operatorname{erf}%
\left(  -x\right)  =\frac{2}{\sqrt{\pi}}\int_{0}^{-x}e^{-u^{2}}%
du=-\operatorname{erf}\left(  x\right)  $%
\[
E_{\frac{1}{2}}\left(  -x\right)  =e^{x^{2}}\frac{2}{\sqrt{\pi}}\int
_{x}^{\infty}e^{-t^{2}}dt=e^{x^{2}}\text{erfc}\left(  x\right)
\]
is again (\ref{E_1/2,1}), because $\operatorname{erf}\left(  -x\right)
=\frac{2}{\sqrt{\pi}}\int_{0}^{-x}e^{-u^{2}}du=-\operatorname{erf}\left(
x\right)  $.

\medskip\refstepcounter{article}{\noindent\textbf{\thearticle. }%
}\ignorespaces\label{art_integral_multiplication_formula} \emph{Integral
multiplication formula for }$E_{a,b}\left(  z\right)  $. The method of
\textbf{art.} \ref{art_integral_duplication_formula} is readily generalized.
Invoking Gauss's multiplication formula (\ref{Gamma_multiplication_formula})
into the definition (\ref{Mittag_Leffler_function_Eab}) yields%
\[
E_{a,b}\left(  z\right)  =\sum_{k=0}^{\infty}\frac{z^{k}}{\Gamma\left(
b+ak\right)  }=\left(  2\pi\right)  ^{\frac{1}{2}\left(  1-n\right)
}n^{nb-\frac{1}{2}}\sum_{k=0}^{\infty}\frac{%
{\displaystyle\prod\limits_{j=1}^{n-1}}
\Gamma\left(  b+ak+\frac{j}{n}\right)  }{\Gamma\left(  nb+nak\right)  }\left(
n^{na}z\right)  ^{k}%
\]
We introduce the Mellin transform (\ref{Melling_prod_Gamma}) of a product of
Gamma functions for $\operatorname{Re}\left(  s\right)  >0$,%
\[%
{\displaystyle\prod\limits_{j=1}^{n-1}}
\Gamma\left(  s+\frac{j}{n}\right)  =\int_{0}^{\infty}u^{s-1}h_{n}\left(
u\right)  du
\]
where the inverse function $h_{n}\left(  u\right)  $ is specified in
(\ref{Mellin_inverse_prod_Gamma_multiplication_Taylor_series}) in
\textbf{art.} \ref{art_Mellin_transform_product_Gamma} as a Taylor series in
$u$, and we obtain%
\[
E_{a,b}\left(  z\right)  =\left(  2\pi\right)  ^{\frac{1}{2}\left(
1-n\right)  }n^{nb-\frac{1}{2}}\int_{0}^{\infty}h_{n}\left(  u\right)
u^{b-1}\sum_{k=0}^{\infty}\frac{\left(  u^{a}n^{na}z\right)  ^{k}}%
{\Gamma\left(  nb+nak\right)  }du
\]
Thus, for any integer $n$, we arrive at an integral multiplication formula for
the Mittag-Leffler functions,%
\begin{equation}
E_{a,b}\left(  z\right)  =\left(  2\pi\right)  ^{\frac{1}{2}\left(
1-n\right)  }n^{nb-\frac{1}{2}}\int_{0}^{\infty}h_{n}\left(  u\right)
u^{b-1}E_{na,nb}\left(  u^{a}n^{na}z\right)  du
\label{E_a,b_multiplciation_Gauss_all_n}%
\end{equation}
The companion of (\ref{E_a,b_multiplciation_Gauss_all_n}) follows similarly
from $E_{a,b}\left(  z\right)  =\sum_{k=0}^{\infty}\frac{z^{k}}{\left(
b-1+ak\right)  \Gamma\left(  b-1+ak\right)  }$ as%
\begin{equation}
E_{a,b}\left(  z\right)  =\left(  2\pi\right)  ^{\frac{1}{2}\left(
1-n\right)  }n^{n\left(  b-1\right)  +\frac{1}{2}}\int_{0}^{\infty}%
h_{n}\left(  u\right)  u^{b-2}E_{na,n\left(  b-1\right)  +1}\left(
u^{a}n^{na}z\right)  du \label{E_a,b_multiplciation_Gauss_all_n_companion}%
\end{equation}
which directly reduces for $b=1$ to the Mittag-Leffler function $E_{a}\left(
z\right)  =E_{a,1}\left(  z\right)  $,%
\[
E_{a}\left(  z\right)  =\left(  2\pi\right)  ^{\frac{1}{2}\left(  1-n\right)
}n^{\frac{1}{2}}\int_{0}^{\infty}h_{n}\left(  u\right)  u^{-1}E_{na}\left(
u^{a}n^{na}z\right)  du
\]

\medskip\refstepcounter{article}{\noindent\textbf{\thearticle. }%
}\ignorespaces\label{art_integral_multiplication_formula_special}
\emph{Special cases of the integral multiplication formula for }%
$E_{a,b}\left(  z\right)  $. The case $n=3$ in
(\ref{E_a,b_multiplciation_Gauss_all_n}) with $h_{3}\left(  u\right)  $ in
(\ref{h_3_BesselK}) expressed in terms of the modified Bessel function
$K_{\nu}\left(  z\right)  $, defined in (\ref{def_BesselK}), becomes%
\begin{equation}
E_{a,b}\left(  z\right)  =\frac{3^{3b-\frac{1}{2}}}{\pi}\int_{0}^{\infty
}K_{\frac{1}{3}}\left(  2\sqrt{u}\right)  u^{b-\frac{1}{2}}E_{3a,3b}\left(
3^{3a}u^{a}z\right)  du \label{E_a,b_BesselK_triple_order}%
\end{equation}
The companion of (\ref{E_a,b_BesselK_triple_order}) follows from
(\ref{E_a,b_multiplciation_Gauss_all_n_companion}) as%
\begin{equation}
E_{a,b}\left(  z\right)  =\frac{3^{3\left(  b-1\right)  +\frac{1}{2}}}{\pi
}\int_{0}^{\infty}K_{\frac{1}{3}}\left(  2\sqrt{u}\right)  u^{b-\frac{3}{2}%
}E_{3a,3b-2}\left(  3^{3a}u^{a}z\right)  du
\label{E_a,b_BesselK_triple_order_variant}%
\end{equation}

In \cite[eq. (25), (28) and (31)]{Apelblat2020}, Apelblat has recently
presented three remarkable integral functional relations,
\begin{align}
E_{a}\left(  t^{a}\right)   &  =\frac{1}{\sqrt{\pi t}}\int_{0}^{\infty
}e^{-\frac{u^{2}}{4t}}E_{2a}\left(  u^{2a}\right)  du\label{Apelblat_1}\\
E_{a}\left(  t^{a}\right)   &  =\frac{1}{\pi}\int_{0}^{\infty}\sqrt{\frac
{u}{t}}K_{\frac{1}{3}}\left(  \frac{2u^{\frac{3}{2}}}{\sqrt{27t}}\right)
E_{3a}\left(  u^{3a}\right)  du\label{Apelblat_2}\\
\frac{E_{a}\left(  t^{a}\right)  -1}{\sqrt{t}}  &  =\int_{0}^{\infty}%
\frac{J_{1}\left(  2\sqrt{tu}\right)  }{\sqrt{u}}E_{a}\left(  u^{a}\right)  du
\label{Apelblat_3}%
\end{align}
which he has skilfully derived by manipulations of Laplace transforms. Apart
from the last relation (\ref{Apelblat_3}), where $J_{p}\left(  z\right)  $ is
the Bessel function of order $p$, the first two are special cases of the
multiplication formula for\emph{ }$E_{a,b}\left(  z\right)  $ in \textbf{art.}
\ref{art_integral_multiplication_formula}. Indeed, after substitution of
\thinspace$x=\frac{1}{4st}u^{2}$ in the first integral (\ref{Apelblat_1}), we
obtain%
\[
E_{a,1}\left(  t^{a}\right)  =\frac{s^{\frac{1}{2}}}{\sqrt{\pi}}\int
_{0}^{\infty}x^{-\frac{1}{2}}e^{-sx}E_{2a,1}\left(  \left(  4st\right)
^{a}x^{a}\right)  dx
\]
which is a special case of (\ref{Laplace_recursion_E_a,b}), more easily
noticed from its generalization (\ref{Laplace_recursion_type_b+1}) below, for
$a\rightarrow2a$, $b=1$ and $x=\left(  4st\right)  ^{a}$. With $z=t^{a}$ and
after substituting $x=\frac{u^{3}}{27t}$ in
(\ref{E_a,b_BesselK_triple_order_variant}), we obtain%
\[
E_{a,b}\left(  t^{a}\right)  =\frac{t^{\frac{1}{2}-b}}{\pi}\int_{0}^{\infty
}K_{\frac{1}{3}}\left(  2\sqrt{\frac{u^{3}}{27t}}\right)  u^{3b-\frac{5}{2}%
}E_{3a,3b-2}\left(  u^{3a}\right)  dx
\]
which leads to (\ref{Apelblat_2}) by choosing $b=1$.

\medskip\refstepcounter{article}{\noindent\textbf{\thearticle. }%
}\ignorespaces\label{art_Laplace_transform} \emph{Laplace transform of
}$E_{a,b}\left(  xz^{\beta}\right)  $. The Laplace transform of $t^{\gamma
-1}E_{a,b}\left(  xt^{\beta}\right)  $ is%
\[
\int_{0}^{\infty}e^{-st}t^{\gamma-1}E_{a,b}\left(  xt^{\beta}\right)
dt=\int_{0}^{\infty}e^{-st}\sum_{k=0}^{\infty}\frac{x^{k}t^{\gamma-1+\beta k}%
}{\Gamma\left(  b+ak\right)  }dt=\sum_{k=0}^{\infty}\frac{x^{k}}{\Gamma\left(
b+ak\right)  }\int_{0}^{\infty}e^{-st}t^{\gamma+\beta k-1}dt
\]
Fubini's theorem states that the summation and integration can be reversed
provided the integrals $\int_{0}^{\infty}e^{-st}t^{\gamma}E_{a,b}\left(
xt^{\beta}\right)  dt$ and $\int_{0}^{\infty}e^{-st}t^{\gamma+\beta
k-1}dt=\frac{\Gamma\left(  \gamma+\beta k\right)  }{s^{\gamma+k\beta}}$ exist
and the resulting series converges, leading to%
\begin{equation}
\int_{0}^{\infty}e^{-st}t^{\gamma-1}E_{a,b}\left(  xt^{\beta}\right)
dt=\sum_{k=0}^{\infty}\frac{\Gamma\left(  \gamma+\beta k\right)  }%
{\Gamma\left(  b+ak\right)  }\frac{x^{k}}{s^{\gamma+k\beta}}
\label{Laplace_transform_general}%
\end{equation}
valid for $\left\vert \beta\right\vert \leq\left\vert a\right\vert $.

A first choice is $\beta=a$ and $\gamma=b$ in (\ref{Laplace_transform_general}%
), which restricts $s$ so that $\left\vert s^{a}\right\vert >\left\vert
x\right\vert $,%
\begin{equation}
\int_{0}^{\infty}e^{-st}t^{b-1}E_{a,b}\left(  xt^{a}\right)  dt=\frac{1}%
{s^{b}}\sum_{k=0}^{\infty}\left(  \frac{x}{s^{a}}\right)  ^{k}=\frac{s^{a-b}%
}{s^{a}-x} \label{Laplace_transform_positive_argument}%
\end{equation}
The Laplace transform (\ref{Laplace_transform_positive_argument}) plays a key
role in the theory of the Mittag-Leffler function.

Two other choices in (\ref{Laplace_transform_general}) follow after the
introduction of the duplication formula of the Gamma function%
\[
\frac{\Gamma\left(  \gamma+\beta k\right)  }{\Gamma\left(  b+ak\right)
}=\frac{\sqrt{\pi}}{2^{b-1}}\frac{\Gamma\left(  \gamma+\beta k\right)
}{2^{ak}\Gamma\left(  \frac{b}{2}+\frac{a}{2}k\right)  \Gamma\left(
\frac{b+1}{2}+\frac{a}{2}k\right)  }%
\]
as $\gamma=\frac{b}{2}$ and $\beta=\frac{a}{2}$ and $\gamma=\frac{b+1}{2}$ and
$\beta=\frac{a}{2}$, respectively. If $\gamma=\frac{b}{2}$ and $\beta=\frac
{a}{2}$, then we find
\[
\int_{0}^{\infty}e^{-st}t^{\frac{b}{2}-1}E_{a,b}\left(  xt^{\frac{a}{2}%
}\right)  dt=\frac{\sqrt{\pi}}{2^{b-1}s^{\frac{b}{2}}}\sum_{k=0}^{\infty}%
\frac{1}{\Gamma\left(  \frac{b+1}{2}+\frac{a}{2}k\right)  }\left(  \frac
{x}{2^{a}s^{\frac{a}{2}}}\right)  ^{k}%
\]
resulting, with the definition (\ref{Mittag_Leffler_function_Eab}), in%
\begin{equation}
\int_{0}^{\infty}e^{-st}t^{\frac{b}{2}-1}E_{a,b}\left(  xt^{\frac{a}{2}%
}\right)  dt=\frac{\sqrt{\pi}}{2^{b-1}s^{\frac{b}{2}}}E_{\frac{a}{2}%
,\frac{b+1}{2}}\left(  \frac{x}{\left(  4s\right)  ^{\frac{a}{2}}}\right)
\label{Laplace_recursion_type_b+1}%
\end{equation}
while the third choice $\gamma=\frac{b+1}{2}$ and $\beta=\frac{a}{2}$ leads to%
\begin{equation}
\int_{0}^{\infty}e^{-st}t^{\frac{b+1}{2}-1}E_{a,b}\left(  xt^{\frac{a}{2}%
}\right)  dt=\frac{\sqrt{\pi}}{2^{b-1}s^{\frac{b+1}{2}}}E_{\frac{a}{2}%
,\frac{b}{2}}\left(  \frac{x}{\left(  4s\right)  ^{\frac{a}{2}}}\right)
\label{Laplace_recursion_type_b}%
\end{equation}
Both (\ref{Laplace_recursion_type_b+1}) and (\ref{Laplace_recursion_type_b})
are slightly more general than and reduce to (\ref{Laplace_recursion_E_a,b})
and (\ref{integral_recursion_E_a,b(-z)}), respectively, for $x=-1$.

\medskip\refstepcounter{article}{\noindent\textbf{\thearticle. }%
}\ignorespaces\label{art_generalized_integration} \emph{Generalized
integration.} By using a variation of the integral of the Beta-function
\cite[6.2]{Abramowitz}, $\int_{0}^{x}u^{z-1}\left(  x-u\right)  ^{w-1}%
du=x^{z+w-1}\frac{\Gamma\left(  z\right)  \Gamma\left(  w\right)  }%
{\Gamma\left(  z+w\right)  }$ for real $x>0$, $\operatorname{Re}\left(
z\right)  >0$ and $\operatorname{Re}\left(  w\right)  >0$, we obtain a
generalized integral variant of the Mittag-Leffler function $E_{a,b}\left(
z\right)  $ in (\ref{Mittag_Leffler_function_Eab}),%
\begin{align*}
\frac{1}{\Gamma\left(  w\right)  }\int_{0}^{x}\left(  x-u\right)
^{w-1}u^{\gamma-1}E_{a,b}\left(  \lambda u^{\beta}\right)  du  &  =\frac
{1}{\Gamma\left(  w\right)  }\sum_{k=0}^{\infty}\frac{\lambda^{k}}%
{\Gamma\left(  b+ak\right)  }\int_{0}^{x}u^{\beta k+\gamma-1}\left(
x-u\right)  ^{w-1}du\\
&  =x^{\gamma+w-1}\sum_{k=0}^{\infty}\frac{\lambda^{k}x^{\beta k}}%
{\Gamma\left(  b+ak\right)  }\frac{\Gamma\left(  \gamma+\beta k\right)
}{\Gamma\left(  \gamma+w+\beta k\right)  }%
\end{align*}
which reduces for $\beta=a$ and $\gamma=b$ to%
\begin{equation}
\frac{1}{\Gamma\left(  w\right)  }\int_{0}^{x}\left(  x-u\right)
^{w-1}u^{b-1}E_{a,b}\left(  \lambda u^{a}\right)  du=x^{b-1+w}E_{a,b+w}\left(
\lambda x^{a}\right)  \label{Mittag_Leffler_generalized_integration}%
\end{equation}
The $m$-fold integral\footnote{The $n$-fold integral is%
\[
F_{n}\left(  x,a\right)  =\int_{a}^{x}du_{1}\int_{a}^{u_{1}}du_{2}\ldots
\int_{a}^{u_{n-1}}du_{n}f\left(  u_{n}\right)  =\frac{1}{\left(  n-1\right)
!}\int_{a}^{x}\left(  x-u\right)  ^{n-1}f\left(  u\right)  du
\]
\par
The generalization towards fractional calculus, where the integer $n$ is
extended to a real number, is treated in \cite{Gorenflo_Mainardi_2000}.}
(\ref{Mittag_Leffler_generalized_integration}) for $w=m$ and $\lambda=1$
possesses the same form as the $m$-fold differentiation in
(\ref{Dm_(z^(b-1)E_a,b(z^a))})%
\[
\frac{d^{m}}{dx^{m}}\left\{  x^{b-1}E_{a,b}\left(  x^{a}\right)  \right\}
=x^{b-1-m}E_{a,b-m}\left(  x^{a}\right)
\]
that is better recognized with Cauchy's integral for the $m$-th derivative of
an analytic function, $\left.  \frac{d^{m}f(z)}{dz^{m}}\right\vert _{z=z_{0}%
}=\frac{m!}{2\pi i}\int_{C(z_{0})}\frac{f(\omega)\;d\omega}{(\omega
-z_{0})^{m+1}}$, where $C\left(  z_{0}\right)  $ is a closed contour around
$z_{0}$,%
\[
\frac{\left(  -1\right)  ^{m+1}m!}{2\pi i}\int_{C(x)}(x-\omega)^{-m-1}%
\omega^{b-1}E_{a,b}\left(  \omega^{a}\right)  \;d\omega=x^{b-1-m}%
E_{a,b-m}\left(  x^{a}\right)
\]
Hence, (\ref{Mittag_Leffler_generalized_integration}) written with the
reflection formula (\ref{Gamma_reflection_formula}), suggests that%
\[
\frac{\frac{\sin\pi w}{\pi}\Gamma\left(  1-w\right)  }{2\pi i}\int
_{C(x)}(x-\omega)^{w-1}\omega^{b-1}E_{a,b}\left(  \omega^{a}\right)
\;d\omega=x^{b-1-w}E_{a,b-w}\left(  x^{a}\right)
\]
holds for any negative real $w$, leading to fractional derivatives. The
Laplace transform (\ref{Laplace_transform_positive_argument}) easier connects
to fractional derivatives, avoiding the contour integral. Applications of the
Mittag-Leffler function to fractional calculus are amply illustrated in
\cite{Gorenflo_2014}.

\textbf{Example} For $a=b=1$ in (\ref{Mittag_Leffler_generalized_integration})
and, next, $x=1$, $\lambda=z$ and $w+1\rightarrow b$, we obtain, for
$\operatorname{Re}\left(  b\right)  >1$,%
\[
E_{1,b}\left(  z\right)  =\frac{1}{\Gamma\left(  b-1\right)  }\int_{0}%
^{1}\left(  1-u\right)  ^{b-2}e^{zu}du
\]
Let $t=1-u$, then the incomplete Gamma function appears%
\[
E_{1,b}\left(  z\right)  =\frac{e^{z}}{\Gamma\left(  b-1\right)  }\int_{0}%
^{1}t^{b-2}e^{-zt}dt=\frac{z^{1-b}e^{z}}{\Gamma\left(  b-1\right)  }\int
_{0}^{z}u^{b-2}e^{-u}du
\]
which is again equal to (\ref{E_1,b}).

\medskip\refstepcounter{article}{\noindent\textbf{\thearticle. }%
}\ignorespaces\label{art_integration_product_Mittag_leffler_fnctions}
\emph{Integration of a product of Mittag-Leffler functions.} We extend the
idea in \textbf{art}. \ref{art_generalized_integration} and consider, for
$x\geq0$,%
\begin{align*}
L  &  =\frac{1}{\Gamma\left(  w\right)  }\int_{0}^{x}\left(  x-u\right)
^{w-1}u^{\beta-1}E_{a,b}\left(  \lambda u^{\alpha}\right)  E_{c,d}\left(
\mu\left(  x-u\right)  ^{\gamma}\right)  du\\
&  =\frac{1}{\Gamma\left(  w\right)  }\sum_{k=0}^{\infty}\sum_{m=0}^{\infty
}\frac{\lambda^{k}}{\Gamma\left(  b+ak\right)  }\frac{\mu^{m}}{\Gamma\left(
d+cm\right)  }\int_{0}^{x}\left(  x-u\right)  ^{w+m\gamma-1}u^{\beta
+k\alpha-1}du
\end{align*}
Introducing the integral of the Beta-function results in%
\[
L=\frac{x^{\beta+w-1}}{\Gamma\left(  w\right)  }\sum_{k=0}^{\infty}\sum
_{m=0}^{\infty}\frac{\left(  \lambda x^{\alpha}\right)  ^{k}}{\Gamma\left(
b+ak\right)  }\frac{\left(  \mu x^{\gamma}\right)  ^{m}}{\Gamma\left(
d+cm\right)  }\frac{\Gamma\left(  \beta+k\alpha\right)  \Gamma\left(
w+m\gamma\right)  }{\Gamma\left(  \beta+k\alpha+w+m\gamma\right)  }%
\]
After the choice $\alpha=a$, $\beta=b$, $\gamma=c$ and $d=w$, the double sum
simplifies to%
\[
L=\frac{x^{b+w-1}}{\Gamma\left(  w\right)  }\sum_{k=0}^{\infty}\sum
_{m=0}^{\infty}\frac{\left(  \lambda x^{a}\right)  ^{k}\left(  \mu
x^{c}\right)  ^{m}}{\Gamma\left(  b+w+ka+mc\right)  }%
\]
Further computations require the choice $c=a$,
\[
L=\frac{x^{b+w-1}}{\Gamma\left(  w\right)  }\sum_{k=0}^{\infty}\sum
_{m=0}^{\infty}\frac{\left(  \lambda x^{a}\right)  ^{k}\left(  \mu
x^{a}\right)  ^{m}}{\Gamma\left(  b+w+a\left(  k+m\right)  \right)  }%
\]
Let $q=k+m$, then $0\leq q$ and $m=q-k\geq0$, while $k\geq0$, thus%
\[
L=\frac{x^{b+w-1}}{\Gamma\left(  w\right)  }\sum_{q=0}^{\infty}\sum_{k=0}%
^{q}\frac{\lambda^{k}\mu^{q-k}x^{aq}}{\Gamma\left(  b+w+aq\right)  }%
=\frac{x^{b+w-1}}{\Gamma\left(  w\right)  }\sum_{q=0}^{\infty}\frac{\mu
^{q}x^{aq}}{\Gamma\left(  b+w+aq\right)  }\sum_{k=0}^{q}\left(  \frac{\lambda
}{\mu}\right)  ^{k}%
\]
Executing the finite geometric series $\sum_{k=0}^{q}\left(  \frac{\lambda
}{\mu}\right)  ^{k}=\frac{\left(  \frac{\lambda}{\mu}\right)  ^{q+1}-1}%
{\frac{\lambda}{\mu}-1}$ leads to%
\[
L=\frac{x^{b+w-1}}{\Gamma\left(  w\right)  \left(  \lambda-\mu\right)
}\left\{  \lambda\sum_{q=0}^{\infty}\frac{\left(  \lambda x^{a}\right)  ^{q}%
}{\Gamma\left(  b+w+aq\right)  }-\mu\sum_{q=0}^{\infty}\frac{\left(  \mu
x^{a}\right)  ^{q}}{\Gamma\left(  b+w+aq\right)  }\right\}
\]
Finally, we arrive for $\beta>0$ and $w>0$ at%
\begin{equation}
\int_{0}^{x}\left(  x-u\right)  ^{w-1}u^{b-1}E_{a,b}\left(  \lambda
u^{a}\right)  E_{a,w}\left(  \mu\left(  x-u\right)  ^{a}\right)
du=x^{b+w-1}\frac{\lambda E_{a,b+w}\left(  \lambda x^{a}\right)  -\mu
E_{a,b+w}\left(  \mu x^{a}\right)  }{\lambda-\mu}
\label{Mittag_Leffler_integration_product_Ea_Ea}%
\end{equation}
Clearly\footnote{In case $\mu\rightarrow\lambda$, then, after using de
l'Hospital's rule,%
\[
\lim_{\mu\rightarrow\lambda}\frac{\lambda E_{a,b+w}\left(  \lambda
x^{a}\right)  -\mu E_{a,b+w}\left(  \mu x^{a}\right)  }{\lambda-\mu}%
=E_{a,b+w}\left(  \lambda x^{a}\right)  +\lambda x^{a}\left.  \frac
{dE_{a,b+w}\left(  z\right)  }{dz}\right\vert _{z=\lambda x^{a}}%
\]
and the differential rule (\ref{D_E_a,b(z)}), formula
(\ref{Mittag_Leffler_integration_product_Ea_Ea}) becomes
\[
\int_{0}^{x}\left(  x-u\right)  ^{w-1}u^{b-1}E_{a,b}\left(  \lambda u^{\alpha
}\right)  E_{a,w}\left(  \lambda\left(  x-u\right)  ^{a}\right)
du=\frac{x^{b+w-1}}{a}\left\{  E_{a,b+w-1}\left(  \lambda x^{a}\right)
+\left(  a+1-\left(  b+w\right)  \right)  E_{a,b+w}\left(  \lambda
x^{a}\right)  \right\}
\]
}, since $\lim_{\mu\rightarrow0}E_{a,w}\left(  \mu\left(  x-u\right)
^{a}\right)  =\frac{1}{\Gamma\left(  w\right)  }$, then
(\ref{Mittag_Leffler_integration_product_Ea_Ea}) reduces to
(\ref{Mittag_Leffler_generalized_integration}).

\medskip\refstepcounter{article}{\noindent\textbf{\thearticle. }%
}\ignorespaces\label{art_Beta_function_like_asymptotic_result_Gauss} \emph{An
asymptotic result}. We invoke the device\footnote{Starting from the above
Beta-integral, it holds that%
\[
\Gamma\left(  z\right)  \lim_{x\rightarrow\infty}\frac{x^{z}\Gamma\left(
x+1\right)  }{\Gamma\left(  z+x+1\right)  }=\lim_{x\rightarrow\infty}\int
_{0}^{x}u^{z-1}\left(  1-\frac{u}{x}\right)  ^{x}du=\int_{0}^{\infty}%
u^{z-1}e^{-u}du=\Gamma\left(  z\right)
\]
from which%
\[
\frac{\Gamma\left(  z+x\right)  }{\Gamma\left(  x\right)  }\sim x^{z}%
\]
}, used by Gauss in his grand treatise \cite[p. 146]{Gauss1813} on the
hypergeometric function \cite[p. 74]{Klein_1933} to deduce the Euler integral
for the Gamma function from the Beta-integral based on%
\[
\lim_{n\rightarrow\infty}\left(  1-\frac{t}{n}\right)  ^{xn}=e^{-xt}%
\]
and start from (\ref{Mittag_Leffler_generalized_integration})%
\[
\int_{0}^{x}\left(  1-\frac{u}{x}\right)  ^{w-1}u^{b-1}E_{a,b}\left(  \lambda
u^{a}\right)  du=\Gamma\left(  w\right)  x^{b}E_{a,b+w}\left(  \lambda
x^{a}\right)
\]
Let $w=1+xs$ with $\operatorname{Re}\left(  s\right)  >0$ and $x$ is real,
then%
\[
\int_{0}^{x}\left(  1-\frac{u}{x}\right)  ^{xs}u^{b-1}E_{a,b}\left(  \lambda
u^{a}\right)  du=\Gamma\left(  xs+1\right)  x^{b}E_{a,b+xs+1}\left(  \lambda
x^{a}\right)
\]
and, after taking the limit for $x\rightarrow\infty$%
\[
\int_{0}^{\infty}e^{-su}u^{b-1}E_{a,b}\left(  \lambda u^{a}\right)
du=\lim_{x\rightarrow\infty}\Gamma\left(  xs+1\right)  x^{b}E_{a,b+xs+1}%
\left(  \lambda x^{a}\right)
\]
Comparison with the Laplace transform
(\ref{Laplace_transform_positive_argument}) for $\left\vert s^{a}\right\vert
>\left\vert \lambda\right\vert $ shows that%
\[
\lim_{x\rightarrow\infty}\Gamma\left(  xs+1\right)  x^{b}E_{a,b+xs+1}\left(
\lambda x^{a}\right)  =\frac{s^{a-b}}{s^{a}-\lambda}%
\]
and, for large real $x$,%
\[
E_{a,b+xs}\left(  \lambda x^{a}\right)  \sim\frac{1}{1-\frac{\lambda}{s^{a}}%
}\frac{1}{\Gamma\left(  xs\right)  \left(  sx\right)  ^{b}}%
\]
Let $z=\lambda x^{a}$, obeying $\left\vert xs\right\vert ^{a}>\left\vert
z\right\vert $, then for real $r=x\left\vert s\right\vert \rightarrow\infty$%
\[
E_{a,b+re^{i\theta}}\left(  z\right)  \sim\frac{r^{-b}e^{-bi\theta}}%
{1-zr^{-a}e^{-ai\theta}}\frac{1}{\Gamma\left(  re^{i\theta}\right)  }\sim
\frac{1}{\sqrt{2\pi}}\frac{r^{\frac{1}{2}-b}e^{-bi\theta}}{1-zr^{-a}%
e^{-ai\theta}}e^{-r\left(  \ln\left(  r\right)  -1\right)  \cos\theta+\theta
r\sin\theta}\left(  1+O\left(  \frac{1}{r}\right)  \right)
\]
where in the last step (\ref{asymptotic_Gamma(b+rexp(it))}) is used. Thus for
$a>0$ and for $-\frac{\pi}{2}<\theta<\frac{\pi}{2}$, we evidently find that
$E_{a,b}\left(  z\right)  \rightarrow0$ for $\operatorname{Re}\left(
b\right)  \rightarrow\infty$ and for $\frac{\pi}{2}<\theta<\frac{3\pi}{2}$
that $E_{a,b}\left(  z\right)  \rightarrow\infty$ for $\operatorname{Re}%
\left(  b\right)  \rightarrow-\infty$. For $\theta=\pm\frac{\pi}{2}$, it holds
that $E_{a,b+ir}\left(  z\right)  \sim\frac{1}{\sqrt{2\pi}}\frac{r^{\frac
{1}{2}-b}e^{-bi\frac{\pi}{2}}}{1-zr^{-a}e^{-ai\frac{\pi}{2}}}e^{\frac{\pi}%
{2}r}$, from which $\left\vert E_{a,\pm ib}\left(  z\right)  \right\vert
\sim\frac{b^{\frac{1}{2}}}{\sqrt{2\pi}}$ for real $b\rightarrow\infty$.

\medskip\refstepcounter{article}{\noindent\textbf{\thearticle. }%
}\ignorespaces\label{art_Apelblat_series} \emph{Apelblat series.} Inspired by
series in \cite{Apelblat2020}, we call
\[
g_{a,b}\left(  t\right)  =\sum_{k=0}^{\infty}f_{k}t^{b+k-1}E_{a,b+k}\left(
xt^{a}\right)
\]
an Apelblat series, where the Taylor series of the function $f\left(
z\right)  =\sum_{k=0}^{\infty}f_{k}z^{k}$ around $z_{0}=0$ converges for
$\left\vert z\right\vert \leq R$. Evidently, if $f_{k}=0$ for $k>n$, then
$f\left(  z\right)  $ is a polynomial of order $n$ in $z$ and $f\left(
z\right)  $ is an entire function. We generalize the method of Apelblat
\cite{Apelblat2020}. We take the Laplace transform $\mathcal{L}\left[
.\right]  $ of both sides and use (\ref{Laplace_transform_positive_argument}),%
\[
\mathcal{L}\left[  g_{a,b}\left(  t\right)  \right]  =\sum_{k=0}^{\infty}%
f_{k}\mathcal{L}\left[  t^{b+k-1}E_{a,b+k}\left(  xt^{a}\right)  \right]
=\sum_{k=0}^{\infty}f_{k}\frac{s^{a-b-k}}{s^{a}-x}=\frac{s^{a-b+\beta}}%
{s^{a}-x}\sum_{k=0}^{\infty}f_{k}s^{-k-\beta}%
\]
and obtain a product of Laplace transformed functions%
\begin{equation}
\mathcal{L}\left[  g_{a,b}\left(  t\right)  \right]  =\mathcal{L}\left[
t^{b-\beta-1}E_{a,b-\beta}\left(  xt^{a}\right)  \right]  \frac{1}{s^{\beta}%
}f\left(  \frac{1}{s}\right)  \label{Apelblat_series_Laplace_transform}%
\end{equation}
The convolution theorem for the Laplace transform suggests us to find the
inverse Laplace transform $\mathcal{L}^{-1}\left[  \frac{1}{s^{\beta}}f\left(
\frac{1}{s}\right)  \right]  $ of $\frac{1}{s^{\beta}}f\left(  \frac{1}%
{s}\right)  $. Apelblat \cite{Apelblat2020} observes and demonstrates that
elegant series follow if a closed form for $\mathcal{L}^{-1}\left[  \frac
{1}{s^{\beta}}f\left(  \frac{1}{s}\right)  \right]  $ exist, else we can
proceed with (\ref{inverse_Laplace}) and Hankel's integral
(\ref{Hankel_e^(aw)_a^(1-z)_op_Gamma(z)})%
\begin{align*}
\mathcal{L}^{-1}\left[  \frac{1}{s^{\beta}}f\left(  \frac{1}{s}\right)
\right]   &  =\frac{1}{2\pi i}\int_{c-i\infty}^{c+i\infty}\frac{1}{s^{\beta}%
}f\left(  \frac{1}{s}\right)  e^{st}ds=\frac{1}{2\pi i}\sum_{k=0}^{\infty
}f_{k}\int_{c-i\infty}^{c+i\infty}s^{-k-\beta}e^{st}ds\\
&  =t^{\beta-1}\sum_{k=0}^{\infty}f_{k}\frac{t^{k}}{\Gamma\left(
\beta+k\right)  }%
\end{align*}
After taking the inverse Laplace transform of both sides in
(\ref{Apelblat_series_Laplace_transform}), we formally arrive at the Apelblat
series, for the free parameter $\operatorname{Re}\left(  \beta\right)  >0$,
\begin{equation}
\sum_{k=0}^{\infty}f_{k}t^{b+k-1}E_{a,b+k}\left(  xt^{a}\right)  =\int_{0}%
^{t}\left(  t-u\right)  ^{b-\beta-1}E_{a,b-\beta}\left(  x\left(  t-u\right)
^{a}\right)  u^{\beta-1}\sum_{k=0}^{\infty}f_{k}\frac{u^{k}}{\Gamma\left(
\beta+k\right)  }du \label{series_Apelblat}%
\end{equation}
which directly follows from (\ref{Mittag_Leffler_generalized_integration}) for
$w=\beta+k$ and $b\rightarrow b-\beta$. Hence, the property
(\ref{series_Apelblat}) of the Apelblat series is a consequence of the
generalized integration property in \textbf{art}.
\ref{art_generalized_integration}.

\textbf{Examples} \textbf{1.} The Taylor series of $\left(  1+cz\right)
^{p}=\sum_{k=0}^{\infty}\binom{p}{k}c^{k}z^{k}$ converges for all complex $p$
provided $\left\vert z\right\vert <1$. The right-hand series in
(\ref{series_Apelblat}) becomes with $f_{k}=\binom{p}{k}c^{k}=\frac
{\Gamma\left(  p+1\right)  }{k!\Gamma\left(  p-k+1\right)  }c^{k}=\frac
{\Gamma\left(  -p+k\right)  }{k!\Gamma\left(  -p\right)  }\left(  -c\right)
^{k}$%
\[
\sum_{k=0}^{\infty}f_{k}\frac{u^{k}}{\Gamma\left(  \beta+k\right)  }=\frac
{1}{\Gamma\left(  -p\right)  }\sum_{k=0}^{\infty}\frac{\Gamma\left(
-p+k\right)  }{k!\Gamma\left(  \beta+k\right)  }\left(  -cu\right)  ^{k}%
\]
which reduces if we choose $\beta=-p$ to $\sum_{k=0}^{\infty}f_{k}\frac{u^{k}%
}{\Gamma\left(  \beta+k\right)  }=\frac{1}{\Gamma\left(  -p\right)  }e^{-cu}$.
The Apelblat series (\ref{series_Apelblat}) then becomes, for $q=-p$ and
$0<\operatorname{Re}\left(  q\right)  <\operatorname{Re}\left(  b\right)  $,
\[
\sum_{k=0}^{\infty}\binom{-q}{k}c^{k}t^{b+k-1}E_{a,b+k}\left(  xt^{a}\right)
=\frac{1}{\Gamma\left(  q\right)  }\int_{0}^{t}u^{b-q-1}E_{a,b-q}\left(
xu^{a}\right)  \left(  t-u\right)  ^{q-1}e^{-c\left(  t-u\right)  }du
\]
For $q=1$, $b=2$, $c=x=1$, we retrieve the series \cite[eq. (64)]%
{Apelblat2020} and for $q=2$, $b=2$, $c=x=1$, the series \cite[eq.
(68)]{Apelblat2020}.

\textbf{2.} The Taylor series of the Bessel function $\left(  \frac{z}%
{2}\right)  ^{-p}J_{p}\left(  z\right)  =\sum_{k=0}^{\infty}\frac{\left(
-1\right)  ^{k}}{k!\Gamma\left(  p+1+k\right)  }\left(  \frac{z}{2}\right)
^{2k}$ in \cite[9.1.10]{Abramowitz} has only even Taylor coefficients
$j_{2k}=\frac{\left(  -\frac{1}{4}\right)  ^{k}}{k!\Gamma\left(  p+1+k\right)
}$ and the odd $j_{2k+1}=0$. If $f_{2k+1}=0$, then%
\[
\sum_{k=0}^{\infty}f_{k}\frac{u^{k}}{\Gamma\left(  \beta+k\right)  }%
=\sum_{k=0}^{\infty}f_{2k}\frac{u^{2k}}{\Gamma\left(  \beta+2k\right)  }%
=\frac{\sqrt{\pi}}{2^{\beta-1}}\sum_{k=0}^{\infty}f_{2k}\frac{\left(  \frac
{u}{2}\right)  ^{2k}}{\Gamma\left(  \frac{\beta}{2}+k\right)  \Gamma\left(
\frac{\beta+1}{2}+k\right)  }%
\]
where the duplication formula of Gamma function is used. If $f_{2k}=\binom
{p}{k}c^{2k}=\frac{\Gamma\left(  -p+k\right)  }{k!\Gamma\left(  -p\right)
}\left(  -c^{2}\right)  ^{k}$ then
\[
\sum_{k=0}^{\infty}f_{k}\frac{u^{k}}{\Gamma\left(  \beta+k\right)  }%
=\frac{\sqrt{\pi}}{2^{\beta-1}\Gamma\left(  -p\right)  }\sum_{k=0}^{\infty
}\frac{\Gamma\left(  -p+k\right)  \left(  -1\right)  ^{k}\left(  \frac{cu}%
{2}\right)  ^{2k}}{k!\Gamma\left(  \frac{\beta}{2}+k\right)  \Gamma\left(
\frac{\beta+1}{2}+k\right)  }%
\]
Choosing $-p=\frac{\beta}{2}$ yields%
\[
\sum_{k=0}^{\infty}f_{k}\frac{u^{k}}{\Gamma\left(  \beta+k\right)  }%
=\frac{\sqrt{\pi}}{2^{\beta-1}\Gamma\left(  \frac{\beta}{2}\right)  }%
\sum_{k=0}^{\infty}\frac{\left(  -1\right)  ^{k}\left(  \frac{cu}{2}\right)
^{2k}}{k!\Gamma\left(  \frac{\beta-1}{2}+1+k\right)  }=\frac{\sqrt{\pi}\left(
\frac{cu}{2}\right)  ^{-\frac{\beta-1}{2}}}{2^{\beta-1}\Gamma\left(
\frac{\beta}{2}\right)  }J_{\frac{\beta-1}{2}}\left(  cu\right)
\]
and%
\[
\sum_{k=0}^{\infty}\binom{-\frac{\beta}{2}}{k}c^{2k}t^{b+2k-1}E_{a,b+2k}%
\left(  xt^{a}\right)  =\frac{\sqrt{\pi}\left(  \frac{2}{c}\right)
^{\frac{\beta-1}{2}}}{2^{\beta-1}\Gamma\left(  \frac{\beta}{2}\right)  }%
\int_{0}^{t}\left(  t-u\right)  ^{b-\beta-1}E_{a,b-\beta}\left(  x\left(
t-u\right)  ^{a}\right)  u^{\frac{\beta-1}{2}}J_{\frac{\beta-1}{2}}\left(
cu\right)  du
\]
while choosing $-p=\frac{\beta+1}{2}$ yields%
\[
\sum_{k=0}^{\infty}f_{k}\frac{u^{k}}{\Gamma\left(  \beta+k\right)  }%
=\frac{\sqrt{\pi}}{2^{\beta-1}\Gamma\left(  \frac{\beta+1}{2}\right)  }%
\sum_{k=0}^{\infty}\frac{\left(  -1\right)  ^{k}\left(  \frac{cu}{2}\right)
^{2k}}{k!\Gamma\left(  \frac{\beta}{2}+k\right)  }=\frac{\sqrt{\pi}\left(
\frac{cu}{2}\right)  ^{-\frac{\beta-2}{2}}}{2^{\beta-1}\Gamma\left(
\frac{\beta+1}{2}\right)  }J_{\frac{\beta-2}{2}}\left(  cu\right)
\]
and the Apelblat series (\ref{series_Apelblat}) then becomes%
\[
\sum_{k=0}^{\infty}\binom{-\frac{\beta+1}{2}}{k}c^{2k}t^{b+2k-1}%
E_{a,b+2k}\left(  xt^{a}\right)  =\frac{\sqrt{\pi}\left(  \frac{c}{2}\right)
^{1-\frac{\beta}{2}}}{2^{\beta-1}\Gamma\left(  \frac{\beta+1}{2}\right)  }%
\int_{0}^{t}\left(  t-u\right)  ^{b-\beta-1}E_{a,b-\beta}\left(  x\left(
t-u\right)  ^{a}\right)  u^{\frac{\beta}{2}}J_{\frac{\beta}{2}-1}\left(
cu\right)  du
\]
For $\beta=0$, $b=1$, $c=x=1$, we retrieve the series \cite[eq. (74)]%
{Apelblat2020},
\[
\sum_{k=0}^{\infty}\binom{-\frac{1}{2}}{k}t^{2k}E_{a,1+2k}\left(
t^{a}\right)  =-\int_{0}^{t}E_{a}\left(  u^{a}\right)  J_{1}\left(
t-u\right)  du
\]

\textbf{3.} Let us now consider the hypergeometric function \cite[15.1]%
{Abramowitz} with Taylor series around the origin,%
\begin{equation}
F\left(  p,q;r;z\right)  =\frac{\Gamma\left(  r\right)  }{\Gamma\left(
p\right)  \Gamma\left(  q\right)  }\sum_{k=0}^{\infty}\frac{\Gamma\left(
p+k\right)  \Gamma\left(  q+k\right)  }{\Gamma\left(  r+k\right)  k!}z^{k}
\label{hypergeometric_z_Gauss}%
\end{equation}
With (\ref{hypergeometric_z_Gauss}), the Apelblat series
(\ref{series_Apelblat}) at the right-hand side becomes%
\[
\sum_{k=0}^{\infty}f_{k}\frac{c^{k}u^{k}}{\Gamma\left(  \beta+k\right)
}=\frac{\Gamma\left(  r\right)  }{\Gamma\left(  p\right)  \Gamma\left(
q\right)  }\sum_{k=0}^{\infty}\frac{\Gamma\left(  p+k\right)  \Gamma\left(
q+k\right)  }{\Gamma\left(  r+k\right)  \Gamma\left(  \beta+k\right)  }%
\frac{c^{k}u^{k}}{k!}%
\]
and, if we choose $q$ equal to $\beta$, then%
\[
\sum_{k=0}^{\infty}f_{k}\frac{c^{k}u^{k}}{\Gamma\left(  \beta+k\right)
}=\frac{\Gamma\left(  r\right)  }{\Gamma\left(  p\right)  \Gamma\left(
\beta\right)  }\sum_{k=0}^{\infty}\frac{\Gamma\left(  p+k\right)  }%
{\Gamma\left(  r+k\right)  }\frac{\left(  cu\right)  ^{k}}{k!}=\frac{1}%
{\Gamma\left(  \beta\right)  }M\left(  p,r,cu\right)
\]
where $M\left(  a,b,z\right)  =\frac{\Gamma\left(  b\right)  }{\Gamma\left(
a\right)  }\sum_{k=0}^{\infty}\frac{\Gamma\left(  a+k\right)  }{\Gamma\left(
b+k\right)  }\frac{u^{k}}{k!}$ is Kummer's confluent hypergeometric function
\cite[13.1.2]{Abramowitz}. The Apelblat series (\ref{series_Apelblat}) thus
becomes%
\[
\frac{\Gamma\left(  r\right)  }{\Gamma\left(  p\right)  }\sum_{k=0}^{\infty
}\frac{\Gamma\left(  p+k\right)  \Gamma\left(  \beta+k\right)  }{\Gamma\left(
r+k\right)  k!}c^{k}t^{b+k-1}E_{a,b+k}\left(  xt^{a}\right)  =\int_{0}%
^{t}\left(  t-u\right)  ^{b-\beta-1}E_{a,b-\beta}\left(  x\left(  t-u\right)
^{a}\right)  u^{\beta-1}M\left(  p,r,cu\right)  du
\]
In order to use the property \cite[13.3.2]{Abramowitz} of the Kummer
function,
\[
\lim_{a\rightarrow\infty}M\left(  a,b,-\frac{z}{a}\right)  =\Gamma\left(
b\right)  z^{\frac{1-b}{2}}J_{b-1}\left(  2\sqrt{z}\right)
\]
we first choose $c=-\frac{1}{p}$ and then take the limit $p\rightarrow\infty$
of both sides becomes, with $\lim_{p\rightarrow\infty}\frac{\Gamma\left(
p+k\right)  }{\Gamma\left(  p\right)  p^{k}}=1$ (see \cite[6.1.47]%
{Abramowitz}),%
\[
\sum_{k=0}^{\infty}\frac{\Gamma\left(  \beta+k\right)  }{\Gamma\left(
r+k\right)  k!}\left(  -1\right)  ^{k}t^{b+k-1}E_{a,b+k}\left(  xt^{a}\right)
=\int_{0}^{t}\left(  t-u\right)  ^{b-\beta-1}E_{a,b-\beta}\left(  x\left(
t-u\right)  ^{a}\right)  u^{\beta+\frac{1-r}{2}-1}J_{r-1}\left(  2\sqrt
{u}\right)  du
\]
Let $\beta=r$, then $r=b-1$, we have
\[
\sum_{k=0}^{\infty}\left(  -1\right)  ^{k}\frac{t^{b+k-1}}{k!}E_{a,b+k}\left(
xt^{a}\right)  =\int_{0}^{t}E_{a}\left(  x\left(  t-u\right)  ^{a}\right)
u^{\frac{b-2}{2}}J_{b-2}\left(  2\sqrt{u}\right)  du
\]
which simplified for $b=2$ to%
\[
\sum_{k=1}^{\infty}\frac{\left(  -1\right)  ^{k-1}t^{k}}{\left(  k-1\right)
!}E_{a,k+1}\left(  xt^{a}\right)  =\int_{0}^{t}E_{a}\left(  xu^{a}\right)
J_{0}\left(  2\sqrt{t-u}\right)  du
\]

\medskip\refstepcounter{article}{\noindent\textbf{\thearticle. }%
}\ignorespaces\label{art_Berberan_Santos} \emph{Berberan-Santos' integral for
the double argument}. Berberan-Santos \cite{Berberan-SantosII_2005} applied a
simplified form (\ref{inverse_Laplace_transform_Berberan_Santos}) of the
inverse Laplace contour integral (\ref{inverse_Laplace}), which he deduced in
\cite{Berberan-SantosI_2005}, to the Mittag-Leffler function $E_{a,b}\left(
z\right)  $. First, let the Laplace transform of a real function $g\left(
u\right)  $ be equal to $E_{a,b}\left(  -z\right)  $, hence, $\int_{0}%
^{\infty}e^{-su}g\left(  u\right)  du=E_{a,b}\left(  -s\right)  $. From the
Laplace transform with $s=\sigma+iT$
\[
\left\vert E_{a,b}\left(  -s\right)  \right\vert =\left\vert \int_{0}^{\infty
}e^{-\sigma u}e^{-iTu}g\left(  u\right)  du\right\vert \leq\int_{0}^{\infty
}e^{-\sigma t}g\left(  t\right)  dt
\]
it follows that $\lim_{s\rightarrow\infty}\int_{0}^{\infty}e^{-st}g\left(
t\right)  dt=0$ for any direction in which $s$ with $\operatorname{Re}\left(
s\right)  >0$ tends to infinity. Since $E_{a,b}\left(  z\right)  $ is entire
function of order $\frac{1}{a}$, $E_{a,b}\left(  z\right)  =O\left(
e^{z^{\frac{1}{a}}}\right)  $, we have that $E_{a,b}\left(  -z\right)
=O\left(  e^{r^{\frac{1}{a}}e^{i\frac{\theta+\pi}{a}}}\right)  =O\left(
e^{r^{\frac{1}{a}}\cos\frac{\theta+\pi}{a}}\right)  $ and $\left\vert
E_{a,b}\left(  -s\right)  \right\vert \rightarrow0$ provided that $\cos
\frac{\theta+\pi}{a}<0$ or $\frac{\pi a}{2}-\pi<\theta<\frac{3\pi a}{2}-\pi$.
Since $s=re^{i\theta}$ with $-\frac{\pi}{2}\leq\theta\leq\frac{\pi}{2}$
because $\operatorname{Re}\left(  s\right)  >0$, the limit $\lim
_{s\rightarrow\infty}\left\vert E_{a,b}\left(  -s\right)  \right\vert =0$
requires that $0\leq a\leq1$. Incidentally, the definition of $c$ in Appendix
\ref{sec_inverse_Laplace_transform} indicates that $c=0$ and Berberan-Santos'
inverse Laplace transform (\ref{inverse_Laplace_transform_Berberan_Santos})
then yields%
\[
g(t)=\frac{2}{\pi}\int_{0}^{\infty}\operatorname{Re}\left(  E_{a,b}\left(
-iw\right)  \right)  \cos twdw\hspace{1cm}\text{for }t>0
\]
Next, Berberan-Santos observed from (\ref{Mittag_Leffler_arg_to_sqr(arg)})
that $E_{a,b}\left(  -iw\right)  =E_{2a,b}\left(  -w^{2}\right)
-iwE_{2a,b+2a}\left(  -w^{2}\right)  $ so that, for real $w$,%
\[
\operatorname{Re}\left(  E_{a,b}\left(  -iw\right)  \right)  =E_{2a,b}\left(
-w^{2}\right)
\]
Hence, the inverse Laplace transform becomes%
\[
g(t)=\frac{2}{\pi}\int_{0}^{\infty}E_{2a,b}\left(  -w^{2}\right)  \cos
twdw\hspace{1cm}\text{for }t>0
\]
Finally, taking the Laplace transform of both sides and reversing the
integrals
\[
E_{a,b}\left(  -s\right)  =\frac{2}{\pi}\int_{0}^{\infty}E_{2a,b}\left(
-w^{2}\right)  \left(  \int_{0}^{\infty}e^{-st}\cos twdt\right)  dw
\]
results in Berberan-Santos' remarkable integral for the double argument in $a$
(not $b$)%
\begin{equation}
E_{a,b}\left(  -s\right)  =\frac{2s}{\pi}\int_{0}^{\infty}\frac{E_{2a,b}%
\left(  -w^{2}\right)  }{s^{2}+w^{2}}dw\hspace{1cm}\text{for }0\leq
a\leq1\text{, }\operatorname{Re}\left(  s\right)  >0
\label{Berberan_Santos_integral_double_argument}%
\end{equation}

\textbf{Example} For $a=\frac{1}{2}$, $b=1$ and with $E_{1}\left(  z\right)
=e^{z}$, Berberan-Santos' integral
(\ref{Berberan_Santos_integral_double_argument}) yields again (\ref{E_1/2,1}),%
\[
E_{\frac{1}{2}}\left(  -x\right)  =\frac{2x}{\pi}\int_{0}^{\infty}%
\frac{e^{-t^{2}}}{t^{2}+x^{2}}dt=e^{x^{2}}\text{erfc}\left(  x\right)
\]

\medskip\refstepcounter{article}{\noindent\textbf{\thearticle. }%
}\ignorespaces\label{art_Berberan_Santos_eigen_bewijs}\emph{Another proof of
Berberan-Santos' integral (\ref{Berberan_Santos_integral_double_argument})}.
Application of Theorem \ref{theorem_entire_function_Cauchy} in Appendix
\ref{sec_complex_function_theorema} to $f\left(  z\right)  =E_{a,b}\left(
s-z\right)  $, only valid for $0\leq a\leq1$ because then $\lim_{r\rightarrow
\infty}\frac{E_{a,b}\left(  s-ir\right)  }{r^{2}}=0$ (see \textbf{art.}
\ref{art_Berberan_Santos}), yields%
\[
E_{a,b}\left(  -s\right)  =\frac{s}{\pi}\int_{-\infty}^{\infty}\frac
{E_{a,b}(-iw)}{s^{2}+w^{2}}dw
\]
Using (\ref{Mittag_Leffler_arg_to_sqr(arg)}), $E_{a,b}\left(  -iw\right)
=E_{2a,b}\left(  -w^{2}\right)  -iwE_{2a,b+2a}\left(  -w^{2}\right)  $ and the
fact that $wE_{2a,b+2a}\left(  -w^{2}\right)  $ is odd and $E_{2a,b}\left(
-w^{2}\right)  $ is even, again leads to
(\ref{Berberan_Santos_integral_double_argument}).

Using $\int_{0}^{\infty}e^{-t\left(  s^{2}+w^{2}\right)  }dt$ in
(\ref{Berberan_Santos_integral_double_argument}) and reversing the integrals,
justified by absolute convergence, gives%
\[
E_{a,b}\left(  -s\right)  =\frac{2s}{\pi}\int_{0}^{\infty}e^{-ts^{2}}\left\{
\int_{0}^{\infty}e^{-tw^{2}}E_{2a,b}\left(  -w^{2}\right)  dw\right\}
dt=\frac{s}{\pi}\int_{0}^{\infty}e^{-ts^{2}}\left\{  \int_{0}^{\infty}%
e^{-tx}x^{-\frac{1}{2}}E_{2a,b}\left(  -x\right)  dx\right\}  dt
\]
Thus, Berberan-Santos' integral
(\ref{Berberan_Santos_integral_double_argument}) is equivalent to the
statement that the Laplace transform of the Laplace transform of $x^{-\frac
{1}{2}}E_{2a,b}\left(  -x\right)  $ equals $\frac{\pi E_{a,b}\left(  -\sqrt
{s}\right)  }{\sqrt{s}}$. In fact, the latter property follows from the
Stieltjes transform, which is an iteration of the Laplace transform and which
is treated in the book by Widder \cite[Chapter VIII]{Widder}. Appendix
\ref{sec_inverse_Laplace_transform} discusses Gross' Laplace transform pair
\cite{Gross_1950}, in which the inverse Laplace transform
(\ref{inverse_Laplace}) is of the same form as the Laplace transform
(\ref{def_pgf_continuousrv}) itself.

\medskip\refstepcounter{article}{\noindent\textbf{\thearticle. }%
}\ignorespaces\label{art_Euler-Maclaurin_summation} \emph{Euler-Maclaurin
summation}. The Euler-Maclaurin summation formula \cite[p. 14]{Rademacher} is%
\begin{equation}
\sum_{n=\alpha+1}^{\beta}f(n)=\int_{\alpha}^{\beta}f(t)\,dt+\sum_{n=1}%
^{N}(-1)^{n}\frac{B_{n}}{n!}\left[  f^{(n-1)}(\beta)-f^{(n-1)}(\alpha)\right]
+R_{N} \label{sum_EulerMaclaurin}%
\end{equation}
with remainder term%
\[
R_{N}=\frac{(-1)^{N-1}}{N!}\int_{\alpha}^{\beta}B_{N}(u-\left[  u\right]
)\,f^{(N)}(u)\,du
\]
where $B_{n}$ and $B_{n}(x)$ are the Bernoulli numbers and the Bernoulli
polynomials defined in \cite[Chapter 23]{Abramowitz}.

The right-hand side summation in the Euler-Maclaurin summation formula
(\ref{sum_EulerMaclaurin}) requires higher order derivatives for the function
$f\left(  w\right)  =\frac{z^{w}}{\Gamma\left(  b+aw\right)  }$. We invoke
Leibniz' rule
\[
f^{(k)}(x)=\left.  \frac{d^{k}}{dw^{k}}\frac{z^{w}}{\Gamma\left(  b+aw\right)
}\right\vert _{w=x}=\sum_{j=0}^{k}\binom{k}{j}\left(  \log z\right)
^{k-j}z^{x}\left.  \frac{d^{j}}{dw^{j}}\frac{1}{\Gamma\left(  b+aw\right)
}\right\vert _{w=x}%
\]
which simplifies, with $\frac{d^{j}}{dw^{j}}\frac{1}{\Gamma\left(
b+aw\right)  }=a^{j}\left.  \frac{d^{j}}{dy^{j}}\frac{1}{\Gamma\left(
y\right)  }\right\vert _{y=b+ax}$, to%
\[
\frac{f^{(k)}(x)}{k!}=z^{x}\sum_{j=0}^{k}\frac{a^{j}}{j!}\left.  \frac{d^{j}%
}{dy^{j}}\frac{1}{\Gamma\left(  y\right)  }\right\vert _{y=b+ax}\frac{\left(
\log z\right)  ^{k-j}}{\left(  k-j\right)  !}%
\]
Since $\frac{1}{\Gamma\left(  z\right)  }$ is an entire function, the $j$-sum
converges when $k\rightarrow\infty$. Applying the Euler-Maclaurin summation
formula (\ref{sum_EulerMaclaurin}) to the definition
(\ref{Mittag_Leffler_function_Eab}) of $E_{a,b}\left(  z\right)  =\sum
_{k=0}^{\infty}\frac{z^{k}}{\Gamma\left(  b+ak\right)  }$ shows that
$\beta\rightarrow\infty$, but that we can choose the integer $\alpha=m$. The
asymptotic form (\ref{asymptotic_Gamma(b+rexp(it))}) indicates that
$\lim_{\beta\rightarrow\infty}\left.  \frac{d^{j}}{dw^{j}}\frac{1}%
{\Gamma\left(  b+aw\right)  }\right\vert _{w=\beta}=0$ for all $j$ and
(\ref{sum_EulerMaclaurin}) becomes
\begin{equation}
\sum_{k=m+1}^{\infty}\frac{z^{k}}{\Gamma\left(  b+ak\right)  }=\int
_{m}^{\infty}\frac{z^{t}}{\Gamma\left(  b+at\right)  }\,dt+z^{m}\sum
_{n=0}^{N-1}\frac{(-1)^{n}B_{n+1}}{n+1}\sum_{j=0}^{n}\frac{a^{j}}{j!}\left.
\frac{d^{j}}{dy^{j}}\frac{1}{\Gamma\left(  y\right)  }\right\vert
_{y=b+am}\frac{\left(  \log z\right)  ^{n-j}}{\left(  n-j\right)  !}+R_{N}
\label{Euler_Maclaurin_Ea,b_met_alfa}%
\end{equation}
The integral is written in terms of the integral $I_{a,b}\left(  z\right)  $,
defined in (\ref{def_integral_Ia,b(z)}), as%
\[
\int_{m}^{\infty}\frac{z^{t}}{\Gamma\left(  b+at\right)  }\,dt=z^{m}\int
_{0}^{\infty}\frac{z^{u}}{\Gamma\left(  b+am+au\right)  }\,du=z^{m}%
I_{a,b+ma}\left(  z\right)
\]

In summary, the Euler-Maclaurin expansion (\ref{Euler_Maclaurin_Ea,b_met_alfa}%
) is, with $\sum_{k=m+1}^{\infty}\frac{z^{k}}{\Gamma\left(  b+ak\right)
}=\sum_{k=0}^{\infty}\frac{z^{k+m+1}}{\Gamma\left(  b+a(m+1)+ak\right)  }$,%
\begin{equation}
z^{m+1}E_{a,b+(m+1)a}\left(  z\right)  =z^{m}I_{a,b+ma}\left(  z\right)
+z^{m}S_{N}+R_{N} \label{Euler_Maclaurin_kort}%
\end{equation}
where
\[
S_{N}=\sum_{n=0}^{N-1}(-1)^{n}\frac{B_{n+1}}{n+1}\sum_{j=0}^{n}\frac{a^{j}%
}{j!}\left.  \frac{d^{j}}{dy^{j}}\frac{1}{\Gamma\left(  y\right)  }\right\vert
_{y=b+ma}\frac{\left(  \log z\right)  ^{n-j}}{\left(  n-j\right)  !}%
\]

\medskip\refstepcounter{article}{\noindent\textbf{\thearticle. }%
}\ignorespaces\label{art_Euler-Maclaurin_remainder} \emph{Euler-Maclaurin sum
}$S_{N}$. We concentrate on $S_{N}$ and reverse the summations,%
\[
S_{N}=\sum_{j=0}^{N-1}\frac{a^{j}}{j!}\left.  \frac{d^{j}}{dy^{j}}\frac
{1}{\Gamma\left(  y\right)  }\right\vert _{y=b+ma}\sum_{n=j}^{N-1}%
(-1)^{n}\frac{B_{n+1}}{n+1}\frac{\left(  \log z\right)  ^{n-j}}{\left(
n-j\right)  !}%
\]
Let $x=\log z$ and use $\frac{d^{j}}{dx^{j}}x^{n}=\frac{n!}{\left(
n-j\right)  !}x^{n-j}$, then
\[
\sum_{n=j}^{N-1}(-1)^{n}\frac{B_{n+1}}{n+1}\frac{x^{n-j}}{\left(  n-j\right)
!}=\sum_{n=0}^{N-1}(-1)^{n}\frac{B_{n+1}}{n+1}\frac{1}{n!}\frac{d^{j}}{dx^{j}%
}x^{n}%
\]
and%
\begin{align*}
\lim_{N\rightarrow\infty}\sum_{n=j}^{N-1}(-1)^{n}\frac{B_{n+1}}{n+1}%
\frac{x^{n-j}}{\left(  n-j\right)  !}  &  =\frac{d^{j}}{dx^{j}}\sum
_{n=0}^{\infty}\frac{B_{n+1}}{\left(  n+1\right)  !}\left(  -x\right)  ^{n}\\
&  =\frac{d^{j}}{dx^{j}}\frac{1}{x}\sum_{n=1}^{\infty}\frac{B_{n}}{n!}\left(
-x\right)  ^{n}=\frac{d^{j}}{dx^{j}}\frac{1}{x}\left(  \sum_{n=0}^{\infty
}\frac{B_{n}}{n!}\left(  -x\right)  ^{n}-1\right)
\end{align*}
The generating function (\ref{gf_Bernoullinumbers}) of the Bernoulli numbers
$B_{n}$,%
\[
\frac{t\,}{e^{t}-1}=\sum_{n=0}^{\infty}B_{n}\,\frac{t^{n}}{n!}\hspace
{1cm}\text{convergent for }\left\vert t\right\vert <2\pi
\]
leads, for $\left\vert x\right\vert <2\pi$, to%
\[
\sum_{n=j}^{\infty}(-1)^{n}\frac{B_{n+1}}{n+1}\frac{x^{n-j}}{\left(
n-j\right)  !}=\frac{d^{j}}{dx^{j}}\frac{1}{x}\left(  \frac{-x\,}{e^{-x}%
-1}-1\right)  =\frac{d^{j}}{dx^{j}}\left(  \frac{1\,}{1-e^{-x}}-x^{-1}\right)
\]
For $N\rightarrow\infty$ and $\left\vert \log z\right\vert <2\pi$, we obtain
\begin{align*}
S_{\infty}  &  =\sum_{j=0}^{\infty}\frac{a^{j}}{j!}\left.  \frac{d^{j}}%
{dy^{j}}\frac{1}{\Gamma\left(  y\right)  }\right\vert _{y=b+ma}\sum
_{n=j}^{\infty}(-1)^{n}\frac{B_{n+1}}{n+1}\frac{\left(  \log z\right)  ^{n-j}%
}{\left(  n-j\right)  !}\\
&  =\sum_{j=0}^{\infty}\frac{a^{j}}{j!}\left.  \frac{d^{j}}{dy^{j}}\frac
{1}{\Gamma\left(  y\right)  }\right\vert _{y=b+ma}\left.  \frac{d^{j}}{dx^{j}%
}\left(  \frac{1\,}{1-e^{-x}}-x^{-1}\right)  \right\vert _{x=\log z}%
\end{align*}
Using the Fermi-Dirac integrals in \textbf{art}.
\ref{art_Taylor_series_to_Ea,b}, $F_{-j-1}(y)=\frac{d^{j}}{dy^{j}}\left(
\frac{1}{1+e^{-y}}\right)  $ and $\frac{d^{j}}{dx^{j}}x^{-1}=\left(
-1\right)  ^{j}j!x^{-1-j}$, we have%
\[
S_{\infty}=\sum_{j=0}^{\infty}\frac{a^{j}}{j!}\left.  \frac{d^{j}}{dy^{j}%
}\frac{1}{\Gamma\left(  y\right)  }\right\vert _{y=b+ma}F_{-j-1}(\log
z+i\pi)-\frac{1}{\log z}\sum_{j=0}^{\infty}\left.  \frac{d^{j}}{dy^{j}}%
\frac{1}{\Gamma\left(  y\right)  }\right\vert _{y=b+ma}\left(  -\frac{a}{\log
z}\right)  ^{j}%
\]
We recognize from (\ref{Ea,b(-z)_Taylor_in_a_around_b_in_plusa}) that, for
$\left\vert z\right\vert <1$,%
\[
E_{a,b}\left(  -ze^{i\pi}\right)  =E_{a,b}\left(  z\right)  =\frac{1}%
{\Gamma\left(  b\right)  }-\sum_{j=0}^{\infty}\frac{a^{j}}{j!}\left.
\frac{d^{j}}{dy^{j}}\frac{1}{\Gamma\left(  y\right)  }\right\vert
_{y=b}F_{-j-1}(\log z+i\pi)
\]
and from the series (\ref{I_a,b(z)_series}) for $I_{a,b}\left(  z\right)  $,
we obtain, for $\left\vert \frac{a}{\log z}\right\vert <1$,%
\begin{equation}
S_{\infty}=-E_{a,b+ma}\left(  z\right)  +\frac{1}{\Gamma\left(  b+ma\right)
}+I_{a,b+ma}\left(  z\right)  \label{Euler-Maclaurin_S_infy}%
\end{equation}
which is thus valid for $\left\vert z\right\vert <1$ and $\left\vert \frac
{a}{\log z}\right\vert <1$.

The Euler-Maclaurin expansion (\ref{Euler_Maclaurin_kort}) becomes, for
$N\rightarrow\infty$ and with (\ref{Euler-Maclaurin_S_infy}) and assuming that
$R_{N}\rightarrow0$,%
\[
z^{m+1}E_{a,b+(m+1)a}\left(  z\right)  +z^{m}E_{a,b+ma}\left(  z\right)
=2z^{m}I_{a,b+ma}\left(  z\right)  +\frac{z^{m}}{\Gamma\left(  b+ma\right)  }%
\]
With (\ref{E_a,ma}), we arrive, for $\left\vert z\right\vert <1$, $\left\vert
\frac{a}{\log z}\right\vert <1$ but all $b$, at%
\begin{equation}
E_{a,b}\left(  z\right)  =z^{m}I_{a,b+ma}\left(  z\right)  +\sum_{l=0}%
^{m}\frac{z^{l}}{\Gamma\left(  b+la\right)  }
\label{Euler_Maclaurin_Ea,b_met_vrije_m}%
\end{equation}

Numerical computations for $m\leq20$ show that
(\ref{Euler_Maclaurin_Ea,b_met_vrije_m}) is increasingly accurate for
increasing $m$ as long as $\left\vert \frac{a}{\log z}\right\vert <1$,
irrespective of $b$. On the other hand, when $\left\vert \frac{a}{\log
z}\right\vert >1$, increasing $m$ (up to 20) show a decreasing accuracy. Since
$\lim_{m\rightarrow\infty}z^{m}I_{a,b+ma}\left(  z\right)  =0$,
(\ref{Euler_Maclaurin_Ea,b_met_vrije_m}) reduces to an identity when
$m\rightarrow\infty$. Anticipating (\ref{I_a,b(z)_integral_Re(z^(1/a))>0}),
the Euler-Maclaurin sum (\ref{Euler_Maclaurin_Ea,b_met_alfa}) becomes%
\begin{align}
\sum_{k=m+1}^{\infty}\frac{z^{k}}{\Gamma\left(  b+ak\right)  }  &
=\frac{z^{\frac{1-b}{a}}}{a}\left\{  e^{z^{\frac{1}{a}}}+\int_{0}^{\infty
}\frac{e^{-z^{\frac{1}{a}}x}}{x^{b+ma}}\left(  \frac{\frac{\sin\left(
b+ma\right)  \pi}{\pi}\ln x+\cos\left(  b+ma\right)  \pi}{\pi^{2}+\left(  \ln
x\right)  ^{2}}\right)  dx\right\} \nonumber\\
&  +z^{m}\sum_{n=0}^{N-1}(-1)^{n}\frac{B_{n+1}}{n+1}\sum_{j=0}^{n}\frac{a^{j}%
}{j!}\left.  \frac{d^{j}}{dy^{j}}\frac{1}{\Gamma\left(  y\right)  }\right\vert
_{y=b+am}\frac{\left(  \log z\right)  ^{n-j}}{\left(  n-j\right)  !}+R_{N}
\label{Euler_MacLaurin_expansion_Ea,b_met_m}%
\end{align}
Numerical computations (for $m=-1$) indicate that the $n$-summation converges
for small $a$ and small $z$, with fast convergence when $z=1$. For $z=1$,
(\ref{Euler_MacLaurin_expansion_Ea,b_met_m}) simplifies to%
\begin{align*}
\sum_{k=m+1}^{\infty}\frac{1}{\Gamma\left(  b+ak\right)  }  &  =\frac{1}%
{a}\left\{  e+\int_{0}^{\infty}\frac{e^{-x}}{x^{b+ma}}\left(  \frac{\frac
{\sin\left(  b+ma\right)  \pi}{\pi}\ln x+\cos\left(  b+ma\right)  \pi}{\pi
^{2}+\left(  \ln x\right)  ^{2}}\right)  dx\right\} \\
&  +\sum_{n=0}^{\infty}\frac{B_{n+1}(-a)^{n}}{\left(  n+1\right)  !}\left.
\frac{d^{n}}{dy^{n}}\frac{1}{\Gamma\left(  y\right)  }\right\vert _{y=b+am}%
\end{align*}
where the latter sum converges for $\left\vert a\right\vert <1$.

\section{Complex integrals for $E_{a,b}\left(  z\right)  $}

\medskip\refstepcounter{article}{\noindent\textbf{\thearticle. }%
}\ignorespaces\label{art_THE_inverse_Laplace_transform} \emph{Basic complex
integral for }$E_{a,b}\left(  z\right)  $. Inverse Laplace transformation
(\ref{inverse_Laplace}) of (\ref{Laplace_transform_positive_argument}) yields,
with $c^{\prime}>\left\vert x\right\vert ^{\frac{1}{a}}$ for
$\operatorname{Re}\left(  b\right)  >0$,%
\[
t^{b-1}E_{a,b}\left(  xt^{a}\right)  =\frac{1}{2\pi i}\int_{c^{\prime}%
-i\infty}^{c^{\prime}+i\infty}\frac{s^{a-b}e^{st}}{s^{a}-x}ds
\]
where the integrand is analytic for $\operatorname{Re}\left(  s\right)
>c^{\prime}$. For real $t\geq0$, we can move the line of integration to
$c>t^{a}\left\vert x\right\vert =t^{a}\left(  c^{\prime}\right)  ^{a}$ and
perform an ordinary substitution $w=st$. Let $z=xt^{a}$ and we arrive at the
\emph{basic complex integral}%
\begin{equation}
E_{a,b}\left(  z\right)  =\frac{1}{2\pi i}\int_{c-i\infty}^{c+i\infty}%
\frac{w^{a-b}e^{w}}{w^{a}-z}dw\hspace{1cm}c>\left\vert z\right\vert
\label{E_a,b(z)_basic_contour_integral}%
\end{equation}

The basic complex integral (\ref{E_a,b(z)_basic_contour_integral}) of
$E_{a,b}\left(  z\right)  $ can also be deduced from Hankel's deformed
integral (\ref{Hankel_contour_phi_integral_1opGamma}) in the power series of
the Mittag-Leffler function (\ref{Mittag_Leffler_function_Eab})
\[
E_{a,b}\left(  z\right)  =\sum_{k=0}^{\infty}\frac{z^{k}}{\Gamma\left(
b+ak\right)  }=\frac{1}{2\pi i}\sum_{k=0}^{\infty}\int_{C_{\phi}}\left(
zw^{-a}\right)  ^{k}w^{-b}e^{w}dw=\frac{1}{2\pi i}\int_{C_{\phi}}\sum
_{k=0}^{\infty}\left(  zw^{-a}\right)  ^{k}w^{-b}e^{w}dw
\]
Only if $\left\vert zw^{-a}\right\vert <1$, implying $\left\vert z\right\vert
<\left\vert w^{a}\right\vert $ for any $w$ along the contour $C_{\phi}$
described in \textbf{art}. \ref{art_Hankel_integral}, then we obtain
\begin{equation}
E_{a,b}\left(  z\right)  =\frac{1}{2\pi i}\int_{C_{\phi}}\frac{w^{-b}e^{w}%
}{1-zw^{-a}}dw \label{E_a,b(z)_basic_contour_integral_Hankel}%
\end{equation}
where the radius $\rho$ of the circle at $w=0$ in the contour $C_{\phi}$ (see
\textbf{art}. \ref{art_Hankel_integral}) can be appropriately chosen to
satisfy $\left\vert z\right\vert <\rho^{\operatorname{Re}a}$ or $\left\vert
z\right\vert ^{\frac{1}{\operatorname{Re}a}}<\rho$ for any $z$. We obtain
again (\ref{E_a,b(z)_basic_contour_integral}) by choosing $\phi=\frac{\pi}{2}$
and $c=\rho$.

\medskip\refstepcounter{article}{\noindent\textbf{\thearticle. }%
}\ignorespaces\label{art_Ea,b_negative_a} \emph{Mittag-Leffler function
}$E_{a,b}\left(  z\right)  $\emph{ for negative real }$a$\emph{. }Although
$E_{a,b}\left(  z\right)  $, defined by the power series
(\ref{Mittag_Leffler_function_Eab}), is not valid for negative
$\operatorname{Re}\left(  a\right)  $ values, the Taylor series
(\ref{Ea,b(-z)_Taylor_in_a_around_b_in_plusa}) of $E_{a,b}\left(  z\right)  $
in $a$ around $a_{0}=0$,%
\[
E_{a,b}\left(  -z\right)  =\frac{1}{\Gamma\left(  b\right)  }-\sum
_{m=0}^{\infty}\frac{1}{m!}\left.  \frac{d^{m}}{du^{m}}\frac{1}{\Gamma\left(
u\right)  }\right\vert _{u=b}F_{-m-1}(\log z)a^{m}%
\]
and its companion (\ref{Ea,b(-z)_Taylor_in_a_around_b_in_mina})
\[
E_{a,b}\left(  -z\right)  =\sum_{m=0}^{\infty}\frac{1}{m!}\left.  \frac{d^{m}%
}{du^{m}}\frac{1}{\Gamma\left(  u\right)  }\right\vert _{u=b}F_{-m-1}(-\log
z)\left(  -a\right)  ^{m}%
\]
indicate that%
\[
E_{a,b}\left(  -\frac{1}{z}\right)  =\frac{1}{\Gamma\left(  b\right)
}-E_{-a,b}\left(  -z\right)
\]
which may be used, by analytic continuation, as a definition of $E_{a,b}%
\left(  z\right)  $ for $\operatorname{Re}\left(  a\right)  <0$. Thus, for
$\operatorname{Re}\left(  a\right)  <0$, the Mittag-Leffler function
$E_{a,b}\left(  z\right)  $ possesses an essential singularity at $z=0$ and is
not an entire function.

The same result is also deduced in \cite[Sec. 4.8, pp. 80-82]{Gorenflo_2014}
from the basic complex integral (\ref{E_a,b(z)_basic_contour_integral_Hankel})
in \textbf{art}. \ref{art_THE_inverse_Laplace_transform}
\[
E_{a,b}\left(  z\right)  =\frac{1}{2\pi i}\int_{C_{\phi}}\frac{w^{a-b}e^{w}%
}{w^{a}-z}dw
\]
is valid for any complex $a,b$ and $z$. Introducing the identity%
\[
\frac{w^{a-b}}{w^{a}-z}=\frac{1}{w^{b}-zw^{b-a}}=\frac{1}{w^{b}}-\frac
{1}{w^{b}-\frac{1}{z}w^{a+b}}%
\]
leads to%
\[
E_{a,b}\left(  z\right)  =\frac{1}{2\pi i}\int_{C_{\phi}}\frac{e^{w}}{w^{b}%
}dw-\frac{1}{2\pi i}\int_{C_{\phi}}\frac{w^{-a-b}e^{w}}{w^{-a}-\frac{1}{z}}dw
\]
Invoking Hankel's integral (\ref{Hankel_contour_long_line}) to the first
integral and the basic complex integral
(\ref{E_a,b(z)_basic_contour_integral_Hankel}) to the last integral leads
again to%
\begin{equation}
E_{a,b}\left(  z\right)  =\frac{1}{\Gamma\left(  b\right)  }-E_{-a,b}\left(
\frac{1}{z}\right)  \label{Mittag_Leffler_negative_a}%
\end{equation}
The Taylor series (\ref{Mittag_Leffler_function_Eab}) of $E_{a,b}\left(
z\right)  $ then shows that%
\[
E_{-a,b}\left(  z\right)  =\frac{1}{\Gamma\left(  b\right)  }-E_{a,b}\left(
\frac{1}{z}\right)  =-\sum_{k=1}^{\infty}\frac{z^{-k}}{\Gamma\left(
b+ak\right)  }%
\]

\medskip\refstepcounter{article}{\noindent\textbf{\thearticle. }%
}\ignorespaces\label{art_evaluation_basic_complex_integral} \emph{Evaluation
of the basic complex integral for }$E_{a,b}\left(  z\right)  $. We will now
evaluate the integral (\ref{E_a,b(z)_basic_contour_integral}) by closing the
contour over the negative $\operatorname{Re}\left(  w\right)  $-plane.
Consider the closed loop integral
\[
\frac{1}{2\pi i}\int_{L}\frac{w^{a-b}e^{w}}{w^{a}-z}dw
\]
where the contour $L$ consists of the vertical line at $w=c+it$, the circle
segment at infinity turning from $\frac{\pi}{2}$ towards $\pi-\varepsilon$,
the line above the negative real axis, the small circle with radius
$\varepsilon$ turning around the origin from $\pi-\varepsilon$ to
$-\pi+\varepsilon$, the line from the origin just below the negative real
axis, over the circle segment with infinite radius from $-\pi+\varepsilon$ to
$-\frac{\pi}{2}$ and ending at the negative side of the vertical line. The
integrand has only singularities of $\left(  w^{a}-z\right)  ^{-1}$, which are
simple poles because the zeros of $w^{a}-z=0$ are all simple and lie at
$w_{k}=\left\vert z\right\vert ^{\frac{1}{a}}e^{i\frac{\arg z}{a}}%
e^{i\frac{2\pi k}{a}}$ for each $k\in\mathbb{Z}$. If $a$ is irrational, there
are infinitely many zeros. The contour $L$ only encloses the poles with
argument $\theta_{k}=\frac{\arg z}{a}+\frac{2\pi k}{a}$ between $-\pi
\leq\theta_{k}\leq\pi$. Thus, the integer $k$ ranges from%
\[
-\left\lfloor \frac{a}{2}+\frac{\arg z}{2\pi}\right\rfloor \leq k\leq
\left\lfloor \frac{a}{2}-\frac{\arg z}{2\pi}\right\rfloor
\]
and there are precisely $\left\lfloor a\right\rfloor $ enclosed poles, where
$\left\lfloor v\right\rfloor $ is the integer smaller than or equal to $v$.
Hence, $a$ must be at least equal to 1, else no singularities are enclosed. By
Cauchy's residue theorem, we obtain%
\begin{align*}
\frac{1}{2\pi i}\int_{L}\frac{w^{a-b}e^{w}}{w^{a}-z}dw  &  =\sum_{k\in\left[
-\left\lfloor \frac{a}{2}+\frac{\arg z}{2\pi}\right\rfloor ,\left\lfloor
\frac{a}{2}-\frac{\arg z}{2\pi}\right\rfloor \right]  }\lim_{w\rightarrow
w_{k}}\frac{w-w_{k}}{w^{a}-z}w^{a-b}e^{w}\\
&  =\frac{1}{a}\sum_{k\in\left[  -\left\lfloor \frac{a}{2}+\frac{\arg z}{2\pi
}\right\rfloor ,\left\lfloor \frac{a}{2}-\frac{\arg z}{2\pi}\right\rfloor
\right]  }w_{k}^{1-b}e^{w_{k}}\hspace{1cm}\text{with }w_{k}=\left\vert
z\right\vert ^{\frac{1}{a}}e^{i\frac{\arg z}{a}}e^{i\frac{2\pi k}{a}}%
\end{align*}
We can always choose $0<\varepsilon<r$ small enough, provided that
$\operatorname{Re}\left(  a-b\right)  >-1$. Evaluation of the contour $L$
yields%
\begin{align*}
\frac{1}{2\pi i}\int_{L}\frac{w^{a-b}e^{w}}{w^{a}-z}dw  &  =\frac{1}{2\pi
i}\int_{c-i\infty}^{c+i\infty}\frac{w^{a-b}e^{w}}{w^{a}-z}dw+\frac{1}{2\pi
i}\left\{  \int_{\infty}^{0}\frac{\left(  ye^{i\pi}\right)  ^{a-b}e^{-y}%
}{\left(  ye^{i\pi}\right)  ^{a}-z}dy+\int_{0}^{\infty}\frac{\left(
ye^{-i\pi}\right)  ^{a-b}e^{-y}}{\left(  ye^{-i\pi}\right)  ^{a}-z}dy+\right\}
\\
&  =E_{a,b}\left(  z\right)  +\frac{1}{2\pi i}\left\{  \int_{0}^{\infty
}\left(  \frac{e^{-i\pi\left(  a-b\right)  }}{y^{a}e^{-i\pi a}-z}%
dy-\frac{e^{i\pi\left(  a-b\right)  }}{y^{a}e^{i\pi a}-z}\right)
y^{a-b}e^{-y}dy\right\}
\end{align*}
Hence, we have%
\[
E_{a,b}\left(  z\right)  =-\frac{1}{2\pi i}\left\{  \int_{0}^{\infty}\left(
\frac{e^{-\left(  a-b\right)  i\pi}}{e^{-ai\pi}y^{a}-z}-\frac{e^{\left(
a-b\right)  i\pi}}{e^{ai\pi}y^{a}-z}\right)  y^{a-b}e^{-y}dy\right\}
+\frac{1}{a}\sum_{k\in\left[  -\left\lfloor \frac{a}{2}+\frac{\arg z}{2\pi
}\right\rfloor ,\left\lfloor \frac{a}{2}-\frac{\arg z}{2\pi}\right\rfloor
\right]  }w_{k}^{1-b}e^{w_{k}}%
\]
Finally for $\operatorname{Re}\left(  a\right)  >\operatorname{Re}\left(
b\right)  -1$, we arrive at%
\begin{align}
E_{a,b}\left(  z\right)   &  =z\frac{\sin\left(  a-b\right)  \pi}{\pi}\int
_{0}^{\infty}\frac{y^{a-b}e^{-y}dy}{y^{2a}-2zy^{a}\cos a\pi+z^{2}}+\frac
{\sin\pi b}{\pi}\int_{0}^{\infty}\frac{y^{2a-b}e^{-y}}{y^{2a}-2zy^{a}\cos
a\pi+z^{2}}dy\nonumber\\
&  +\frac{1}{a}\sum_{k\in\left[  -\left\lfloor \frac{a}{2}+\frac{\arg z}{2\pi
}\right\rfloor ,\left\lfloor \frac{a}{2}-\frac{\arg z}{2\pi}\right\rfloor
\right]  }w_{k}^{1-b}e^{w_{k}}\hspace{1cm}\text{with }w_{k}=\left\vert
z\right\vert ^{\frac{1}{a}}e^{i\frac{\arg z}{a}}e^{i\frac{2\pi k}{a}}
\label{Mittag_Leffler_function_as_residues}%
\end{align}
The last residue sum illustrates again that $E_{a,b}\left(  z\right)  $ is an
entire function of order $\frac{1}{a}$. The sign of $w_{k}$ is determined by
$\cos\left(  \frac{\arg z+2\pi k}{a}\right)  $. Consequently for large
$\left\vert z\right\vert $, in the sectors $-\frac{\pi a}{2}-2\pi k<\arg
z<\frac{\pi a}{2}-2\pi k$, the function $E_{a,b}\left(  z\right)  $ tends to
infinity, while in the complementary sectors $\frac{\pi a}{2}-2\pi k<\arg
z<\frac{3\pi a}{2}-2\pi k$, $E_{a,b}\left(  z\right)  \rightarrow0$.

For real $z=-x^{a}=x^{a}e^{i\pi}$ in
(\ref{Mittag_Leffler_function_as_residues}) and $w_{k}=xe^{i\frac{\left(
2k+1\right)  \pi}{a}}$, we find%
\begin{align*}
E_{a,b}\left(  -x^{a}\right)   &  =-x^{a}\frac{\sin\left(  a-b\right)  \pi
}{\pi}\int_{0}^{\infty}\frac{y^{a-b}e^{-y}dy}{y^{2a}+2x^{a}y^{a}\cos
a\pi+x^{2a}}+\frac{\sin\pi b}{\pi}\int_{0}^{\infty}\frac{y^{2a-b}e^{-y}%
}{y^{2a}+2x^{a}y^{a}\cos a\pi+z^{2}}dy\\
&  +\frac{x^{1-b}}{a}e^{i\left\{  \frac{\pi\left(  1-b\right)  }{a}\right\}
}\sum_{k=0}^{\left\lfloor a\right\rfloor }e^{x\cos\frac{\left(  2k-1\right)
\pi}{a}}e^{i\left\{  -x\sin\frac{\left(  2k-1\right)  \pi}{a}-\frac
{2k\pi\left(  1-b\right)  }{a}\right\}  }%
\end{align*}
For real $x$, $b=1$ and $0<a<1$, the condition to enclose a pole is $-a\pi
\leq\arg z+2\pi k\leq a\pi$. Since $\arg z=\pi\notin\left[  -a\pi,a\pi\right]
$ for $0<a<1$, the above sum disappears and we obtain%
\begin{equation}
E_{a}\left(  -x^{a}\right)  =x^{a}\frac{\sin a\pi}{\pi}\int_{0}^{\infty}%
\frac{y^{a-1}e^{-y}dy}{y^{2a}+2x^{a}y^{a}\cos a\pi+x^{2a}}\hspace
{1cm}\text{for }0<a<1 \label{Mittag_Leffler_E_a(-x^a)_0<a<1_integral}%
\end{equation}
from which the asymptotic behavior for large $x$, with $y^{2a}+2x^{a}y^{a}\cos
a\pi+x^{2a}\sim x^{2a}$, are
\[
E_{a}\left(  -x^{a}\right)  \sim\frac{\sin a\pi}{\pi}\frac{\Gamma\left(
a\right)  }{x^{a}}%
\]
For small, real $x$, the series (\ref{Mittag_Leffler_function_Eab}) gives
\[
E_{a}\left(  -x^{a}\right)  =\sum_{k=0}^{\infty}\frac{\left(  -1\right)
^{k}x^{ak}}{\Gamma\left(  1+ak\right)  }=1-\frac{x^{a}}{\Gamma\left(
1+a\right)  }+O\left(  x^{2a}\right)  =\exp\left(  -\frac{x^{a}}{\Gamma\left(
1+a\right)  }\right)  +O\left(  x^{2a}\right)
\]
Hence, $E_{a}\left(  -x^{a}\right)  $ for $0<a<1$ is said to \textquotedblleft
interpolate\textquotedblright\ between an exponential (for small $x$) and a
power law (for large $x$). These two regimes have been studied by Mainardi
\cite{Mainardi_2014}, who illustrated their accuracy with several plots for
$a=0.25,0.5,0.75,0.9$ and $0.99$. \textbf{Art}.
\ref{art_1_op_sinpiw_integral_for_0<a<1} presents another derivation in
(\ref{E_a(-x)_integral_a_in(0,1]}) that is equivalent to
(\ref{Mittag_Leffler_E_a(-x^a)_0<a<1_integral}).

\medskip\refstepcounter{article}{\noindent\textbf{\thearticle. }%
}\ignorespaces\label{art_Mittag_Leffler_contour_integral}
\emph{Mittag-Leffler's contour integral}. In a rather long article of 1905,
Mittag-Leffler \cite[at p. 133-135]{Mittag_Leffler_1905} proceeds one step
further and substitutes $w=t^{\frac{1}{a}}$ in the basic complex integral
(\ref{E_a,b(z)_basic_contour_integral}). The map $w\rightarrow t^{\frac{1}{a}%
}$ is multi-valued. For $w=\left\vert w\right\vert e^{i\theta_{w}}$ with
$-\pi\leq\theta_{w}\leq\pi$, the inverse map $t\rightarrow w^{a}$ shows that
$t=\left\vert t\right\vert e^{i\theta_{t}}=\left\vert w\right\vert
^{a}e^{ia\theta_{w}+2\pi iak}$, from which the argument $\theta_{t}%
=a\theta_{w}+2\pi ak$ for any integer $k$. The contour $C_{\phi}$ in the
$w$-plane requires that $\frac{\pi}{2}<\left\vert \phi\right\vert <\pi$,
because $e^{w}\rightarrow0$ for large $w$ at the contour $C_{\phi}$.
Similarly, the transformed contour $C_{a\phi}$ in the $t$-plane requires that
$\frac{a\pi}{2}<\left\vert \arg t\right\vert <a\pi$ in order that
$e^{t^{\frac{1}{a}}}\rightarrow0$ along the straight lines towards infinity.
Thus, the map $w\rightarrow t^{\frac{1}{a}}$ changes the angles from $\phi$ to
$a\phi$ of the straight lines in the contour (Fig. 43 in Bieberbach's book
\cite[p. 273]{Bieberbach}). Moreover, we must choose the branch (i.e. the
appropriate integer $k$ in $\theta_{t}=a\theta_{w}+2\pi ak$) of the function
$t^{\frac{1}{a}}$ that is positive for positive $t$, i.e. along the positive
real $t$-axis, because the same holds for along the positive real $w$-axis.
Performing the substitution $w=t^{\frac{1}{a}}$ in
(\ref{E_a,b(z)_basic_contour_integral}) leads, for $\left\vert z\right\vert
<\left\vert t\right\vert $ and $\left\vert \arg z\right\vert <\frac{a\pi}{2}$,
to Mittag-Leffler's contour integral%
\begin{equation}
E_{a,b}\left(  z\right)  =\frac{1}{2\pi ia}\int_{C_{a\phi}}\frac{t^{\frac
{1-b}{a}}e^{t^{\frac{1}{a}}}}{t-z}dt\hspace{1cm}\text{with }\frac{\pi}%
{2}<\left\vert \phi\right\vert <\pi\label{Mittag_Leffler_Bieberbach_integral}%
\end{equation}
Mittag-Leffler's integral (\ref{Mittag_Leffler_Bieberbach_integral}) for $b=1$
is actually more elegant than the basic complex integral
(\ref{E_a,b(z)_basic_contour_integral}) at the expense of a more complicated
contour $C_{a\phi}$.

\medskip\refstepcounter{article}{\noindent\textbf{\thearticle. }%
}\ignorespaces\label{art_deductions_Mittag_Leffler_contour_integral}
\emph{Deductions from Mittag-Leffler's contour integral
(\ref{Mittag_Leffler_Bieberbach_integral}) for }$E_{a,b}\left(  z\right)  $.
Based on (\ref{Mittag_Leffler_Bieberbach_integral}), we follow Bieberbach
\cite[p. 273]{Bieberbach}, who deduces two bounds for $0<a<2$ and $b=1$.
First, using%
\[
\frac{1}{t-z}=-\frac{1}{z}+\frac{t}{z\left(  t-z\right)  }%
\]
in (\ref{Mittag_Leffler_Bieberbach_integral}), we have%
\begin{align*}
E_{a,b}\left(  z\right)   &  =-\frac{1}{z}\frac{1}{2\pi ia}\int_{C_{a\phi}%
}t^{\frac{1-b}{a}}e^{t^{\frac{1}{a}}}dt+\frac{1}{2\pi ia}\int_{C_{a\phi}}%
\frac{t^{\frac{1-b+a}{a}}e^{t^{\frac{1}{a}}}}{z\left(  t-z\right)  }dt\\
&  =-\frac{1}{z}\frac{1}{2\pi i}\int_{C_{\phi}}w^{a-b}e^{w}dw+\frac{1}{2\pi
ia}\int_{C_{a\phi}}\frac{t^{\frac{a-b-1}{a}}e^{t^{\frac{1}{a}}}}{z\left(
t-z\right)  }dt
\end{align*}
With Hankel's contour integral (\ref{Hankel_contour_phi_integral_1opGamma}),%
\[
E_{a,b}\left(  z\right)  =-\frac{1}{z}\frac{1}{\Gamma\left(  b-a\right)
}+\frac{1}{2\pi ia}\int_{C_{a\phi}}\frac{t^{\frac{a-b+1}{a}}e^{t^{\frac{1}{a}%
}}}{z\left(  t-z\right)  }dt
\]
The remaining integral is upper bounded by%
\[
\left\vert \frac{1}{2\pi ia}\int_{C_{a\phi}}\frac{t^{\frac{a-b+1}{a}%
}e^{t^{\frac{1}{a}}}}{z\left(  t-z\right)  }dt\right\vert <\frac{1}{2\pi
a\left\vert z\right\vert }\int_{C_{a\phi}}\frac{\left\vert t\right\vert
^{\frac{a-b+1}{a}}\left\vert e^{t^{\frac{1}{a}}}\right\vert }{\left\vert
z\right\vert \left\vert 1-\frac{t}{z}\right\vert }\left\vert dt\right\vert
=\frac{c}{\left\vert z\right\vert ^{2}}%
\]
because $\frac{1}{2\pi a}\int_{C_{a\phi}}\frac{\left\vert t\right\vert
^{\frac{a-b-1}{a}}\left\vert e^{t^{\frac{1}{a}}}\right\vert }{\left\vert
1-\frac{t}{z}\right\vert }\left\vert dt\right\vert $ converges for
$t\rightarrow\infty$ if $\left\vert e^{t^{\frac{1}{a}}}\right\vert
=e^{\left\vert t\right\vert ^{\frac{1}{a}}\cos\frac{\arg t}{a}}\rightarrow0$
for large $t$, which requires that $\frac{\arg t}{a}\geq\frac{\pi}{2}$. In
addition, we must prevent that $\frac{t}{z}$ for large $z$ can tend
arbitrarily close to 1, which is guaranteed if $\left\vert \arg z\right\vert
\notin\left(  \frac{a\pi}{2},a\pi\right)  $, because $\left\vert \arg
t\right\vert \in\left(  \frac{a\pi}{2},a\pi\right)  $. Hence, for $\left\vert
\arg z\right\vert \notin\left(  \frac{a\pi}{2},a\pi\right)  $, we arrive at%
\[
\left\vert E_{a,b}\left(  z\right)  +\frac{1}{z\Gamma\left(  b-a\right)
}\right\vert <\frac{c}{\left\vert z\right\vert ^{2}}%
\]

For the second bound, Bieberbach \cite[p. 275]{Bieberbach} cleverly observes
that a similar integration path $C_{a\phi}^{\prime}$ as in Mittag-Leffler's
integral (\ref{Mittag_Leffler_Bieberbach_integral}) can be followed with the
only difference that the circular part of the path now has a radius smaller
than $\left\vert z\right\vert $. In other words, while $\left\vert
z\right\vert <\left\vert t\right\vert $ in Mittag-Leffler's integral
(\ref{Mittag_Leffler_Bieberbach_integral}), the path $C_{a\phi}^{\prime}$ now
turns over an angle $-a\phi$ to $a\phi$ with the radius smaller than
$\left\vert z\right\vert $. The closed contour (see also \cite[p. 346, Fig.
6.13-2]{Sansone}), that first follows the Mittag-Leffler path $C_{a\phi}$ and
returns via the path $C_{a\phi}^{\prime}$, encloses the point $t=z$ as the
only singularity, provided $-a\phi<\arg z<a\phi$. Hence, by Cauchy's residue
theorem, it holds that%
\[
\frac{1}{2\pi ia}\int_{C_{a\phi}}\frac{t^{\frac{1-b}{a}}e^{t^{\frac{1}{a}}}%
}{t-z}dt-\frac{1}{2\pi ia}\int_{C_{a\phi}^{\prime}}\frac{t^{\frac{1-b}{a}%
}e^{t^{\frac{1}{a}}}}{t-z}dt=\frac{1}{a}z^{\frac{1-b}{a}}e^{z^{\frac{1}{a}}}%
\]
and%
\[
E_{a,b}\left(  z\right)  =\frac{1}{2\pi ia}\int_{C_{a\phi}^{\prime}}%
\frac{t^{\frac{1-b}{a}}e^{t^{\frac{1}{a}}}}{t-z}dt+\frac{1}{a}z^{\frac{1-b}%
{a}}e^{z^{\frac{1}{a}}}%
\]
from which%
\[
\left\vert E_{a,b}\left(  z\right)  -\frac{1}{a}z^{\frac{1-b}{a}}%
e^{z^{\frac{1}{a}}}\right\vert \leq\frac{1}{2\pi a\left\vert z\right\vert
}\int_{C_{a\phi}^{\prime}}\frac{\left\vert t^{\frac{1-b}{a}}\right\vert
\left\vert e^{t^{\frac{1}{a}}}\right\vert }{\left\vert \frac{t}{z}%
-1\right\vert }dt=\frac{c^{\prime}}{\left\vert z\right\vert }%
\]
by the same convergence argument as above, where now on the circular segment
$\left\vert t\right\vert <\left\vert z\right\vert $ can be chosen small
enough. In summary, the second bound for $\left\vert \arg z\right\vert
\in\left(  \frac{a\pi}{2},a\pi\right)  $ is
\[
\left\vert E_{a,b}\left(  z\right)  -\frac{1}{a}z^{\frac{1-b}{a}}%
e^{z^{\frac{1}{a}}}\right\vert <\frac{c^{\prime}}{\left\vert z\right\vert }%
\]
The derivation illustrates why Bieberbach considers $0<a<2$, because
$\left\vert \arg z\right\vert \in\left(  0,2\pi\right)  $.

The second bound shows that $aE_{a,b}\left(  z^{a}\right)  \approx
z^{1-b}e^{z}$ is independent of $a$ so that $aE_{a,b}\left(  z^{a}\right)
\approx\frac{1}{a}E_{\frac{1}{a},b}\left(  z^{\frac{1}{a}}\right)  $. Hence,
we are led for non-negative real $z$ and $\left\vert z\right\vert
>\varepsilon$ to%
\[
E_{\frac{1}{a},b}\left(  z\right)  \approx a^{2}E_{a,b}\left(  z^{a^{2}%
}\right)
\]
whose exact corresponding relation (\ref{Integral_Ia,b_map_1/a_to_a}) for the
associated integral $I_{a,b}\left(  z\right)  =\int_{0}^{\infty}\frac{z^{u}%
}{\Gamma\left(  b+au\right)  }\,du$, explored in Section
\ref{sec_integral_Iab(z)}, is $I_{\frac{1}{a},b}\left(  z\right)
=a^{2}I_{a,b}\left(  z^{a^{2}}\right)  $.

\medskip\refstepcounter{article}{\noindent\textbf{\thearticle. }%
}\ignorespaces\label{art_complex_integral_cosec_type} \emph{Deductions from
Cauchy's residue theorem}. If $f\left(  z\right)  $ is an entire function and
$\lim_{r\rightarrow\infty}\left\vert \frac{f\left(  re^{i\theta}\right)
}{\sin\pi re^{i\theta}}\right\vert \rightarrow0$ in a semicircle with either
$-\frac{\pi}{2}<\theta<\frac{\pi}{2}$ or $\frac{\pi}{2}<\theta<\frac{3\pi}{2}%
$, then it follows directly from Cauchy's residue theorem
\cite{Titchmarshfunctions} that%
\[
\frac{1}{2\pi i}\int_{c-i\infty}^{c+i\infty}\frac{\pi}{\sin\pi w}f\left(
w\right)  dw=-1_{\left\{  -\frac{\pi}{2}<\theta<\frac{\pi}{2}\right\}  }%
\sum_{k=0}^{\infty}\left(  -1\right)  ^{k}f\left(  k\right)  \hspace
{1cm}\text{for }0<c<1
\]
The definition (\ref{Mittag_Leffler_function_Eab}) of $E_{a,b}\left(
z\right)  =\sum_{k=0}^{\infty}\frac{z^{k}}{\Gamma\left(  b+ak\right)  }$
contains $f\left(  w\right)  =\frac{z^{w}}{\Gamma\left(  b+aw\right)  }$,
which is an entire function in $w$. From the asymptotic behavior
(\ref{asymptotic_Gamma(b+rexp(it))}) of $\frac{1}{\left\vert \Gamma\left(
b+are^{i\theta}\right)  \right\vert }$ in \textbf{art}.
\ref{art_asymptotic_behavior_LogGamma}, it follows that the above contour can
be closed over the positive $\operatorname{Re}\left(  w\right)  $-plane,
resulting in\footnote{The integral (\ref{complex_integral_cosec_type}) is
rewritten with the reflection formula of the Gamma function as a Barnes-Mellin
type integral%
\begin{equation}
E_{a,b}\left(  -z\right)  =-\frac{1}{2\pi i}\int_{c+\infty e^{-i\theta}%
}^{c+\infty e^{i\theta}}\frac{\Gamma\left(  1-w\right)  \Gamma\left(
w\right)  }{\Gamma\left(  b+aw\right)  }z^{w}dw\hspace{1cm}\text{for
}-1<c<0\text{ and }0<\theta<\frac{\pi}{2} \label{complex_Barnes_type}%
\end{equation}
}%
\begin{equation}
E_{a,b}\left(  z\right)  =-\frac{1}{2\pi i}\int_{c+\infty e^{-i\theta}%
}^{c+\infty e^{i\theta}}\frac{\pi}{\sin\pi w}\frac{\left(  -z\right)  ^{w}%
}{\Gamma\left(  b+aw\right)  }dw\hspace{1cm}\text{for }-1<c<0
\label{complex_integral_cosec_type}%
\end{equation}
where the path above and below the real axis follow the lines $c+re^{\pm
i\theta}$, where $0<\theta<\frac{\pi}{2}$. Thus, the line of integration
cannot be parallel with the imaginary axis, unless $a<2$. If $a=b=1$, the
reflection formula (\ref{Gamma_reflection_formula}) leads to Mellin's integral
$\frac{1}{2\pi i}\int_{c+\infty}^{c+\infty}\Gamma\left(  w\right)  \left(
-z\right)  ^{w}dw=e^{z}$.

Let us consider the contour $\mathcal{C}$, consisting of the line
$c+re^{i\theta}$ with $0\leq r\leq T$, the line parallel to the real axis from
$c+Te^{i\theta}$ to the left at $c-m+Te^{i\theta}$, the line back to the real
axis at $c-m$ (and the reflection of this parallelogram around the real axis).
The path parallel to the real axis%
\[
\int_{c+Te^{-i\theta}}^{c-m+Te^{i\theta}}\frac{\pi}{\sin\pi w}\frac{\left(
-z\right)  ^{w}}{\Gamma\left(  b+aw\right)  }dw=\int_{c}^{c-m}\frac{\pi}%
{\sin\pi\left(  x+Te^{-i\theta}\right)  }\frac{\left(  -z\right)
^{x+Te^{-i\theta}}}{\Gamma\left(  b+ax+aTe^{-i\theta}\right)  }dx
\]
vanishes by (\ref{asymptotic_Gamma(b+rexp(it))}) for $T\rightarrow\infty$
provided $0<\theta<\frac{\pi}{2}$. The contour $\mathcal{C}$ encloses the
poles $\frac{\pi}{\sin\pi w}$ at $w=-k$ from $k=1$ to $m$ with residue
$\left(  -1\right)  ^{k}$. Hence, by shifting the lines $c+re^{\pm i\theta}$
to $c-m+re^{\pm i\theta}$, maintaining the angle $0<\theta<\frac{\pi}{2}$, we
deform the integral (\ref{complex_integral_cosec_type}) into%
\[
E_{a,b}\left(  z\right)  =\sum_{k=1}^{m}\frac{1}{\Gamma\left(  b-ka\right)
}\frac{1}{z^{k}}-\frac{1}{2\pi i}\int_{c^{\prime}+\infty e^{-i\theta}%
}^{c^{\prime}+\infty e^{i\theta}}\frac{\pi}{\sin\pi w}\frac{\left(  -z\right)
^{w}}{\Gamma\left(  b+aw\right)  }dw\hspace{1cm}\text{for }-1-m<c^{\prime}<m
\]
For complex $z$, it is generally for any $a$ complicated to bound the integral
for large $\left\vert z\right\vert $ to deduce an asymptotic expansion.

\medskip\refstepcounter{article}{\noindent\textbf{\thearticle. }%
}\ignorespaces
We evaluate (\ref{complex_integral_cosec_type}) along the line
$w=c+re^{-i\theta}$ below the real $w$-axis to the line $w=c+re^{i\theta}$
above the real $w$-axis for $0<\theta<\frac{\pi}{2}$,%
\begin{align*}
E_{a,b}\left(  -z\right)   &  =\frac{e^{-i\theta}}{2i}\int_{0}^{\infty}%
\frac{z^{c+re^{-i\theta}}}{\Gamma\left(  b+ac+are^{-i\theta}\right)  }%
\frac{dr}{\sin\pi\left(  c+re^{-i\theta}\right)  }-\frac{e^{i\theta}}{2i}%
\int_{0}^{\infty}\frac{z^{c+re^{i\theta}}}{\Gamma\left(  b+ac+are^{i\theta
}\right)  }\frac{dr}{\sin\pi\left(  c+re^{i\theta}\right)  }\\
&  =\frac{1}{2}\int_{0}^{\infty}\frac{z^{c+re^{-i\theta}}e^{-i\left(
\frac{\pi}{2}+\theta\right)  }}{\Gamma\left(  b+ac+are^{-i\theta}\right)
}\frac{dr}{\sin\pi\left(  c+re^{-i\theta}\right)  }+\frac{1}{2}\int
_{0}^{\infty}\frac{z^{c+re^{i\theta}}e^{i\left(  \frac{\pi}{2}+\theta\right)
}}{\Gamma\left(  b+ac+are^{i\theta}\right)  }\frac{dr}{\sin\pi\left(
c+re^{i\theta}\right)  }%
\end{align*}
Since the integrand $\frac{\pi}{\sin\pi w}\frac{\left(  -z\right)  ^{w}%
}{\Gamma\left(  b+aw\right)  }$ is real on the real axis, the reflection
principle \cite[p. 155]{Titchmarshfunctions} leads, for real $z=x$, to%
\[
E_{a,b}\left(  -x\right)  =x^{c}\int_{0}^{\infty}\operatorname{Re}\left\{
\frac{x^{re^{i\theta}}e^{i\left(  \frac{\pi}{2}+\theta\right)  }}%
{\Gamma\left(  b+ac+are^{i\theta}\right)  \sin\pi\left(  c+re^{i\theta
}\right)  }\right\}  dr
\]
which slightly simplifies, after choosing $c=-\frac{1}{2}$, to%
\[
E_{a,b}\left(  -x\right)  =-\frac{1}{\sqrt{x}}\int_{0}^{\infty}%
\operatorname{Re}\left\{  \frac{e^{i\left(  \frac{\pi}{2}+\theta\right)
}e^{\left(  r\log x\right)  e^{i\theta}}}{\Gamma\left(  b-\frac{a}%
{2}+are^{i\theta}\right)  \cos\pi re^{i\theta}}\right\}  dr
\]

Choosing $\theta=\frac{\pi}{2}$, and thus restricting $0<a<2$, yields%
\[
E_{a,b}\left(  -x\right)  =\frac{1}{\sqrt{x}}\int_{0}^{\infty}%
\operatorname{Re}\left\{  \frac{e^{ir\log x}}{\Gamma\left(  b-\frac{a}%
{2}+air\right)  }\right\}  \frac{dr}{\cosh\pi r}%
\]
Rather than continuing with this integral, we present another approach.

\medskip\refstepcounter{article}{\noindent\textbf{\thearticle. }%
}\ignorespaces\label{art_1_op_sinpiw_integral_for_0<a<1}If we choose
$\theta=\frac{\pi}{2}$ in (\ref{complex_integral_cosec_type}), then we must
restrict $0<a<2$. In that case, we evaluate the contour along the line
$w=c+it$,%
\[
E_{a,b}\left(  -z\right)  =-\frac{z^{c}}{2}\int_{-\infty}^{\infty}\frac
{1}{\sin\pi\left(  c+it\right)  }\frac{z^{it}}{\Gamma\left(  b+ac+iat\right)
}dt\hspace{1cm}\text{for }-1<c<0
\]
We limit ourselves to real $z=x$. Since $c<0$, the reflection formula
(\ref{Gamma_reflection_formula}) yields
\[
\frac{1}{\Gamma\left(  b-a\left\vert c\right\vert +ait\right)  }%
=-\Gamma\left(  1-b+a\left\vert c\right\vert -ait\right)  \frac{\sin\pi\left(
a\left\vert c\right\vert -b-ait\right)  }{\pi}%
\]
We change the sign of $c$ and obtain%
\[
E_{a,b}\left(  -x\right)  =\frac{1}{2\pi}\int_{-\infty}^{\infty}\frac{\sin
\pi\left(  a\left(  c-it\right)  -b\right)  }{\sin\pi\left(  c-it\right)
}\frac{\Gamma\left(  ac+1-b-ait\right)  }{x^{c-it}}dt\hspace{1cm}\text{for
}0<c<1
\]
In order to invoke Euler's integral (\ref{Euler_integral_Gamma_functie}), we
must require that $\operatorname{Re}\left(  ac+1-b\right)  >0$%
\[
E_{a,b}\left(  -x\right)  =\frac{x^{\frac{1-b}{a}}}{2\pi}\int_{-\infty
}^{\infty}dt\frac{\sin\pi\left(  a\left(  c-it\right)  -b\right)  }{\sin
\pi\left(  c-it\right)  }\int_{0}^{\infty}due^{-x^{\frac{1}{a}}u}u^{ac-b-ait}%
\]
We can interchange the integrals by absolute convergence requiring that
$0<a<1$,%
\begin{equation}
E_{a,b}\left(  -x\right)  =\frac{x^{\frac{1-b}{a}}}{2\pi}\int_{0}^{\infty
}du\frac{e^{-x^{\frac{1}{a}}u}}{u^{b}}\int_{-\infty}^{\infty}dt\frac{\sin
\pi\left(  a\left(  c-it\right)  -b\right)  }{\sin\pi\left(  c-it\right)
}u^{a\left(  c-it\right)  } \label{E_a,b(-x)_double_integral}%
\end{equation}
In passing, we check for $a\rightarrow1$ and $b=1$ that, by using the Dirac
impulse function $\delta\left(  x\right)  $,%
\[
E_{1}\left(  -x\right)  =\frac{1}{2\pi}\int_{0}^{\infty}due^{-xu}u^{c-1}%
\int_{-\infty}^{\infty}dte^{-it\log u}=\int_{0}^{\infty}due^{-xu}u^{c-1}%
\delta\left(  \log u\right)  =e^{-x}%
\]

The latter integral in (\ref{E_a,b(-x)_double_integral}) can be recasted as%
\begin{equation}
Q_{a,b}\left(  u\right)  =-\frac{1}{2\pi}\int_{-\infty}^{\infty}dt\frac
{\sin\pi\left(  a\left(  c-it\right)  -b\right)  }{\sin\pi\left(  c-it\right)
}u^{a\left(  c-it\right)  }=\frac{1}{2\pi i}\int_{c-i\infty}^{c+i\infty}%
\frac{\sin\pi\left(  aw-b\right)  }{\sin\pi w}\left(  u^{a}\right)
^{w}dw\hspace{1cm}\text{for }0<c<1 \label{defQ_integral}%
\end{equation}
In \textbf{art}. \ref{art_integral_Q_a,b(u)} below, we prove for real $u$ and
$0<a<1$ that%
\begin{equation}
Q_{a,b}\left(  u\right)  =\frac{1}{2\pi i}\int_{c-i\infty}^{c+i\infty}%
\frac{\sin\pi\left(  aw-b\right)  }{\sin\pi w}\left(  u^{a}\right)
^{w}dw=\frac{1}{\pi}\left(  \frac{u^{a}\sin\pi\left(  a-b\right)  -u^{2a}%
\sin\pi b}{1+2u^{a}\cos\pi a+u^{2a}}\right)  \label{Integral_minusQ_exact}%
\end{equation}
Substitution of (\ref{Integral_minusQ_exact}) into
(\ref{E_a,b(-x)_double_integral}) yields, for $0<a\leq1$ and
$a+1>\operatorname{Re}\left(  b\right)  $,%
\begin{align}
E_{a,b}\left(  -x\right)   &  =-\frac{x^{\frac{1-b}{a}}\sin\pi\left(
a-b\right)  }{\pi}\int_{0}^{\infty}\frac{u^{a-b}}{1+2u^{a}\cos\pi a+u^{2a}%
}e^{-x^{\frac{1}{a}}u}du\nonumber\\
&  +\frac{x^{\frac{1-b}{a}}\sin\pi b}{\pi}\int_{0}^{\infty}\frac{u^{2a-b}%
}{1+2u^{a}\cos\pi a+u^{2a}}e^{-x^{\frac{1}{a}}u}du
\label{E_a,b(-x)_integral_a_in(0,1]}%
\end{align}
which equals (\ref{Mittag_Leffler_function_as_residues}) without residue sum
(after the substitution $y=x^{\frac{1}{a}}u$). If $b=1$, then
(\ref{E_a,b(-x)_integral_a_in(0,1]}) simplifies to
\begin{equation}
E_{a}\left(  -x\right)  =\frac{\sin\pi a}{\pi}\int_{0}^{\infty}\frac{u^{a-1}%
}{1+2u^{a}\cos\pi a+u^{2a}}e^{-x^{\frac{1}{a}}u}du
\label{E_a(-x)_integral_a_in(0,1]}%
\end{equation}
which is listed in \cite[eq. (34)]{Berberan-SantosII_2005} and deducible from
(\ref{Mittag_Leffler_E_a(-x^a)_0<a<1_integral}).

\medskip\refstepcounter{article}{\noindent\textbf{\thearticle. }%
}\ignorespaces\label{art_integral_Q_a,b(u)} \emph{The integral }%
$Q_{a,b}\left(  u\right)  $. The contour of the integral $Q_{a,b}\left(
u\right)  $ in (\ref{defQ_integral}) for real $u$ can be closed over either
half-plane. If $\left\vert u\right\vert <1$, then we close the contour over
$\operatorname{Re}\left(  w\right)  >0$-plane and obtain%
\begin{align*}
\frac{1}{2\pi i}\int_{c-i\infty}^{c+i\infty}dt\frac{\sin\pi\left(
aw-b\right)  }{\sin\pi w}\left(  u^{a}\right)  ^{w}  &  =-\frac{1}{\pi}%
\sum_{k=1}^{\infty}\left(  -1\right)  ^{k}\sin\pi\left(  ak-b\right)  u^{ak}\\
&  =\frac{1}{2i\pi}\left(  e^{i\pi\left(  a-b\right)  }u^{a}\sum_{k=0}%
^{\infty}\left(  -e^{i\pi a}u^{a}\right)  ^{k}-e^{-i\pi\left(  a-b\right)
}u^{a}\sum_{k=0}^{\infty}\left(  -e^{-i\pi a}u^{a}\right)  ^{k}\right) \\
&  =\frac{u^{a}}{2i\pi}\left(  \frac{e^{i\pi\left(  a-b\right)  }}{1+e^{i\pi
a}u^{a}}-\frac{e^{-i\pi\left(  a-b\right)  }}{1+e^{-i\pi a}u^{a}}\right) \\
&  =\frac{1}{\pi}\frac{u^{a}\sin\pi\left(  a-b\right)  -u^{2a}\sin\pi
b}{1+2u^{a}\cos\pi a+u^{2a}}%
\end{align*}
If $\left\vert u\right\vert >1$, then we close the contour of the
$\operatorname{Re}\left(  w\right)  <0$ plane,%
\begin{align*}
\frac{1}{2\pi i}\int_{c-i\infty}^{c+i\infty}dt\frac{\sin\pi\left(
aw-b\right)  }{\sin\pi w}\left(  u^{a}\right)  ^{w}  &  =\frac{1}{\pi}%
\sum_{k=0}^{\infty}\left(  -1\right)  ^{k}\sin\pi\left(  -ak-b\right)
u^{-ak}\\
&  =\frac{-1}{2i\pi}\left(  \frac{e^{i\pi b}}{1+e^{i\pi a}u^{-a}}%
-\frac{e^{-i\pi b}}{1+e^{-i\pi a}u^{-a}}\right) \\
&  =\frac{1}{\pi}\left(  \frac{u^{a}\sin\pi\left(  a-b\right)  -u^{2a}\sin\pi
b}{1+2u^{a}\cos\pi a+u^{2a}}\right)
\end{align*}
Both the case $\left\vert u\right\vert <1$ and $\left\vert u\right\vert >1$
lead to the same result! In summary, for any real $u$ and for $0<a<1$, we have
proved (\ref{Integral_minusQ_exact}).

The contour in (\ref{defQ_integral}) can be rewritten as%
\[
Q_{a,b}\left(  u\right)  =\frac{u^{\frac{a}{2}}}{\pi}\left\{  \sin\pi\left(
\frac{a}{2}-b\right)  \int_{0}^{\infty}\frac{\cosh\pi at}{\cosh\pi t}%
\cos\left(  at\log u\right)  dt-\cos\pi\left(  \frac{a}{2}-b\right)  \int
_{0}^{\infty}\frac{\sinh\pi at}{\cosh\pi t}\sin\left(  at\log u\right)
dt\right\}
\]
Since $Q_{a,1}\left(  u\right)  =Q_{a,1}\left(  \frac{1}{u}\right)  $, as
follows from (\ref{Integral_minusQ_exact}), it holds that%

\[
\int_{0}^{\infty}\frac{\sinh\pi at}{\cosh\pi t}\sin\left(  at\log u\right)
dt=\frac{u^{a}-1}{u^{a}+1}\tan\left(  \frac{\pi a}{2}\right)  \int_{0}%
^{\infty}\frac{\cosh\pi at}{\cosh\pi t}\cos\left(  at\log u\right)  dt
\]
so that%
\[
Q_{a,b}\left(  u\right)  =\frac{u^{\frac{a}{2}}}{\pi}\sin\pi\left(  \frac
{a}{2}-b\right)  \left\{  1-\frac{u^{a}-1}{u^{a}+1}\cot\pi\left(  \frac{a}%
{2}-b\right)  \tan\left(  \frac{\pi a}{2}\right)  \right\}  \int_{0}^{\infty
}\frac{\cosh\pi at}{\cosh\pi t}\cos\left(  at\log u\right)  dt
\]
It follows then from (\ref{Integral_minusQ_exact}) that%
\[
\int_{0}^{\infty}\frac{\cosh\pi at}{\cosh\pi t}\cos\left(  at\log u\right)
dt=\frac{u^{\frac{a}{2}}\left(  \sin\pi\left(  a-b\right)  -u^{a}\sin\pi
b\right)  }{\left(  1+2u^{a}\cos\pi a+u^{2a}\right)  \left\{  \sin\pi\left(
\frac{a}{2}-b\right)  +\frac{1-u^{a}}{1+u^{a}}\cos\pi\left(  \frac{a}%
{2}-b\right)  \tan\left(  \frac{\pi a}{2}\right)  \right\}  }%
\]
whose right-hand is independent of $b$, so that $b$ can be chosen at will. The
simplest choice\footnote{It is readily verified that%
\[
\frac{\sin\pi\left(  a-b\right)  -u^{a}\sin\pi b}{\sin\pi\left(  \frac{a}%
{2}-b\right)  +\frac{1-u^{a}}{1+u^{a}}\cos\pi\left(  \frac{a}{2}-b\right)
\tan\left(  \frac{\pi a}{2}\right)  }=\cos\left(  \frac{\pi a}{2}\right)
\left(  1+u^{a}\right)
\]
} is $b=\frac{a}{2}$, then%
\[
\int_{0}^{\infty}\frac{\cosh\pi at}{\cosh\pi t}\cos\left(  at\log u\right)
dt=\frac{u^{\frac{a}{2}}\left(  1+u^{a}\right)  \cos\frac{\pi a}{2}}%
{1+2u^{a}\cos\pi a+u^{2a}}%
\]

Suppose that we ignore the restriction that
\[
\frac{1}{2\pi i}\int_{c-i\infty}^{c+i\infty}dt\frac{\sin\pi\left(
aw-b\right)  }{\sin\pi w}\left(  u^{a}\right)  ^{w}=\frac{1}{\pi}\sum
_{k=1}^{\infty}\left(  -1\right)  ^{k-1}\sin\pi\left(  ak-b\right)  u^{ak}%
\]
only converges for $u<1$ and that we substitute the series formally back in
(\ref{E_a,b(-x)_double_integral}) and change the order of integration and
summation. Then, we find%
\[
E_{a,b}\left(  -x\right)  =\frac{1}{\pi}\sum_{k=1}^{\infty}\frac{\sin
\pi\left(  ak-b\right)  \Gamma\left(  ak-b+1\right)  }{\left(  -x\right)
^{k}}=-\sum_{k=1}^{\infty}\frac{1}{\left(  -x\right)  ^{k}\Gamma\left(
b-ak\right)  }=\frac{1}{\Gamma\left(  b\right)  }-E_{-a,b}\left(  -\frac{1}%
{x}\right)
\]
which is precisely equal to (\ref{Mittag_Leffler_negative_a}), in spite of the
divergence of the series!

\medskip\refstepcounter{article}{\noindent\textbf{\thearticle. }%
}\ignorespaces\label{art_alternative_integrals_0<a<1}\emph{Deductions from the
integral (\ref{E_a(-x)_integral_a_in(0,1]}).} Berberan-Santos \cite[eq.
(35)]{Berberan-SantosII_2005} mentions\footnote{The integral in
Berberan-Santos \cite[eq. (35)]{Berberan-SantosII_2005}, $E_{a}\left(
-x\right)  =1-\frac{1}{2a}+\frac{x^{\frac{1}{a}}}{\pi}\int_{0}^{\infty}%
\arctan\left(  \frac{u^{a}+\cos\left(  \pi a\right)  }{\sin\left(  \pi
a\right)  }\right)  e^{-x^{\frac{1}{a}}u}du$, misses a factor of $\frac{1}{a}$
before the integral.} that partial integration of
(\ref{E_a(-x)_integral_a_in(0,1]}) results in%
\begin{equation}
E_{a}\left(  -x\right)  =1-\frac{1}{2a}+\frac{x^{\frac{1}{a}}}{\pi a}\int
_{0}^{\infty}\arctan\left(  \frac{u^{a}+\cos\left(  \pi a\right)  }%
{\sin\left(  \pi a\right)  }\right)  e^{-x^{\frac{1}{a}}u}du
\label{E_a(-x)_integral_a_in(0,1]_arctan}%
\end{equation}
Indeed, let $y=u^{a}$ in (\ref{E_a(-x)_integral_a_in(0,1]}), then%
\[
E_{a}\left(  -x\right)  =\frac{\sin\pi a}{\pi a}\int_{0}^{\infty}%
\frac{e^{-x^{\frac{1}{a}}y^{\frac{1}{a}}}}{1+2y\cos\pi a+y^{2}}dy
\]
With $\frac{1}{1+2y\cos\pi a+y^{2}}=\frac{1}{2i\sin\pi a}\left(  \frac
{1}{e^{-i\pi a}+y}-\frac{1}{e^{i\pi a}+y}\right)  $, we have%
\begin{equation}
E_{a}\left(  -x\right)  =\frac{1}{2\pi ia}\int_{0}^{\infty}\left(  \frac
{1}{e^{-i\pi a}+y}-\frac{1}{e^{i\pi a}+y}\right)  e^{-x^{\frac{1}{a}}%
y^{\frac{1}{a}}}dy \label{E_a(-x)_integral_a_in(0,1]_partial_fractions}%
\end{equation}
as well as%
\begin{align*}
\int\frac{dy}{1+2y\cos\pi a+y^{2}}  &  =\frac{1}{2i\sin\pi a}\left(  \log
\frac{e^{-i\pi a}+y}{e^{i\pi a}+y}\right)  =\frac{1}{2i\sin\pi a}\left(
\log\frac{y+\cos\pi a-i\sin\pi a}{y+\cos\pi a+i\sin\pi a}\right) \\
&  =\frac{1}{2i\sin\pi a}\left(  \log\frac{1-i\frac{\sin\pi a}{y+\cos\pi a}%
}{1+i\frac{\sin\pi a}{y+\cos\pi a}}\right)  =-\frac{1}{\sin\pi a}%
\arctan\left(  \frac{\sin\pi a}{y+\cos\pi a}\right)
\end{align*}
because $\arctan z=\frac{i}{2}\log\frac{1-iz}{1+iz}$. Hence, we find the
indefinite integral with a constant $K$,
\begin{equation}
\int\frac{u^{a-1}du}{1+2u^{a}\cos\pi a+u^{2a}}=-\frac{1}{a\sin\pi a}%
\arctan\left(  \frac{\sin\pi a}{u^{a}+\cos\pi a}\right)  +K
\label{onbepaalde_integral}%
\end{equation}
With this preparation, partial integration of
(\ref{E_a(-x)_integral_a_in(0,1]}) yields%
\begin{align*}
E_{a}\left(  -x\right)   &  =\frac{1}{a\pi}\arctan\left(  \frac{\sin\pi
a}{\cos\pi a}\right)  -\frac{x^{\frac{1}{a}}}{a\pi}\int_{0}^{\infty}%
\arctan\left(  \frac{\sin\pi a}{u^{a}+\cos\pi a}\right)  e^{-x^{\frac{1}{a}}%
u}du\\
&  =1-\frac{x^{\frac{1}{a}}}{a\pi}\int_{0}^{\infty}\arctan\left(  \frac
{\sin\pi a}{u^{a}+\cos\pi a}\right)  e^{-x^{\frac{1}{a}}u}du
\end{align*}
After invoking $\arctan x=\frac{\pi}{2}-\arctan\frac{1}{x}$ for $x>0$, we
arrive at (\ref{E_a(-x)_integral_a_in(0,1]_arctan}).

Further, using the integral for $\arctan z=\int_{0}^{z}\frac{dt}{1+t^{2}}$,
\[
1-E_{a}\left(  -x\right)  =\frac{x^{\frac{1}{a}}}{a\pi}\int_{0}^{\infty}%
du\int_{0}^{\frac{\sin\pi a}{u^{a}+\cos\pi a}}dt\frac{e^{-x^{\frac{1}{a}}u}%
}{1+t^{2}}%
\]
and reverse the integrals, provided that $0<a\leq\frac{1}{2}$,%
\[
1-E_{a}\left(  -x\right)  =\frac{x^{\frac{1}{a}}}{a\pi}\int_{0}^{\tan\pi
a}\frac{\int_{0}^{\left(  \frac{\sin\pi a}{t}-\cos\pi a\right)  ^{\frac{1}{a}%
}}e^{-x^{\frac{1}{a}}u}du}{1+t^{2}}dt
\]
Hence, we obtain, for $0<a\leq\frac{1}{2}$,%
\begin{equation}
E_{a}\left(  -x\right)  =\frac{1}{a\pi}\int_{0}^{\tan\pi a}\frac{e^{-\left(
x\left(  \frac{\sin\pi a}{t}-\cos\pi a\right)  \right)  ^{\frac{1}{a}}}%
}{1+t^{2}}dt \label{E_a(-x)_integral_a_in(0,1/2]_finite}%
\end{equation}
After letting $u^{a}=\frac{\sin\pi a}{t}-\cos\pi a$ in
(\ref{E_a(-x)_integral_a_in(0,1/2]_finite}), we retrieve
(\ref{E_a(-x)_integral_a_in(0,1]}) again.

\medskip\refstepcounter{article}{\noindent\textbf{\thearticle. }%
}\ignorespaces\label{art_bounds_integrals_0<a<1}\emph{Bounds from the integral
(\ref{E_a(-x)_integral_a_in(0,1]}).} We split the integration interval in
(\ref{E_a(-x)_integral_a_in(0,1]}) into two parts,%
\[
E_{a}\left(  -x\right)  =\frac{\sin\pi a}{\pi}\int_{1}^{\infty}\frac
{u^{a-1}\left(  e^{-x^{\frac{1}{a}}u^{-1}}+e^{-x^{\frac{1}{a}}u}\right)
}{1+2u^{a}\cos\pi a+u^{2a}}du
\]
First,%
\begin{align*}
\int_{1}^{\infty}\frac{u^{a-1}e^{-x^{\frac{1}{a}}u^{-1}}}{1+2u^{a}\cos\pi
a+u^{2a}}du  &  =e^{-x^{\frac{1}{a}}}\int_{1}^{\infty}\frac{u^{a-1}%
e^{x^{\frac{1}{a}}\left(  1-u^{-1}\right)  }}{1+2u^{a}\cos\pi a+u^{2a}}du\\
&  \geq e^{-x^{\frac{1}{a}}}\int_{1}^{\infty}\frac{u^{a-1}}{1+2u^{a}\cos\pi
a+u^{2a}}du
\end{align*}
Using (\ref{onbepaalde_integral}) and $\frac{\sin\pi a}{1+\cos\pi a}=\tan
\frac{\pi a}{2}$, we arrive at the lower bound,%
\[
\frac{\sin\pi a}{\pi}\int_{1}^{\infty}\frac{u^{a-1}e^{-x^{\frac{1}{a}}u^{-1}}%
}{1+2u^{a}\cos\pi a+u^{2a}}du\geq\frac{1}{\pi a}e^{-x^{\frac{1}{a}}}%
\arctan\left(  \tan\frac{\pi a}{2}\right)  =\frac{1}{2}e^{-x^{\frac{1}{a}}}%
\]
Second, after partial integration and again using (\ref{onbepaalde_integral}),
we obtain
\begin{align*}
\frac{\sin\pi a}{\pi}\int_{1}^{\infty}\frac{u^{a-1}e^{-x^{\frac{1}{a}}u}%
}{1+2u^{a}\cos\pi a+u^{2a}}du  &  =\frac{1}{2}e^{-x^{\frac{1}{a}}}%
-\frac{x^{\frac{1}{a}}}{\pi a}\int_{1}^{\infty}e^{-x^{\frac{1}{a}}u}%
\arctan\left(  \frac{\sin\pi a}{u^{a}+\cos\pi a}\right)  du\\
&  =\frac{1}{2}e^{-x^{\frac{1}{a}}}-\frac{x^{\frac{1}{a}}e^{-x^{\frac{1}{a}}}%
}{\pi a}\int_{1}^{\infty}e^{-x^{\frac{1}{a}}\left(  u-1\right)  }%
\arctan\left(  \frac{\sin\pi a}{u^{a}+\cos\pi a}\right)  du
\end{align*}
so that%
\[
E_{a}\left(  -x\right)  \geq e^{-x^{\frac{1}{a}}}-\frac{x^{\frac{1}{a}%
}e^{-x^{\frac{1}{a}}}}{\pi a}\int_{1}^{\infty}e^{-x^{\frac{1}{a}}\left(
u-1\right)  }\arctan\left(  \frac{\sin\pi a}{u^{a}+\cos\pi a}\right)  du
\]
Further, we may bound the latter integral,%
\begin{align*}
\int_{1}^{\infty}e^{-x^{\frac{1}{a}}\left(  u-1\right)  }\arctan\left(
\frac{\sin\pi a}{u^{a}+\cos\pi a}\right)  du  &  <\int_{1}^{\infty
}e^{-x^{\frac{1}{a}}\left(  u-1\right)  }\arctan\left(  \frac{\sin\pi
a}{1+\cos\pi a}\right)  du\\
&  =\frac{\pi a}{2}\int_{1}^{\infty}e^{-x^{\frac{1}{a}}\left(  u-1\right)
}du=\frac{\pi a}{2x^{\frac{1}{a}}}%
\end{align*}
which leads to the lower bound, for $0<a<1$,%
\[
E_{a}\left(  -x\right)  \geq\frac{1}{2}e^{-x^{\frac{1}{a}}}%
\]

It follows directly from (\ref{E_a(-x)_integral_a_in(0,1]}) that%
\begin{align*}
E_{a}\left(  -x\right)   &  =\frac{\sin\pi a}{\pi}\int_{0}^{\infty}\left(
1-\frac{2u^{a}\cos\pi a+u^{2a}}{1+2u^{a}\cos\pi a+u^{2a}}\right)
u^{a-1}e^{-x^{\frac{1}{a}}u}du\\
&  =\frac{\sin\pi a}{\pi}\frac{\Gamma\left(  a\right)  }{x}-\frac{\sin\pi
a}{\pi}\int_{0}^{\infty}\frac{\left(  2\cos\pi a+u^{a}\right)  }{1+2u^{a}%
\cos\pi a+u^{2a}}u^{2a-1}e^{-x^{\frac{1}{a}}u}du
\end{align*}
Further, for $0<a<\frac{1}{2}$ (because then $\cos\pi a\geq0$),%
\begin{align*}
E_{a}\left(  -x\right)   &  \leq\frac{\sin\pi a}{\pi}\frac{\Gamma\left(
a\right)  }{x}-\frac{\sin\pi a}{\pi}\int_{0}^{1}\frac{\left(  2\cos\pi
a+u^{a}\right)  }{1+2u^{a}\cos\pi a+u^{2a}}u^{2a-1}e^{-x^{\frac{1}{a}}u}du\\
&  \leq\frac{\sin\pi a}{\pi}\frac{\Gamma\left(  a\right)  }{x}-\frac{\sin\pi
a}{\pi}\frac{2\cos\pi a}{2+2\cos\pi a}e^{-x^{\frac{1}{a}}}\int_{0}^{1}%
u^{2a-1}du
\end{align*}
and%
\[
E_{a}\left(  -x\right)  \leq\frac{\sin\pi a}{\pi}\frac{\Gamma\left(  a\right)
}{x}-\frac{e^{-x^{\frac{1}{a}}}}{4a\pi}\frac{\sin2\pi a}{1+\cos\pi a}%
\]
which illustrates, for $0<a<\frac{1}{2}$, that%
\begin{equation}
\frac{\sin\pi a}{\pi}\frac{\Gamma\left(  a\right)  }{x}-\frac{e^{-x^{\frac
{1}{a}}}}{4a\pi}\frac{\sin2\pi a}{1+\cos\pi a}\geq E_{a}\left(  -x\right)
\geq\frac{1}{2}e^{-x^{\frac{1}{a}}} \label{upper_and_lower_bound_Ea_0<a<1/2}%
\end{equation}

\section{The Mittag-Leffler function in probability theory}

\label{sec_Mittag_Leffler_probability_theory}

\refstepcounter{article}{\noindent\textbf{\thearticle.}}%
\ignorespaces\label{art_monotonicity_Mittag_Leffler_function}
\emph{Monotonicity of }$E_{a,b}\left(  -x\right)  $\emph{ for }$0\leq a\leq1$.
Widder \cite[Chapter IV]{Widder} devotes an entire chapter to absolutely and
completely monotone functions. The
Hausdorff--Bernstein--Widder\footnote{Widder \cite[p. 144]{Widder} mentions
that Hausdorff, Bernstein and himself have independently proved the theorem.}
theorem \cite[p. 161, Theorem 12a]{Widder} states that a necessary and
sufficient condition that a function $\varphi\left(  x\right)  $ on
$[0,\infty)$ should be completely monotonic is that there exists a bounded and
non-decreasing function $f\left(  u\right)  $ such that the integral%
\[
\varphi\left(  x\right)  =\int_{0}^{\infty}e^{-xu}df\left(  u\right)
\]
converges for all real $x\geq0$. In other words, a function $\varphi\left(
x\right)  $ is completely monotonic on $[0,\infty)$ if and only if
$\varphi\left(  x\right)  $ is a Laplace transform of a bounded and
non-decreasing measure $f\left(  u\right)  $. The fact that $\left(
-1\right)  ^{n}\frac{d^{n}\varphi\left(  x\right)  }{dx^{n}}=\int_{0}^{\infty
}x^{n}e^{-xu}df\left(  u\right)  \geq0$ for all non-negative integers $n$ is a
direct consequence of the Hausdorff--Bernstein--Widder theorem, but the
condition $\left(  -1\right)  ^{n}\frac{d^{n}\varphi\left(  x\right)  }%
{dx^{n}}\geq0$ for all non-negative integers $n$ is also sufficient \cite[p.
56-59]{Bernstein_1929} and thus an equivalent statement for complete
monotonicity of a function $\varphi\left(  x\right)  $.

Since the integrand in (\ref{E_a(-x)_integral_a_in(0,1]}) is positive for
$0\leq a\leq1$ because
\[
\left(  1-u^{a}\right)  ^{2}\leq1+2u^{a}\cos\pi a+u^{2a}\leq\left(
1+u^{a}\right)  ^{2}%
\]
(\ref{E_a(-x)_integral_a_in(0,1]}) shows that $E_{a}\left(  -x\right)  >0$ for
$0<a\leq1$. Hence, the integral (\ref{E_a(-x)_integral_a_in(0,1]}) directly
demonstrates that $E_{a}\left(  -x\right)  >0$ for $0<a\leq1$ is completely
monotonic. Also, the case for $b=a$ in (\ref{E_a,b(-x)_integral_a_in(0,1]})
reduces to%
\begin{equation}
x^{1-\frac{1}{a}}E_{a,a}\left(  -x\right)  =\frac{\sin\pi a}{\pi}\int
_{0}^{\infty}\frac{u^{a}}{1+2u^{a}\cos\pi a+u^{2a}}e^{-x^{\frac{1}{a}}u}du
\label{E_a,a(-x)_integral_monotonic}%
\end{equation}
illustrating that also $x^{1-\frac{1}{a}}E_{a,a}\left(  -x\right)  $ is
complete monotonic. The monotonicity of $E_{a}\left(  -x\right)  >0$ for
$0<a\leq1$ was first conjectured by Feller \cite[Section 7]{Feller_TAMS1949}
and later proved by Pollard \cite{Pollard_1948}.

Pollard\footnote{Pollard's PhD advisor was D. V. Widder.} \cite{Pollard_1948}
introduces $\frac{1}{t+z}=\int_{0}^{\infty}e^{-s\left(  t+z\right)  }ds$ in
Mittag-Leffler's integral (\ref{Mittag_Leffler_Bieberbach_integral}) and
obtains a Laplace transform,%
\[
E_{a,b}\left(  -z\right)  =\int_{0}^{\infty}e^{-sz}\left(  \frac{1}{2\pi
ia}\int_{C_{a\phi}}t^{\frac{1-b}{a}}e^{-st+t^{\frac{1}{a}}}dt\right)  ds
\]
from which he proves\footnote{Berberan-Santos \cite{Berberan-SantosII_2005}
claims that monotonicity, defined by $\left(  -1\right)  ^{n}\frac{d^{n}%
}{dx^{n}}E_{a}\left(  -x\right)  \geq0$ for all $n$, follows from $g\left(
t\right)  >0$ in \textbf{art}. \ref{art_Berberan_Santos} for $0\leq a\leq1$,
but I find his argument circular.} that $E_{a}\left(  -x\right)  $ is
completely monotonic for real $x\geq0$ and $0\leq a\leq1$, in the sense that%
\[
E_{a}\left(  -x\right)  =\int_{0}^{\infty}e^{-xt}dF_{a}\left(  t\right)
\]
where $F_{a}\left(  t\right)  $ is nondecreasing, bounded and a probability
distribution. In other words, $E_{a}\left(  z\right)  $ for $0\leq a\leq1$ has
no zeros on the negative real axis. Pollard \cite{Pollard_1948} also
explicitly determined the non-negative function $F_{a}^{\prime}\left(
t\right)  $. However, his proof is not easy and, therefore, omitted, but
replaced by our derivation.

\medskip\refstepcounter{article}{\noindent\textbf{\thearticle.}}%
\ignorespaces\label{art_E_a(-x)_probability_theory} $E_{a}\left(  -x\right)
$\emph{ with } $0<a<1$ \emph{in probability theory.} We construct two
probability density functions from the Mittag-Leffler function $E_{a}\left(
-x\right)  $ with $0<a<1$ and show that $E_{a}\left(  -s^{a}\right)  $ is both
a generating function and a probability distribution.

First, after replacing $x$ by $s^{a}$ in (\ref{E_a(-x)_integral_a_in(0,1]}),
the Laplace transform
\[
E_{a}\left(  -s^{a}\right)  =\frac{\sin\pi a}{\pi}\int_{0}^{\infty}%
\frac{t^{a-1}}{1+2t^{a}\cos\pi a+t^{2a}}e^{-st}dt
\]
indicates that%
\begin{equation}
f_{a}\left(  t\right)  =\frac{\sin\pi a}{\pi}\frac{t^{a-1}}{1+2t^{a}\cos\pi
a+t^{2a}} \label{pdf_fa(t)_pgf_Ea(-s^a)}%
\end{equation}
is a probability density function (pdf) for $t>0$ and for $0<a\leq1$. Indeed,
the Laplace transform (\ref{def_pgf_continuousrv}) of the non-negative
function $f_{a}\left(  t\right)  $ is \cite{PVM_PAComplexNetsCUP}\ a
probability generating function (pgf) $\varphi_{X}(z)=E\left[  e^{-zX}\right]
$ of a random variable $X$, provided $\varphi_{X}(0)=1$. Moreover, with
$f_{a}\left(  1\right)  =\frac{1}{2\pi}\tan\frac{\pi a}{2}$, the pdf
(\ref{pdf_fa(t)_pgf_Ea(-s^a)}) obeys the functional equation -- recall that
$Q_{a,1}\left(  u\right)  =Q_{a,1}\left(  \frac{1}{u}\right)  $ in
(\ref{Integral_minusQ_exact}) --
\[
tf_{a}\left(  t\right)  =\frac{1}{t}f_{a}\left(  \frac{1}{t}\right)
\]
Its companion, that follows similarly from (\ref{E_a,a(-x)_integral_monotonic}%
) as%
\[
s^{a-1}E_{a,a}\left(  -s^{a}\right)  =\frac{\sin\pi a}{\pi}\int_{0}^{\infty
}\frac{t^{a}}{1+2t^{a}\cos\pi a+t^{2a}}e^{-st}dt
\]
with non-negative function
\[
g_{a}\left(  t\right)  =\frac{\sin\pi a}{\pi}\frac{t^{a}}{1+2t^{a}\cos\pi
a+t^{2a}}=tf_{a}\left(  t\right)  =g_{a}\left(  \frac{1}{t}\right)
\]
is not\footnote{However, insertion into (\ref{pgf_Mittag_Leffler}) below
leads, for $\operatorname{Re}\left(  s\right)  >0$, to%
\[
\int_{0}^{\infty}\frac{u^{a}}{1+2u^{a}\cos\pi a+u^{2a}}\frac{du}{s+u}%
=\frac{\pi}{\sin\pi a}\frac{1}{s^{a}+1}%
\]
} a probability generating function, because $\lim_{s\rightarrow0}%
s^{a-1}E_{a,a}\left(  -s^{a}\right)  =\infty$ for $a<1$. We verify from
Laplace transform theory that
\begin{equation}
-\frac{d}{ds}E_{a}\left(  -s^{a}\right)  =s^{a-1}E_{a,a}\left(  -s^{a}\right)
\label{differentiation_Mittag_Leffler_pgf}%
\end{equation}
in agreement with (\ref{D_(z^(b-1)E_a,b(z^a))}) taken into account
(\ref{E_a,0}).

The logarithmic derivative of the pdf (\ref{pdf_fa(t)_pgf_Ea(-s^a)})%
\[
f_{a}^{\prime}\left(  t\right)  =\frac{f_{a}\left(  t\right)  }{t}\left\{
a-1-2at^{a}\left(  \frac{\cos\pi a+t^{a}}{1+2t^{a}\cos\pi a+t^{2a}}\right)
\right\}
\]
shows that $f_{a}^{\prime}\left(  t\right)  =0$ if%
\[
a=1+2a\left(  \frac{t^{a}}{1+2t^{a}\cos\pi a+t^{2a}}\right)  \left(  \cos\pi
a+t^{a}\right)
\]
from which it follows that $f_{a}\left(  t\right)  $ has no extremum for
$a\leq\frac{1}{2}$, because the last term is positive for all $t\geq0$. Thus,
for $0<a\leq\frac{1}{2}$, the pdf $f_{a}\left(  t\right)  $ monotonously
decreases for all $t$. If $1>a>\frac{1}{2}$, then $\cos\pi a<0$ and
\[
\frac{\frac{1}{a}-1}{2\left(  \frac{t^{a}}{1+2t^{a}\cos\pi a+t^{2a}}\right)
}=-\cos\pi a-t^{a}%
\]
shows that the left-hand side is always positive, whereas the right-hand side
is positive for small $t$, but negative for $t>1$, in fact, for $\left(
-\cos\pi a\right)  ^{\frac{1}{a}}<t$. For small $t$, we can approximate the
equation by a quadratic equation,%
\[
2t^{2a}+2t^{a}\cos\pi a+\left(  \frac{1}{a}-1\right)  \simeq0
\]
with zeros%
\[
t^{a}\simeq-\frac{\cos\pi a\pm\sqrt{\cos^{2}\pi a-2\left(  \frac{1}%
{a}-1\right)  }}{2}%
\]
and observe that two extrema are possible for $t\in\left[  0,1\right]  $.

Second, the Laplace transform (\ref{Laplace_transform_positive_argument}) with
$b=1$ and $x=-1$,%
\[
\int_{0}^{\infty}e^{-st}E_{a}\left(  -t^{a}\right)  dt=\frac{s^{a-1}}{s^{a}+1}%
\]
and with $b=a$ and $x=-1$%
\begin{equation}
\int_{0}^{\infty}e^{-st}t^{a-1}E_{a,a}\left(  -t^{a}\right)  dt=\frac{1}%
{s^{a}+1} \label{pgf_Mittag_Leffler}%
\end{equation}
indicate that the latter is a pgf $\varphi_{M}\left(  s\right)  $, because
$\lim_{s\rightarrow0}\frac{1}{s^{a}+1}=1$ for $a>0$, but the first is not. To
the pgf $\varphi_{M}\left(  s\right)  =E\left[  e^{-sM}\right]  $ corresponds
the pdf
\begin{equation}
f_{M}\left(  t\right)  =t^{a-1}E_{a,a}\left(  -t^{a}\right)  \hspace
{1cm}\text{for }0<a<1 \label{pdf_Mittag_Leffler}%
\end{equation}
of a random variable $M$ and (\ref{differentiation_Mittag_Leffler_pgf})
demonstrates that the corresponding Mittag-Leffler distribution for $0<a<1$ is%
\begin{equation}
F_{M}\left(  t\right)  =\Pr\left[  M\leq t\right]  =\int_{0}^{t}f_{M}\left(
u\right)  du=1-E_{a}\left(  -t^{a}\right)  \label{Mittag_Leffler_distribution}%
\end{equation}
with mean $E\left[  M\right]  =-\varphi_{M}^{\prime}\left(  0\right)
=\lim_{s\rightarrow0}\frac{as^{a-1}}{\left(  s^{a}+1\right)  ^{2}}=\infty$. In
fact, for $0<a<1$, the pgf (\ref{pgf_Mittag_Leffler}) is not analytic at
$s=0$, implying that the Taylor series around $s=0$ does not exist, nor any
derivative. Hence, the Mittag-Leffler random variable $M$, defined by the pgf
(\ref{pgf_Mittag_Leffler}) and pdf (\ref{pdf_Mittag_Leffler}) for $0<a<1$,
does not possess any finite moment $E\left[  M^{k}\right]  $. In the limit
$a\rightarrow1$, the Mittag-Leffler random variable $M$ becomes an exponential
random variable with mean 1.

\medskip\refstepcounter{article}{\noindent\textbf{\thearticle.}}%
\ignorespaces\label{art_prob_prop_Mittag-Leffler_rv} \emph{Probabilistic
properties of the Mittag-Leffler random variable. }The sum $S_{n}=\sum
_{j=1}^{n}M_{j}$ of $n$ i.i.d Mittag-Leffler random variables $M_{1}%
,M_{2},\ldots,M_{n}$, each with the same Mittag-Leffler distribution
(\ref{Mittag_Leffler_distribution}), has the pgf \cite[p. 30]%
{PVM_PAComplexNetsCUP}%
\[
\varphi_{S_{n}}\left(  z\right)  =E\left[  e^{-z\sum_{j=1}^{n}M_{j}}\right]
=\varphi_{M}^{n}\left(  z\right)  =\left(  1+z^{a}\right)  ^{-n}%
\]
If we choose the scaling parameter $\beta$ in $\varphi_{S_{n}}\left(  \beta
z\right)  =\left(  1+\beta^{a}z^{a}\right)  ^{-n}$ equal to $\beta^{a}%
=\frac{1}{n}$, then
\[
\lim_{n\rightarrow\infty}\varphi_{S_{n}}\left(  \frac{z}{n^{\frac{1}{a}}%
}\right)  =\lim_{n\rightarrow\infty}\left(  1+\frac{z^{a}}{n}\right)
^{-n}=e^{-z^{a}}%
\]
Hence, the scaled sum $\beta S_{n}=\sum_{j=1}^{n}n^{-\frac{1}{a}}M_{j}$ tends
for $n\rightarrow\infty$ to a random variable $R$ with pgf equal to
$\varphi_{R}(z)=E\left[  e^{-zR}\right]  =e^{-z^{a}}$, whose form belongs to
the class of stable distributions\footnote{Let $X_{1}$ and $X_{2}$ be i.i.d
random variables, similarly distributed as a random variable $X$. The random
variable $X$ is stable if for any constants $a>0$, $b>0$, $c>0$ and $d$, the
random variable $aX_{1}+bX_{2}$ has the same distribution as $cX+d$, denoted
as $aX_{1}+bX_{2}\overset{d}{=}$ $cX+d$. Another definition \cite[p.
170]{FellerII} states that $X$ is stable if and only if $\sum_{j=1}^{n}%
X_{j}\overset{d}{=}c_{n}X+d_{n}$ for any integer $n>1$ and where the constant
$c_{n}>0$ and $d_{n}\in\mathbb{R}$.
\par
Gorenflo and Mainardi \cite{Gorenflo_Mainardi_2012} discuss fractional
diffusion processes and their relation to Levy stable distributions.}. If
$a=1$, then $R=1$ and not random. Since $e^{-z^{a}}$ with $0<1<a$ is only
analytic for $\operatorname{Re}\left(  z\right)  >0$, inverse Laplace
transform (\ref{inverse_Laplace}) provides us with the pdf%
\begin{align*}
f_{R}\left(  t\right)   &  =\frac{1}{2\pi i}\int_{c-\infty e^{-i\phi}%
}^{c+\infty e^{i\phi}}e^{-z^{a}}e^{zt}dz\hspace{1cm}c>0\text{ and }\phi
>\frac{\pi}{2}\\
&  =\sum_{k=0}^{\infty}\frac{\left(  -1\right)  ^{k}}{k!}\frac{1}{2\pi i}%
\int_{c-\infty e^{-i\phi}}^{c+\infty e^{i\phi}}z^{ka}e^{zt}dz
\end{align*}
Introducing Hankel's integral (\ref{Hankel_e^(aw)_a^(1-z)_op_Gamma(z)})%
\[
f_{R}\left(  t\right)  =\sum_{k=0}^{\infty}\frac{\left(  -1\right)  ^{k}}%
{k!}\frac{t^{-ak-1}}{\Gamma\left(  -ka\right)  }%
\]
and the reflection formula (\ref{Gamma_reflection_formula}) results in%
\begin{equation}
f_{R}\left(  t\right)  =\frac{1}{\pi}\sum_{k=1}^{\infty}\frac{\left(
-1\right)  ^{k-1}}{k!}\frac{\Gamma\left(  ka+1\right)  \sin\pi ka}{t^{ak+1}%
}\hspace{1cm}\text{with }0<1<a \label{pdf_stable_R}%
\end{equation}
which is derived in another, more complicated way by Pollard
\cite{Pollard_1946}. Integration leads to the distribution%
\[
1-F_{R}\left(  t\right)  =\Pr\left[  R>t\right]  =\frac{1}{\pi}\sum
_{k=1}^{\infty}\frac{\left(  -1\right)  ^{k-1}}{k!}\frac{\Gamma\left(
ka\right)  \sin\pi ka}{t^{ak}}%
\]
with the interesting result that $\lim_{t\rightarrow0}\Pr\left[  R\geq
t\right]  =1$, while $\lim_{t\rightarrow0}f_{R}\left(  t\right)  =\infty$. The
sum of $n$ i.i.d. random variables with same distribution $f_{R}\left(
t\right)  $ in (\ref{pdf_stable_R}) equals $\sum_{j=1}^{n}R_{j}=n^{\frac{1}%
{a}}R$, because
\[
E\left[  e^{-z\sum_{j=1}^{n}R_{j}}\right]  =\varphi_{R}^{n}\left(  z\right)
=e^{-nz^{a}}=e^{-\left(  n^{\frac{1}{a}}z\right)  ^{a}}=E\left[
e^{-zn^{\frac{1}{a}}R}\right]
\]
Thus, $\sum_{j=1}^{n}R_{j}=n^{\frac{1}{a}}R$ expresses self-similarity: a sum
of random variables maintains the same distribution upon scaling, which is an
alternative description to a \textquotedblleft stable\textquotedblright\ distribution.

We consider the random variable $X=RW$, where $W$ is independent of $R$ and
will be chosen later. The pgf of $X$ is computed by invoking conditional
expectations \cite[p. 32]{PVM_PAComplexNetsCUP},%
\[
\varphi_{X}(z)=E\left[  e^{-zRW}\right]  =E_{W}\left[  E_{R}\left[  \left.
e^{-zRW}\right\vert W\right]  \right]
\]
The inner conditional expectation is a random variable equal to%
\[
E_{R}\left[  \left.  e^{-zRW}\right\vert W\right]  =e^{-z^{a}W^{a}}%
\]
Let us now define the random variable $Y=W^{a}$, then the outer expectation
becomes%
\[
\varphi_{Y}(z)=E_{Y}\left[  e^{-z^{a}Y}\right]  =\varphi_{Y}(z^{a})
\]
If $Y$ is an exponential random variable with mean $\mu$, then $\varphi
_{Y}\left(  z\right)  =E\left[  e^{-zY}\right]  =\frac{1}{z+\mu}$. Hence,
choosing the mean equal $\mu=1$, then shows that%
\[
\varphi_{X}(z)=\frac{1}{z^{a}+1}%
\]
and (\ref{pgf_Mittag_Leffler}) demonstrates that $X=M$ has a Mittag-Leffler
distribution. The random variable $W=Y^{\frac{1}{a}}$ has the distribution
$\Pr\left[  W\leq x\right]  =\Pr\left[  Y^{\frac{1}{a}}\leq x\right]
=\Pr\left[  Y\leq x^{a}\right]  =\left(  1-e^{-x^{a}}\right)  $ with the pdf
\cite[p. 18]{PVM_PAComplexNetsCUP}%
\[
f_{W}\left(  x\right)  =\frac{d\Pr\left[  W\leq x\right]  }{dx}=ax^{a-1}%
e^{-x^{a}}%
\]
which illustrates that $W$ is a Weibull random variable \cite[p.
59]{PVM_PAComplexNetsCUP} with $E\left[  W^{b}\right]  =\Gamma\left(  \frac
{b}{a}+1\right)  $ for any real $b>-a$. Hence, all moments $E\left[
W^{k}\right]  $ for non-negative integer $k$ exist. The Weibull distribution
is one of the three possible limit extremal distributions of a sequence of
i.i.d. random variables \cite[pp. 65-69]{Berger} and reduces for $a=1$ to the
exponential distribution, just like the Mittag-Leffler distribution
(\ref{Mittag_Leffler_distribution}). Here, we have shown for the parameter
$0<a<1$ that the Weibull random variable $W=\frac{M}{R}$ is the quotient of
two random variables, whose moments do not exist.

Based on Laplace transforms, Feller \cite[p. 453]{FellerII} shows that the
distribution $\Pr\left[  R>\frac{t}{x^{1/a}}\right]  =1-F_{R}\left(  \frac
{t}{x^{\frac{1}{\alpha}}}\right)  $ has a Laplace transform equal to
$E_{a}\left(  -st^{a}\right)  =\sum_{k=0}^{\infty}\frac{\left(  -s\right)
^{k}t^{ak}}{\Gamma\left(  1+ka\right)  }$. We present a direct computation.
Partial integration of the pgf $\varphi_{X}(z)=\int_{0}^{\infty}e^{-zt}%
f_{X}\left(  t\right)  dt$,%
\[
\int_{0}^{\infty}e^{-zt}\left(  1-F_{X}\left(  t\right)  \right)
dt=\frac{1-\varphi_{X}(z)}{z}%
\]
is transformed, after letting $z\rightarrow\beta z$ for $\beta>0$, and
substituting $u=\beta t$, into%
\[
\int_{0}^{\infty}e^{-zu}\left(  1-F_{X}\left(  \frac{u}{\beta}\right)
\right)  du=\frac{1}{z}\left(  1-\varphi_{X}(\beta z)\right)
\]
Applied to the stable random variable $R$ with pgf $\varphi_{R}(z)=e^{-z^{a}}%
$,%
\[
\int_{0}^{\infty}e^{-zu}\left(  1-F_{R}\left(  \frac{u}{\beta}\right)
\right)  du=\frac{1}{z}\left(  1-e^{-\beta^{a}z^{a}}\right)
\]
and choosing $x=\beta^{a}>0$ yields%
\[
\int_{0}^{\infty}e^{-zu}\left(  1-F_{R}\left(  \frac{u}{x^{\frac{1}{\alpha}}%
}\right)  \right)  du=\frac{1}{z}\left(  1-e^{-xz^{a}}\right)
\]
Taking the Laplace transform of both sides with respect to $x$,%
\begin{align*}
\int_{0}^{\infty}dxe^{-sx}\int_{0}^{\infty}due^{-zu}\left(  1-F_{R}\left(
\frac{u}{x^{\frac{1}{\alpha}}}\right)  \right)   &  =\frac{1}{z}\int
_{0}^{\infty}dxe^{-sx}\left(  1-e^{-xz^{a}}\right) \\
&  =\frac{1}{z}\left(  \frac{1}{s}-\frac{1}{s+z^{a}}\right)  =\frac{1}%
{s}\left(  \frac{z^{a-1}}{s+z^{a}}\right)
\end{align*}
Introducing the Laplace transform (\ref{Laplace_transform_positive_argument})
with $b=1$ and $x=-s$%
\[
\int_{0}^{\infty}e^{-zt}E_{a}\left(  -st^{a}\right)  dt=\frac{z^{a-1}}%
{z^{a}+s}%
\]
and interchanging the integrals on the left-hand side, allowed by absolute
convergence,%
\[
\int_{0}^{\infty}dte^{-zt}\int_{0}^{\infty}dxe^{-sx}\left(  1-F_{R}\left(
\frac{t}{x^{\frac{1}{\alpha}}}\right)  \right)  =\frac{1}{s}\int_{0}^{\infty
}e^{-zt}E_{a}\left(  -st^{a}\right)  dt
\]
finally leads to%
\[
\int_{0}^{\infty}e^{-sx}\left(  1-F_{R}\left(  \frac{t}{x^{\frac{1}{\alpha}}%
}\right)  \right)  dx=\frac{1}{s}E_{a}\left(  -st^{a}\right)
\]
which is, however, a factor $\frac{1}{s}$ different from Feller's \cite[p.
453]{FellerII} result above\footnote{In the limit $s\rightarrow0$, the
right-hand side diverges and the left-hand side is $\int_{0}^{\infty}\left(
1-F_{R}\left(  \frac{t}{x^{\frac{1}{\alpha}}}\right)  \right)  dx=at^{a}%
\int_{0}^{\infty}\left(  1-F_{R}\left(  u\right)  \right)  u^{-a-1}%
du=at^{a}\int_{0}^{\infty}\Pr\left[  R>u\right]  u^{-a-1}du$ illustrating that
the integrand at the origin is $O\left(  u^{-a-1}\right)  $, leading to a
diverging integral.}.

\section{The integral $I_{a,b}\left(  z\right)  $}

\label{sec_integral_Iab(z)}We will study properties of the integral
$I_{a,b}\left(  z\right)  $ in (\ref{def_integral_Ia,b(z)}), whose integrand
$f\left(  u\right)  =\frac{z^{u}}{\Gamma\left(  b+au\right)  }$ attains a
maximum around $u\approx\frac{z^{\frac{1}{a}}-b}{a}$ for $b>0$. Indeed, the
derivative $f^{\prime}\left(  u\right)  =\frac{z^{u}}{\Gamma\left(
b+au\right)  }\left(  \log z-a\psi\left(  b+au\right)  \right)  $ vanishes
when $\log z=a\psi\left(  b+au\right)  $ and the expression
(\ref{Digamma_Binets_integral}) for the digamma function indicates, for large
$z$, that $\psi\left(  z\right)  \approx\log z$. For negative real $b$, on the
other hand, there may exist more than one extremum.

\medskip\refstepcounter{article}{\noindent\textbf{\thearticle. }%
}\ignorespaces\label{art_Ia,b(z)_for_z_complex} \emph{Complex argument }$z$.
For $z=re^{i\theta}$ with $\theta\in\left[  -\pi,\pi\right]  $, the integral
in (\ref{def_integral_Ia,b(z)}),
\[
I_{a,b}\left(  re^{i\theta}\right)  =\int_{0}^{\infty}\frac{r^{u}e^{i\theta
u}}{\Gamma\left(  b+au\right)  }\,du
\]
is split up in real and imaginary part, assuming that $a$ and $b$ are real, as%
\[
I_{a,b}\left(  re^{i\theta}\right)  =\int_{0}^{\infty}\frac{r^{u}\cos\theta
u}{\Gamma\left(  b+au\right)  }\,du+i\int_{0}^{\infty}\frac{r^{u}\sin\theta
u}{\Gamma\left(  b+au\right)  }\,du
\]
which illustrates that $I_{a,b}\left(  re^{i\theta}\right)  $ is only real if
$\theta=0$, i.e. only when $z$ is a real non-negative number. The point $z=0$
is a singularity as shown in \textbf{art}. \ref{art_integral_Ia,b(z)_series}
below. For $r>0$, both integrals exists and $I_{a,b}\left(  z\right)  $ is
thus defined along the negative real axis, where $z^{u}$ has a branch cut.
Both integrals decrease in $\theta$ for positive $a$ and $b$, and demonstrate
that $\left\vert I_{a,b}\left(  re^{i\theta}\right)  \right\vert $ decreases
with $\theta$. In contrast to the definition
(\ref{Mittag_Leffler_function_Eab}) of the Mittag-Leffler function
$E_{a,b}\left(  z\right)  $,%
\[
E_{a,b}\left(  re^{i\theta}\right)  =\sum_{k=0}^{\infty}\frac{r^{k}%
e^{ik\theta}}{\Gamma\left(  b+ak\right)  }=\sum_{k=0}^{\infty}\frac{r^{k}\cos
k\theta}{\Gamma\left(  b+ak\right)  }+i\sum_{k=0}^{\infty}\frac{r^{k}\sin
k\theta}{\Gamma\left(  b+ak\right)  }%
\]
demonstrates that $E_{a,b}\left(  re^{i\theta}\right)  $ is real along the
entire real axis, because $\sin\theta=0$ for $\theta=m\pi$ for $m\in
\mathbb{Z}$.

Moreover, for $\theta\neq0$, after substitution of $t=\theta u$
\[
I_{a,b}\left(  re^{i\theta}\right)  =\frac{1}{\theta}\int_{0}^{\infty}%
\frac{\left(  r^{\frac{1}{\theta}}\right)  ^{t}\cos t}{\Gamma\left(
b+\frac{a}{\theta}t\right)  }\,dt+\frac{i}{\theta}\int_{0}^{\infty}%
\frac{\left(  r^{\frac{1}{\theta}}\right)  ^{t}\sin t}{\Gamma\left(
b+\frac{a}{\theta}t\right)  }\,dt
\]
we observe that%
\[
I_{a,b}\left(  re^{i\theta}\right)  =\frac{1}{\theta}I_{\frac{a}{\theta}%
,b}\left(  r^{\frac{1}{\theta}}e^{i}\right)
\]
Thus, the integral at any complex number $z=re^{i\theta}$ with for $\theta
\neq0$ can be mapped to an evaluation along the straight line with angle equal
to 1 radius. $\allowbreak$

\medskip\refstepcounter{article}{\noindent\textbf{\thearticle. }%
}\ignorespaces\label{art_scaling_integral_Ia,b(z)} \emph{Functional relations
of the integral }$I_{a,b}\left(  z\right)  $. In contrast to the
Mittag-Leffler function $E_{a,b}\left(  z\right)  $, we can scale the integral
$I_{a,b}\left(  z\right)  $ easily by considering various (real) substitutions
in (\ref{def_integral_Ia,b(z)}). We consider a linear transformation
$b+au=\beta+\alpha v$ with $a>0$ and $\alpha>0$. Thus, let $u=\frac{\beta
-b}{a}+\frac{\alpha}{a}v$ in (\ref{def_integral_Ia,b(z)}), then
\begin{align*}
I_{a,b}\left(  z\right)   &  =\int_{0}^{\infty}\frac{z^{u}}{\Gamma\left(
b+au\right)  }\,du=\frac{\alpha}{a}z^{\frac{\beta-b}{a}}\int_{\frac{b-\beta
}{\alpha}}^{\infty}\frac{z^{\frac{\alpha}{a}v}}{\Gamma\left(  \beta+\alpha
v\right)  }\,dv\\
&  =\frac{\alpha}{a}z^{\frac{\beta-b}{a}}\left(  \int_{0}^{\infty}%
\frac{\left(  z^{\frac{\alpha}{a}}\right)  ^{v}}{\Gamma\left(  \beta+\alpha
v\right)  }\,dv-\int_{0}^{\frac{b-\beta}{\alpha}}\frac{z^{\frac{\alpha}{a}v}%
}{\Gamma\left(  \beta+\alpha v\right)  }\,dv\right)
\end{align*}
Using the definition (\ref{def_integral_Ia,b(z)}) leads to
\begin{equation}
I_{a,b}\left(  z\right)  =\frac{\alpha}{a}z^{\frac{\beta-b}{a}}\left(
I_{\alpha,\beta}\left(  z^{\frac{\alpha}{a}}\right)  -\int_{0}^{\frac{b-\beta
}{\alpha}}\frac{z^{\frac{\alpha}{a}v}}{\Gamma\left(  \beta+\alpha v\right)
}\,dv\right)  \label{linear_map_integral_Ia,b(z)}%
\end{equation}
from which%
\[
\alpha z^{\beta}I_{\alpha,\beta}\left(  z^{\alpha}\right)  -az^{b}%
I_{a,b}\left(  z^{a}\right)  =\alpha z^{\beta}\int_{0}^{\frac{b-\beta}{\alpha
}}\frac{\left(  z^{\alpha}\right)  ^{v}}{\Gamma\left(  \beta+\alpha v\right)
}\,dv
\]
Differentiating $z^{-\frac{\beta-b}{a}}I_{a,b}\left(  z\right)  =\frac{\alpha
}{a}\int_{\frac{b-\beta}{\alpha}}^{\infty}\frac{z^{\frac{\alpha}{a}v}}%
{\Gamma\left(  \beta+\alpha v\right)  }\,dv$ with respect to $b$ yields
$\frac{d}{db}\left(  z^{\frac{b-\beta}{a}}I_{a,b}\left(  z\right)  \right)
=-\frac{1}{a}\frac{z^{\frac{b-\beta}{a}}}{\Gamma\left(  b\right)  }$, which is
independent of $\alpha$ and leads to the linear differential equation in $b$,%
\[
\frac{dI_{a,b}\left(  z\right)  }{db}+\frac{\log z}{a}I_{a,b}\left(  z\right)
=-\frac{1}{a\Gamma\left(  b\right)  }\,
\]

The situation simplifies considerably when we choose $\beta=b$ in
(\ref{linear_map_integral_Ia,b(z)}),%
\[
I_{a,b}\left(  z\right)  =\frac{\alpha}{a}I_{\alpha,b}\left(  z^{\frac{\alpha
}{a}}\right)
\]
In particular, for $a>0$ and $\alpha=\frac{1}{a}$, the map%
\begin{equation}
I_{\frac{1}{a},b}\left(  z\right)  =a^{2}I_{a,b}\left(  z^{a^{2}}\right)
\label{Integral_Ia,b_map_1/a_to_a}%
\end{equation}
suggests a similar relation between $E_{\frac{1}{a},b}\left(  z\right)  $ and
$E_{a,b}\left(  z\right)  $ based on (\ref{Euler_Maclaurin_Ea,b_met_vrije_m})
and Bieberbach's deductions in \textbf{art.}
\ref{art_deductions_Mittag_Leffler_contour_integral}.

For $\beta=0$ and $\alpha=1$ in (\ref{linear_map_integral_Ia,b(z)}) and
$\alpha=\beta=1$, we have respectively
\[
I_{a,b}\left(  z\right)  =\frac{z^{-\frac{b}{a}}}{a}\left(  I_{1,0}\left(
z^{\frac{1}{a}}\right)  -\int_{0}^{b}\frac{\left(  z^{\frac{1}{a}}\right)
^{v}}{\Gamma\left(  v\right)  }\,dv\right)
\]
and%
\[
I_{a,b}\left(  z\right)  =\frac{1}{a}z^{\frac{1-b}{a}}\left(  I_{1,1}\left(
z^{\frac{1}{a}}\right)  -\int_{0}^{b-1}\frac{\left(  z^{\frac{1}{a}}\right)
^{v}}{v\Gamma\left(  v\right)  }\,dv\right)
\]
which bears resemblance to the last expression for $E_{a,b}\left(  z\right)  $
shown in \textbf{art.} \ref{art_fractional_a}. If $b>\beta>0$ and $z$ is real
and positive, then the integral in (\ref{linear_map_integral_Ia,b(z)}) is
non-negative and
\[
I_{a,b}\left(  z\right)  <\frac{\alpha}{a}z^{\frac{\beta-b}{a}}I_{\alpha
,\beta}\left(  z^{\frac{\alpha}{a}}\right)
\]

\medskip\refstepcounter{article}{\noindent\textbf{\thearticle. }%
}\ignorespaces\label{art_integral_Ia,b(z)_recursion_dvgl} \emph{Differential
recursion. }Differentiating the integral $I_{a,b}\left(  z\right)  $ in
(\ref{def_integral_Ia,b(z)}) and using the functional equation
(\ref{func_eq_Gamma}) of the Gamma function,%
\[
\frac{d}{dz}z^{b-1}I_{a,b}\left(  z^{a}\right)  =\int_{0}^{\infty}%
\frac{\left(  au+b-1\right)  z^{au+b-2}}{\Gamma\left(  b+au\right)  }%
\,du=\int_{0}^{\infty}\frac{z^{au+b-2}}{\Gamma\left(  b-1+au\right)  }\,du
\]
leads to a recursion equation in%
\begin{equation}
\frac{d}{dz}\left\{  z^{b-1}I_{a,b}\left(  z^{a}\right)  \right\}
=z^{b-2}I_{a,b-1}\left(  z^{a}\right)  \label{D_(z^(b-1)I_a,b(z^a))}%
\end{equation}
which is precisely the same as for $E_{a,b}\left(  .\right)  $ in
(\ref{D_(z^(b-1)E_a,b(z^a))}).

\medskip\refstepcounter{article}{\noindent\textbf{\thearticle.}}%
\ignorespaces\label{art_integral_Ia,b(z)_complex_integral} \emph{A complex
integral representation for }$I_{a,b}\left(  z\right)  $. We start by
concentrating on the integral%
\[
z^{b-1}I_{a,b}\left(  z^{a}\right)  =\int_{0}^{\infty}\frac{z^{b-1+au}}%
{\Gamma\left(  b+au\right)  }\,du
\]
whose Laplace transform is%
\[
L_{a,b}\left(  s\right)  =\int_{0}^{\infty}\left\{  z^{b-1}I_{a,b}\left(
z^{a}\right)  \right\}  e^{-zs}dz=\int_{0}^{\infty}e^{-zs}\int_{0}^{\infty
}\frac{z^{b-1+au}}{\Gamma\left(  b+au\right)  }du\,dz
\]
The reversal of the integrals is justified by absolute convergence,%
\begin{align*}
L_{a,b}\left(  s\right)   &  =\int_{0}^{\infty}\frac{du}{\Gamma\left(
b+au\right)  }\int_{0}^{\infty}z^{b+au-1}e^{-zs}dz=\int_{0}^{\infty}\frac
{du}{\Gamma\left(  b+au\right)  }\frac{\Gamma\left(  b+au\right)  }{s^{b+au}%
}\\
&  =s^{-b}\int_{0}^{\infty}e^{-au\log s}du=\frac{s^{-b}}{a\log s}%
\end{align*}
The inverse Laplace transform (\ref{inverse_Laplace}) returns a complex
integral,%
\begin{equation}
z^{b-1}I_{a,b}\left(  z^{a}\right)  =\frac{1}{2\pi ia}\int_{c-i\infty
}^{c+i\infty}\frac{e^{zs}}{s^{b}\log s}ds \label{I_ab_inverse_Laplace}%
\end{equation}
where $c>1$ because $L_{a,b}\left(  s\right)  $ is only analytic for
$\operatorname{Re}\left(  s\right)  >1$ due to the pole at $s=1$. Indeed, the
Taylor expansion of $\ln s$ around $s=1$, valid for $0<s<2$,%
\[
\ln s=\ln\left(  1+\left(  s-1\right)  \right)  =\sum_{k=0}^{\infty}%
\frac{\left(  -1\right)  ^{k}\left(  s-1\right)  ^{k+1}}{k+1}=\left(
s-1\right)  \left(  1+O\left(  s-1\right)  \right)
\]
demonstrates that $L_{a,b}\left(  s\right)  =\frac{s^{-b}}{a\log s}$ has a
pole at $s=1$ with residue 1, while $\frac{dL_{a,b}\left(  s\right)  }%
{db}=\frac{-s^{-b}}{a}$ only has the branch cut. After replacing $z^{a}$ by
$z$ in (\ref{I_ab_inverse_Laplace}) and combining with the definition
(\ref{def_integral_Ia,b(z)}), we obtain%
\begin{equation}
I_{a,b}\left(  z\right)  =\int_{0}^{\infty}\frac{e^{u\log z}}{\Gamma\left(
b+au\right)  }\,du=\frac{z^{\frac{1-b}{a}}}{2\pi ia}\int_{c-i\infty
}^{c+i\infty}\frac{e^{z^{\frac{1}{a}}s}}{s^{b}\log s}ds\hspace{1cm}c>1
\label{I_a,b(z)_complex_integral}%
\end{equation}

For $\operatorname{Re}\left(  z^{\frac{1}{a}}\right)  <0$, the contour in
(\ref{I_a,b(z)_complex_integral}) can be closed over the positive
$\operatorname{Re}\left(  s\right)  $-plane, in which the integrand is
analytic and\footnote{Hence, for $\operatorname{Re}\left(  z^{\frac{1}{a}%
}\right)  =\operatorname{Re}\left(  e^{\frac{1}{a}\log z}\right)  <0$, it
holds that%
\[
0=\int_{0}^{\infty}\frac{e^{u\log r}e^{iu\theta}}{\Gamma\left(  b+au\right)
}\,du=\frac{r^{-\frac{b}{a}}}{a}\int_{b}^{\infty}\frac{e^{w\frac{\log r}{a}%
}e^{i\left(  w-b\right)  \frac{\theta}{a}}}{\Gamma\left(  w\right)
}\,du=\frac{r^{-\frac{b}{a}}}{a}\int_{b}^{\infty}\frac{e^{w\frac{\log r}{a}%
}\left\{  \cos\left(  \frac{\theta}{a}\left(  w-b\right)  \right)
+i\sin\left(  \frac{\theta}{a}\left(  w-b\right)  \right)  \right\}  }%
{\Gamma\left(  w\right)  }\,du
\]
implying that%
\[
\int_{b}^{\infty}\frac{e^{w\frac{\log r}{a}}\cos\left(  \frac{\theta}%
{a}\left(  w-b\right)  \right)  }{\Gamma\left(  w\right)  }\,du=\int
_{b}^{\infty}\frac{e^{w\frac{\log r}{a}}\sin\left(  \frac{\theta}{a}\left(
w-b\right)  \right)  }{\Gamma\left(  w\right)  }\,du=0
\]
} $I_{a,b}\left(  z\right)  =0$. Let $z=re^{i\theta}$ with $\theta=\arg
z\in\left[  -\pi,\pi\right]  $ and recalling that $a>0$, then
$\operatorname{Re}\left(  z^{\frac{1}{a}}\right)  =r^{\frac{1}{a}}\cos
\frac{\theta}{a}$ and $\operatorname{Re}\left(  z^{\frac{1}{a}}\right)  <0$
requires that $\cos\frac{\theta}{a}<0$, which is equivalent to $\frac{\pi}%
{2}<\frac{\pm\theta}{a}<\pi$ or $\frac{\pi a}{2}<\left\vert \arg z\right\vert
<\pi a$. The latter condition, combined with $0<\left\vert \arg z\right\vert
<\pi$ is only possible if $0<a<2$. If $1\leq a<2$, then the combined condition
means that $\frac{\pi}{2}<\left\vert \arg z\right\vert <\pi$ or that
$\operatorname{Re}\left(  z\right)  <0$. Only if $0<a<\frac{1}{2}$, then the
combined condition means that $\frac{\pi}{2}<\left\vert \arg z\right\vert
<\frac{\pi a}{2}$ or that $\operatorname{Re}\left(  z\right)  >0$, while for
$\frac{1}{2}\leq a\leq1$, $\operatorname{Re}\left(  z\right)  $ can be either sign.

For $\operatorname{Re}\left(  z^{\frac{1}{a}}\right)  >0$, we close the
contour in (\ref{I_a,b(z)_complex_integral}) over the negative
$\operatorname{Re}\left(  s\right)  $-plane around the branch cut of $s^{b}%
\ln\left(  s\right)  $, which is the negative real axis. Thus, we consider the
contour $C$ that consists of the line at $c>1$, the quarter of a circle with
infinite radius from $\frac{\pi}{2}$ to $\pi-\varepsilon$, the line segment
above the real negative axis from minus infinity to $s=0$, the circle around
the origin $s=0$ from $\pi-\varepsilon$ back to $-\pi-\varepsilon$ with radius
$\delta$, the line segment below the real negative axis from $s=0$ towards
minus infinity, the quarter circle with infinite radius back to close the
contour $C$. This contour encloses the pole at $s=1$, whose residue is
$\lim_{s\rightarrow1}\frac{e^{z^{\frac{1}{a}}s}\left(  s-1\right)  }{s^{b}\ln
s}=e^{z^{\frac{1}{a}}}$. Cauchy's Residue Theorem \cite{Titchmarshfunctions}
results in%
\[
\frac{1}{2\pi i}\int_{C}\frac{e^{z^{\frac{1}{a}}s}}{s^{b}\ln s}ds=e^{z^{\frac
{1}{a}}}%
\]
while the evaluation of the contour $C$ yields%
\[
\frac{1}{2\pi i}\int_{C}\frac{e^{z^{\frac{1}{a}}s}}{s^{b}\ln s}ds=az^{\frac
{b-1}{a}}I_{a,b}\left(  z\right)  +\frac{1}{2\pi i}\int_{\infty}^{0}%
\frac{e^{-z^{\frac{1}{a}}x}d\left(  -x\right)  }{x^{b}e^{ib\left(
\pi-\varepsilon\right)  }\ln\left(  xe^{i\left(  \pi-\varepsilon\right)
}\right)  }+\frac{1}{2\pi i}\int_{0}^{\infty}\frac{e^{-z^{\frac{1}{a}}%
x}d\left(  -x\right)  }{x^{b}e^{ib\left(  -\pi-\varepsilon\right)  }\ln\left(
xe^{i\left(  -\pi-\varepsilon\right)  }\right)  }%
\]
since the parts of $C$ along the circles vanish for $\operatorname{Re}\left(
z\right)  >0$, but for $\delta\rightarrow0$ only provided $\operatorname{Re}%
\left(  b\right)  \leq1$. Hence, we obtain%
\[
az^{\frac{b-1}{a}}I_{a,b}\left(  z\right)  =e^{z^{\frac{1}{a}}}+\frac{1}{2\pi
i}\int_{0}^{\infty}\frac{e^{-z^{\frac{1}{a}}x}}{x^{b}}\left(  \frac{e^{ib\pi}%
}{\ln x-i\pi}-\frac{e^{-ib\pi}}{\ln x+i\pi}\right)  dx
\]
Finally, with $\frac{e^{ib\pi}}{\ln x-i\pi}-\frac{e^{-ib\pi}}{\ln x+i\pi}=2\pi
i\frac{\frac{\sin b\pi}{\pi}\ln x+\cos b\pi}{\left(  \pi^{2}+\left(  \ln
x\right)  ^{2}\right)  }$ and the definition (\ref{def_integral_Ia,b(z)}), we
arrive, for $\operatorname{Re}\left(  b\right)  \leq1$ and $\operatorname{Re}%
\left(  z^{\frac{1}{a}}\right)  >0$, at%
\begin{equation}
I_{a,b}\left(  z\right)  =\int_{0}^{\infty}\frac{z^{u}}{\Gamma\left(
b+au\right)  }\,du=\frac{z^{\frac{1-b}{a}}}{a}\left\{  e^{z^{\frac{1}{a}}%
}+\int_{0}^{\infty}\frac{e^{-z^{\frac{1}{a}}x}}{x^{b}}\left(  \frac{\frac{\sin
b\pi}{\pi}\ln x+\cos b\pi}{\pi^{2}+\left(  \ln x\right)  ^{2}}\right)
dx\right\}  \label{I_a,b(z)_integral_Re(z^(1/a))>0}%
\end{equation}
while, for $\operatorname{Re}\left(  z^{\frac{1}{a}}\right)  <0$,
$I_{a,b}\left(  z\right)  =0$. Expression
(\ref{I_a,b(z)_integral_Re(z^(1/a))>0}) is obvious related to Bierberbach
integral in \textbf{art}. \ref{art_deductions_Mittag_Leffler_contour_integral}%
, written as
\begin{equation}
E_{a,b}\left(  z\right)  =\frac{z^{\frac{1-b}{a}}}{a}\left\{  e^{z^{\frac
{1}{a}}}+\frac{1}{2\pi i}\frac{1}{z}\int_{C_{a\phi}^{\prime}}\frac{\left(
\frac{t}{z}\right)  ^{\frac{1-b}{a}}e^{t^{\frac{1}{a}}}}{\left(  \frac{t}%
{z}\right)  -1}dt\right\}  \label{E_a,b(z)_Bieberbach_growth_z}%
\end{equation}

For $b=1$, (\ref{I_a,b(z)_integral_Re(z^(1/a))>0}) simplifies to%
\[
I_{a,1}\left(  z\right)  =\frac{1}{a}\left\{  e^{z^{\frac{1}{a}}}-\int
_{0}^{\infty}\frac{e^{-z^{\frac{1}{a}}x}}{x\left(  \pi^{2}+\left(  \ln
x\right)  ^{2}\right)  }dx\right\}
\]
where $\int_{0}^{\infty}\frac{e^{-\lambda x}}{x\left(  \pi^{2}+\left(  \ln
x\right)  ^{2}\right)  }dx$ is increasing\footnote{Indeed (see \cite[p.73]%
{PVM_PAComplexNetsCUP}),%
\[
\int_{0}^{\infty}\frac{e^{-\lambda x}}{x\left(  \pi^{2}+\left(  \ln x\right)
^{2}\right)  }dx=\int_{-\infty}^{\infty}\frac{e^{-\lambda e^{t}}}{\pi
^{2}+t^{2}}dt\leq\int_{-\infty}^{\infty}\frac{dt}{\pi^{2}+t^{2}}=1
\]
} in $\operatorname{Re}\left(  \lambda\right)  >0$ from 0 to 1. Hence, for
$\operatorname{Re}\left(  z^{\frac{1}{a}}\right)  >0$, the following lower and
upper bound hold,
\[
\frac{1}{a}e^{z^{\frac{1}{a}}}-\frac{1}{a}<I_{a,1}\left(  z\right)  <\frac
{1}{a}e^{z^{\frac{1}{a}}}%
\]

With $\alpha=\beta=1$ in (\ref{linear_map_integral_Ia,b(z)}) and
$I_{1,1}\left(  z^{\frac{1}{a}}\right)  $ from
(\ref{I_a,b(z)_integral_Re(z^(1/a))>0}), we find for $\operatorname{Re}\left(
z^{\frac{1}{a}}\right)  >0$,%
\[
I_{a,b}\left(  z\right)  =\frac{z^{\frac{1-b}{a}}}{a}\left(  e^{z^{\frac{1}%
{a}}}-\int_{0}^{\infty}\frac{e^{-z^{\frac{1}{a}}x}}{x\left(  \pi^{2}+\left(
\ln x\right)  ^{2}\right)  }dx-\int_{0}^{b-1}\frac{\left(  z^{\frac{1}{a}%
}\right)  ^{x}}{x\Gamma\left(  v\right)  }\,dx\right)
\]
Comparison with (\ref{I_a,b(z)_integral_Re(z^(1/a))>0}) indicates, for
$b\leq1$ and $\operatorname{Re}\left(  z\right)  >0$,
that\footnote{Differentiating with respect to $b$,%
\[
-\frac{z^{b-1}}{\Gamma\left(  b\right)  }=\int_{0}^{\infty}\frac{e^{-zx}}%
{\pi^{2}+\left(  \ln x\right)  ^{2}}\frac{d}{db}\left(  x^{-b}\left(
\frac{\sin b\pi}{\pi}\ln x+\cos b\pi\right)  \right)  dx
\]
leads, using the reflection formula (\ref{Gamma_reflection_formula}), to an
identity%
\[
\frac{z^{b-1}}{\Gamma\left(  b\right)  }=\frac{\sin b\pi}{\pi}\int_{0}%
^{\infty}e^{-zx}x^{-b}dx=\frac{\sin b\pi}{\pi}\frac{\Gamma\left(  1-b\right)
}{z^{1-b}}%
\]
}%
\[
\int_{0}^{\infty}\frac{e^{-zx}}{x^{b}}\left(  \frac{\frac{\sin b\pi}{\pi}\ln
x+\cos b\pi}{\pi^{2}+\left(  \ln x\right)  ^{2}}\right)  dx=-\int_{0}^{\infty
}\frac{e^{-zx}}{x\left(  \pi^{2}+\left(  \ln x\right)  ^{2}\right)  }%
dx-\int_{0}^{b-1}\frac{z^{v}}{\Gamma\left(  v+1\right)  }\,dv
\]

\medskip\refstepcounter{article}{\noindent\textbf{\thearticle.}}%
\ignorespaces\label{art_integral_Ia,b(z)_complex_2} \emph{Another complex
integral representation for }$I_{a,b}\left(  z\right)  $. Another complex
integral follows directly from Hankel's contour
(\ref{Hankel_contour_integral_1opGamma}) as%
\[
I_{a,b}\left(  z\right)  =\,\frac{1}{2\pi i}\int_{C}w^{-b}e^{w}dw\int
_{0}^{\infty}e^{-u\left(  a\log w-\log z\right)  }du
\]
Only if $\operatorname{Re}\left(  a\log w-\log z\right)  >0$, which is
equivalent to $\operatorname{Re}\left(  \log\frac{w^{a}}{z}\right)
=\log\left\vert \frac{w^{a}}{z}\right\vert >0$ and $\left\vert \frac{w^{a}}%
{z}\right\vert >1$, then%
\begin{equation}
I_{a,b}\left(  z\right)  =\,\frac{1}{2\pi i}\int_{C}\frac{w^{-b}e^{w}}{a\log
w-\log z}dw\hspace{1cm}\text{with }\left\vert w\right\vert >\left\vert
z\right\vert ^{\frac{1}{a}} \label{I_a,b(z)_complex_integral_Hankel}%
\end{equation}
where the constraint $\left\vert w\right\vert >\left\vert z\right\vert
^{\frac{1}{a}}$ requires to deform the contour $C$ (as explained in
\textbf{art}. \ref{art_Mittag_Leffler_contour_integral}). The contour integral
(\ref{I_a,b(z)_complex_integral_Hankel}) bears a resemblance to the basic
complex integral (\ref{E_a,b(z)_basic_contour_integral_Hankel}) for
$E_{a,b}\left(  z\right)  $, whereas the above integral
(\ref{I_a,b(z)_complex_integral}) is closer to the Mittag-Leffler integral
(\ref{Mittag_Leffler_Bieberbach_integral}), although a formal substitution
$s=\frac{w^{a}}{z}$ in (\ref{I_a,b(z)_complex_integral_Hankel}) leads to
\[
I_{a,b}\left(  z\right)  =\,\frac{z^{\frac{1-b}{a}}}{2\pi ia}\int_{C^{\prime}%
}\frac{s^{\frac{1-b}{a}-1}e^{\left(  zs\right)  ^{\frac{1}{a}}}}{\log s}ds
\]

\medskip\refstepcounter{article}{\noindent\textbf{\thearticle.}}%
\ignorespaces\label{art_integral_Ia,b(z)_series} \emph{A series for }%
$I_{a,b}\left(  z\right)  $. The Taylor series of $I_{a,b}\left(  z\right)  $
around $z=\zeta$ equals%
\[
I_{a,b}\left(  z\right)  =I_{a,b}\left(  \zeta\right)  +\sum_{k=0}^{\infty
}\frac{1}{k!}\left.  \frac{d^{k}I_{a,b}\left(  z\right)  }{dz^{k}}\right\vert
_{z=\zeta}\left(  z-\zeta\right)  ^{k}%
\]
where the derivative for real $a>0$%
\[
\left.  \frac{d^{k}I_{a,b}\left(  z\right)  }{dz^{k}}\right\vert _{z=\zeta
}=\int_{0}^{\infty}\frac{\left.  \frac{d^{k}}{dz^{k}}z^{u}\right\vert
_{z=\zeta}}{\Gamma\left(  b+au\right)  }\,du=\frac{1}{\zeta^{k}}\int
_{0}^{\infty}\frac{u\left(  u-1\right)  \ldots\left(  u-k+1\right)  \zeta^{u}%
}{\Gamma\left(  b+au\right)  }\,du
\]
converges for all $k$, except when $\zeta=0$. Hence, in contrast to
$E_{a,b}\left(  z\right)  $, the function $I_{a,b}\left(  z\right)  $ is not
entire and has an essential singularity at $z=0$, where none of the
derivatives exists.

We expand the integrand of the integral $I_{a,b}\left(  z\right)  $, defined
in (\ref{def_integral_Ia,b(z)}), in a Taylor series around $b$,%
\[
\frac{1}{\Gamma\left(  b+au\right)  }=\sum_{j=0}^{\infty}\frac{1}{j!}\left.
\frac{d^{j}}{dy^{j}}\frac{1}{\Gamma\left(  y\right)  }\right\vert
_{y=b}\left(  au\right)  ^{j}%
\]
which converges for all $u$ and $b$, and obtain
\[
I_{a,b}\left(  z\right)  =\int_{0}^{\infty}\frac{z^{u}}{\Gamma\left(
b+au\right)  }\,du=\,\sum_{j=0}^{\infty}\frac{a^{j}}{j!}\left.  \frac{d^{j}%
}{dy^{j}}\frac{1}{\Gamma\left(  y\right)  }\right\vert _{y=b}\int_{0}^{\infty
}u^{j}e^{-u\left(  -\log z\right)  }du
\]
Only if $\operatorname{Re}\left(  \log z\right)  <0$ or $0\leq\left\vert
z\right\vert <1$, then we arrive at%
\begin{equation}
I_{a,b}\left(  z\right)  =-\frac{1}{\log z}\,\sum_{j=0}^{\infty}\left.
\frac{d^{j}}{dy^{j}}\frac{1}{\Gamma\left(  y\right)  }\right\vert
_{y=b}\left(  -\frac{a}{\log z}\right)  ^{j} \label{I_a,b(z)_series}%
\end{equation}
but this series only converges when $\left\vert \frac{a}{\log z}\right\vert
<1$. Indeed, alternatively, after $p$-times repeated partial integration, we
obtain%
\[
I_{a,b}\left(  z\right)  =-\frac{1}{\log z}\,\sum_{j=0}^{p-1}\left.
\frac{d^{j}}{dy^{j}}\frac{1}{\Gamma\left(  y\right)  }\right\vert
_{y=b}\left(  -\frac{a}{\log z}\right)  ^{j}+\left(  -\frac{a}{\log z}\right)
^{p}\int_{0}^{\infty}\left.  \frac{d^{p}}{dy^{p}}\frac{1}{\Gamma\left(
y\right)  }\right\vert _{y=b+au}z^{u}\,du
\]
where the last integral exists for all $p$ and $z$. Thus, if $p\rightarrow
\infty$ and $\left\vert \frac{a}{\log z}\right\vert <1$, then repeated partial
integration again produces (\ref{I_a,b(z)_series}). The series
(\ref{I_a,b(z)_series}) indicates that $I_{a,b}\left(  0\right)  =0$.

\medskip\textbf{Acknowledgement} I am grateful to Professor A. Apelblat and
Professor F. Mainardi for their very stimulating discussions. After the
submission of the first version on arXiv:2005.13330, they have contacted me
and pointed me to earlier work.

{\footnotesize
\bibliographystyle{plain}
\bibliography{cac,math,misc,net,pvm,qth,tel}
}

\appendix

\section{The Gamma function $\Gamma\left(  z\right)  $}

\label{sec_theory_Gamma_function}We review properties of the Gamma function.
Besides the basic functions \cite[Chapter IX and X]{Hardy_pure_math} like the
exponential, logarithm, circular or trigonometric functions (as sinus,
cosinus, tangens, etc.), the Gamma function is the next important complex
function. Nearly all books on complex function theory
\cite{Erdelyi_v1,Evgrafov_1965,Markushevich,Sansone,Titchmarshfunctions,Whittaker_Watson}
treat the Gamma function.

The Gamma function $\Gamma\left(  z\right)  $ is an extension of the factorial
$n!=1.2.3\ldots n$ in the integers $n\geq1$ to complex numbers $z$. The
factorial obeys%
\[
n!=n(n-1)!
\]
which directly generalizes to the functional equation $\Pi\left(  z\right)
=z\Pi\left(  z-1\right)  $ with $n!=\Pi\left(  n\right)  $ and $\Pi\left(
0\right)  =1$, in the notation of Gauss in his truly impressive manuscript
\cite[p. 146]{Gauss1813}. Later in 1814, Legendre defined the Gamma function
by its current notation $\Gamma\left(  z\right)  =\Pi\left(  z-1\right)  $,
with $\Gamma\left(  1\right)  =1$, and the functional equation $\Pi\left(
z\right)  =z\Pi\left(  z-1\right)  $ translates to the Gamma function
$\Gamma\left(  z\right)  $ as%
\begin{equation}
\Gamma\left(  z+1\right)  =z\Gamma\left(  z\right)  \label{func_eq_Gamma}%
\end{equation}
The first step in the theory of the Gamma function consists of finding a
solution of the functional equation (\ref{func_eq_Gamma}). Euler has proposed
his famous integral
\begin{equation}
\Gamma\left(  z\right)  =\int_{0}^{\infty}e^{-t}t^{z-1}dt\hspace{1cm}\text{for
}\operatorname{Re}\left(  z\right)  >0 \label{Euler_integral_Gamma_functie}%
\end{equation}
Partial integration of (\ref{Euler_integral_Gamma_functie}) shows that Euler's
integral (\ref{Euler_integral_Gamma_functie}) obeys the functional equation
(\ref{func_eq_Gamma}) and $\Gamma\left(  1\right)  =\int_{0}^{\infty}%
e^{-t}dt=1$. However, since Euler's integral is only valid for
$\operatorname{Re}\left(  z\right)  >0$, other ingenious methods have been
devised that are valid for all complex numbers $z$.

In his beautiful book on the Gamma function \cite{NielsenChelsea}, Nielsen
cites the historic achievements and reviews most contributions before 1906.
Nielsen \cite{NielsenChelsea} starts his book with the functional equation
(\ref{func_eq_Gamma}) of the Gamma function and immediately remarks that any
solution can be multiplied by a periodic function $\omega\left(  z\right)
=\omega\left(  z+1\right)  $ with period $1$. Next, Nielsen
\cite{NielsenChelsea} concentrates on the digamma function, which is the
logarithmic derivative $\psi\left(  z\right)  =\frac{d}{dz}\log\Gamma\left(
z\right)  $ and the functional equation (\ref{func_eq_Gamma}) tells us that%
\begin{equation}
\psi\left(  z+1\right)  =\psi\left(  z\right)  +\frac{1}{z}
\label{func_eq_digamma}%
\end{equation}
After $n$ iterations,%
\[
\psi\left(  z\right)  =\psi\left(  z+n\right)  -\sum_{k=0}^{n-1}\frac{1}{z+k}%
\]
and taking the limit $n\rightarrow\infty$, we formally obtain $\psi\left(
z\right)  =\lim_{n\rightarrow\infty}\psi\left(  z+n\right)  -\lim
_{n\rightarrow\infty}\sum_{k=0}^{n}\frac{1}{z+k}$. However, the latter series
does not converge\footnote{Indeed,
\[
\sum_{k=0}^{n}\frac{1}{z+k}=\frac{1}{z}+\sum_{k=1}^{n_{0}}\frac{1}{k\left(
1+\frac{z}{k}\right)  }+\sum_{k=1+n_{0}}^{n}\frac{1}{k\left(  1+\frac{z}%
{k}\right)  }%
\]
We can choose $n_{0}>\left\vert z\right\vert $, so that $0<\left\vert
1+\frac{z}{k}\right\vert <2$ and $\left\vert \sum_{k=1+n_{0}}^{n}\frac
{1}{k\left(  1+\frac{z}{k}\right)  }\right\vert >\frac{1}{2}\sum_{k=1+n_{0}%
}^{n}\frac{1}{k}\rightarrow\infty$, because the harmonic series $\sum
_{k=1}^{n}\frac{1}{k}$ diverges for $n\rightarrow\infty$.}, implying that
$\lim_{n\rightarrow\infty}\psi\left(  z+n\right)  =\infty$. Nielsen then
applies Weierstrass's factorization theory for entire functions
\cite{Titchmarshfunctions} and deduces Weierstrass's product,%
\begin{equation}
\frac{1}{\Gamma(z+1)}=e^{\gamma\,z}\prod_{n=1}^{\infty}\left(  1+\frac{z}%
{n}\right)  \,e^{-z/n} \label{Weierstrass_product_Gamma_function}%
\end{equation}
which we will derive from Gauss's product (\ref{Gauss_product_Gamma_function})
in \textbf{art. }\ref{art_Weierstrass_product}. Weierstrass created his
magnificent theory for entire functions, a pearl of complex function theory,
inspired by Gauss's product (\ref{Gauss_product_Gamma_function}) and Gauss's
remark on factorization in \cite[p. 146]{Gauss1813}.

\subsection{Gauss's approach}

Iterating the functional equation (\ref{func_eq_Gamma}) $n$-times gives
$\Gamma\left(  z\right)  =\frac{\Gamma\left(  z+n\right)  }{z\left(
z+1\right)  \ldots\left(  z+n-1\right)  }$, but purely iterating
(\ref{func_eq_Gamma}) for non-integer values is not successful. Therefore,
Gauss \cite[p. 144]{Gauss1813} proposes to consider the more general form%
\begin{equation}
\Pi\left(  k,z\right)  =\frac{\Pi\left(  k\right)  \Pi\left(  z\right)  }%
{\Pi\left(  k+z\right)  }k^{z}=\frac{1}{\left(  z+1\right)  }\frac{2}{\left(
z+2\right)  }\frac{3}{\left(  z+3\right)  }\ldots\frac{k}{\left(  z+k\right)
}k^{z} \label{Gauss_finite_factoring_Gamma}%
\end{equation}
which satisfies%
\begin{equation}
\Pi\left(  k,z+1\right)  =\Pi\left(  k,z\right)  \frac{\left(  z+1\right)
}{\left(  1+\frac{z+1}{k}\right)  } \label{func_eq_Pi(k,z_in_z}%
\end{equation}
as well as%
\begin{equation}
\Pi\left(  k+1,z\right)  =\Pi\left(  k,z\right)  \frac{\left(  1+\frac{1}%
{k}\right)  ^{z+1}}{\left(  1+\frac{z+1}{k}\right)  }
\label{func_eq_Pi(k,z)_in_k}%
\end{equation}
Iterating (\ref{func_eq_Pi(k,z)_in_k}), with $\Pi\left(  1,z\right)  =\frac
{1}{z+1}$, results in%
\begin{equation}
\Pi\left(  k,z\right)  =\frac{1}{z+1}\frac{2^{z+1}}{\left(  2+z\right)  }%
\frac{3^{z+1}}{2^{z}\left(  3+z\right)  }\ldots\frac{k^{z+1}}{\left(
k-1\right)  ^{z}\left(  k+z\right)  }=\frac{1}{z+1}\prod_{n=2}^{k}%
\frac{n^{z+1}}{\left(  n-1\right)  ^{z}\left(  n+z\right)  }
\label{Pi(k,z)_iterated}%
\end{equation}
The interesting observation from $\lim_{k\rightarrow\infty}\frac{\left(
1+\frac{1}{k}\right)  ^{z+1}}{\left(  1+\frac{z+1}{k}\right)  }=1$ in
(\ref{func_eq_Pi(k,z)_in_k}) is that $\lim_{k\rightarrow\infty}\Pi\left(
k,z\right)  $ exist for all $z$, which Gauss demonstrates after taking the
logarithm of both sides, while the first functional equation
(\ref{func_eq_Pi(k,z_in_z}) indicates that $\lim_{k\rightarrow\infty}%
\Pi\left(  k,z\right)  $ satisfies the functional equation $\Pi\left(
z+1\right)  =(z+1)\Pi\left(  z\right)  $. Combining both, Gauss is led to%
\begin{equation}
\Pi\left(  z\right)  =\lim_{k\rightarrow\infty}\Pi\left(  k,z\right)
=\lim_{k\rightarrow\infty}\frac{\Pi\left(  k\right)  \Pi\left(  z\right)
}{\Pi\left(  k+z\right)  }k^{z} \label{def_Gamma_Gauss}%
\end{equation}
which is equivalent with $\Gamma\left(  z+1\right)  =\Pi\left(  z\right)  $ to%
\[
\Gamma\left(  z\right)  =\lim_{k\rightarrow\infty}\frac{k!\Gamma\left(
z\right)  }{\Gamma\left(  k+1+z\right)  }k^{z}=\lim_{k\rightarrow\infty}%
\frac{k^{z}k!}{z\left(  z+1\right)  \ldots\left(  z+k\right)  }%
\]
Rewriting Gauss's definition in (\ref{def_Gamma_Gauss}) as $\Gamma\left(
z\right)  =\lim_{k\rightarrow\infty}\Pi\left(  k,z-1\right)  $ and introducing
(\ref{Pi(k,z)_iterated}) yields%
\[
\Gamma\left(  z\right)  =\lim_{k\rightarrow\infty}\frac{1}{z}\prod_{n=2}%
^{k}\frac{n^{z}}{\left(  n-1\right)  ^{z-1}\left(  n+z-1\right)  }%
=\lim_{k\rightarrow\infty}\frac{1}{z}\prod_{n=1}^{k-1}\frac{\left(
n+1\right)  ^{z}}{n^{z-1}\left(  n+z\right)  }%
\]
Finally\footnote{Gauss proceeds further in \cite[p. 148]{Gauss1813} and
derives the reflection formula from his classical result \cite[15.1.20]%
{Abramowitz} for the hypergeometric series at $z=1$, for $c\neq-k$ ($k$
integer) and Re$\left(  c-a-b\right)  >0$,%
\begin{equation}
F\left(  a,b;c;1\right)  =\frac{\Gamma(c)\Gamma(c-a-b)}{\Gamma(c-a)\Gamma
(c-b)}=\frac{\Gamma(c)}{\Gamma(a)\Gamma(b)}\sum_{n=0}^{\infty}\frac
{\Gamma(a+n)\Gamma(b+n)}{\Gamma(c+n)n!} \label{hypergeometric_z=1_Gauss}%
\end{equation}
\par
He also deduces his multiplication formula, compares his theory with Stirling
and Euler's logarithmic expansion in terms of Bernoulli numbers, studies the
digamma function, derives his fractional argument digamma function, deduces an
integral for the digamma function and complements Euler's computations. In
short, an amazing sequence of beautifully derived deep results that constitute
our current basis of the Gamma function. In line with his genius, Gauss even
laid the basis of prime factors of an entire function of which Weierstass has
given the functional theory \cite[Chapter VIII]{Titchmarshfunctions}.}, we
arrive at Gauss's infinite product\footnote{Both Klein \cite[p. 74]%
{Klein_1933} and Whittaker and Watson \cite[p. 237]{Whittaker_Watson} mention
that Euler has given (\ref{Gauss_product_Gamma_function}) in a letter to
Goldbach in 1729, but that Gauss has provided the first rigorous analysis in
\cite{Gauss1813}.} for the Gamma function%
\begin{equation}
\Gamma\left(  z\right)  =\frac{1}{z}\prod_{n=1}^{\infty}\left(  1+\frac{1}%
{n}\right)  ^{z}\left(  1+\frac{z}{n}\right)  ^{-1}
\label{Gauss_product_Gamma_function}%
\end{equation}
which converges for all complex $z$, except for the integers at
$z=0,-1,-2,\ldots$, at which $\Gamma\left(  z\right)  $ has simples poles. The
inverse of the product (\ref{Gauss_product_Gamma_function}) shows that
$\frac{1}{\Gamma\left(  z\right)  }$ is an entire function and that
$\Gamma\left(  z\right)  $ has no zeros in the finite complex plane.

Of course, the proposal of (\ref{Gauss_finite_factoring_Gamma}) by Gauss was
crucial towards his elegant product (\ref{Gauss_product_Gamma_function}).
Gauss posited (\ref{Gauss_finite_factoring_Gamma}) without providing
intuition. Perhaps the most convincing argument for Gauss's starting point
(\ref{Gauss_finite_factoring_Gamma}) is given by Klein \cite[p. 71]%
{Klein_1933}, who gives three definitions of the Gamma function, of which the
third is also discussed by Gauss himself \cite[p. 151]{Gauss1813}. Klein
\cite[p. 74]{Klein_1933} starts from the Beta-integral, studied by Euler and
valid for $\operatorname{Re}\left(  z\right)  >0$ and $\operatorname{Re}%
\left(  q\right)  >0$,%
\[
B\left(  z,q\right)  =\frac{\Gamma\left(  z\right)  \Gamma\left(  q\right)
}{\Gamma\left(  z+q\right)  }=\int_{0}^{1}u^{z-1}\left(  1-u\right)  ^{q-1}du
\]
for $q=k+1$ and makes the substitution $u=\frac{v}{k}$,%
\[
B\left(  z,k+1\right)  =\int_{0}^{1}u^{z-1}\left(  1-u\right)  ^{k}%
du=k^{-z}\int_{0}^{k}v^{z-1}\left(  1-\frac{v}{k}\right)  ^{k}dv
\]
Using $\lim_{k\rightarrow\infty}\left(  1-\frac{v}{k}\right)  ^{k}=e^{-v}$ and
Euler's integral (\ref{Euler_integral_Gamma_functie}), Klein \cite[p.
74]{Klein_1933} arrives for $\operatorname{Re}\left(  z\right)  >0$ at%
\[
\Gamma\left(  z\right)  =\int_{0}^{\infty}v^{z-1}e^{-v}dv=\lim_{k\rightarrow
\infty}\left(  k^{z}B\left(  z,k+1\right)  \right)  =\lim_{k\rightarrow\infty
}\left(  k^{z}\frac{\Gamma\left(  z\right)  \Gamma\left(  k+1\right)  }%
{\Gamma\left(  z+k+1\right)  }\right)
\]
which is Gauss's definition (\ref{def_Gamma_Gauss}). In contrast to Euler's
integral (\ref{Euler_integral_Gamma_functie}) for $\operatorname{Re}\left(
z\right)  >0$, the functional equation (\ref{func_eq_Pi(k,z)_in_k}) is valid
for all $z$ and so is Gauss's product (\ref{Gauss_product_Gamma_function}).

\subsection{Deductions from Gauss's product
(\ref{Gauss_product_Gamma_function}) for $\Gamma\left(  z\right)  $}

\medskip\refstepcounter{article}{\noindent\textbf{\thearticle. }%
}\ignorespaces\label{art_Weierstrass_product} \emph{Weierstrass's product}.
Weierstrass's product (\ref{Weierstrass_product_Gamma_function}) can be
obtained from Gauss's definition $\Gamma\left(  z+1\right)  =\lim
_{k\rightarrow\infty}\Pi\left(  k,z\right)  $. Following Erd\'{e}lyi \emph{et
al}. \cite[p. 2]{Erdelyi_v1}, Gauss's definition of the Gamma function can be
written, with (\ref{Gauss_finite_factoring_Gamma}), as%
\begin{align*}
\frac{1}{\Gamma\left(  z+1\right)  }  &  =\lim_{k\rightarrow\infty}\left(
1+z\right)  \left(  1+\frac{z}{2}\right)  \left(  1+\frac{z}{3}\right)
\ldots\left(  1+\frac{z}{k}\right)  e^{-z\log k}\\
&  =\lim_{k\rightarrow\infty}\left(  1+z\right)  e^{-z}\left(  1+\frac{z}%
{2}\right)  e^{-\frac{z}{2}}\left(  1+\frac{z}{3}\right)  e^{-\frac{z}{3}%
}\ldots\left(  1+\frac{z}{k}\right)  e^{-\frac{z}{k}}e^{z\left(  \sum
_{n=1}^{k}\frac{1}{n}-\log k\right)  }\\
&  =\lim_{k\rightarrow\infty}\prod_{n=1}^{k}\left(  1+\frac{z}{n}\right)
e^{-\frac{z}{n}}\lim_{k\rightarrow\infty}e^{z\left(  \sum_{n=1}^{k}\frac{1}%
{n}-\log k\right)  }%
\end{align*}
Introducing\footnote{Gauss \cite[154, footnote]{Gauss1813} gives the
Euler-Mascheroni constant $\gamma$ up to 40 decimals accurate. Gauss provided
the method (i.e. the power series of the digamma function $\psi\left(
z\right)  $) and Fredericus Bernhardus Gothofredus Nicolai has performed the
computation. Whittaker and Watson \cite[p. 235]{Whittaker_Watson} mention that
J. C. Adams computed $\gamma$ up to 260 decimals accurate.} Euler's constant
\cite[6.1.3]{Abramowitz}%
\begin{equation}
\gamma=\lim_{k\rightarrow\infty}\left(  \sum_{n=1}^{k}\frac{1}{n}-\log
k\right)  =0.57721\ldots\label{def_Euler_constant_gamma}%
\end{equation}
leads to Weierstrass's product (\ref{Weierstrass_product_Gamma_function})
which illustrates that Euler's constant $\gamma$ plays a fundamental role in
the theory of the Gamma function.

Whittacker and Watson \cite[p. 235]{Whittaker_Watson} elegantly demonstrate
that the limit in (\ref{def_Euler_constant_gamma}) exists. They define%
\[
u_{n}=\int_{0}^{1}\frac{t}{n\left(  n+t\right)  }dt=\frac{1}{n}-\log\frac
{n+1}{n}%
\]
With $\sum_{n=1}^{k-1}\log\frac{n+1}{n}=\sum_{n=1}^{k-1}\log\left(
n+1\right)  -\sum_{n=1}^{k-1}\log n=\log\left(  k\right)  $ and
\begin{align}
\sum_{n=1}^{k}\frac{1}{n}-\log k  &  =\sum_{n=1}^{k}\left(  \frac{1}{n}%
-\log\frac{n+1}{n}\right)  +\log\frac{k+1}{k} \label{gamma_finite_k_identity}%
\\
&  =\sum_{n=1}^{k}u_{n}+\log\left(  1+\frac{1}{k}\right) \nonumber
\end{align}
an alternative representation of Euler's constant
(\ref{def_Euler_constant_gamma}) is obtained as $\gamma=\sum_{n=1}^{\infty
}u_{n}$. Since $0<u_{n}=\int_{0}^{1}\frac{t}{n\left(  n+t\right)  }dt<\int
_{0}^{1}\frac{1}{n^{2}}dt=\frac{1}{n^{2}}$ and $\sum_{n=1}^{\infty}\frac
{1}{n^{2}}=\zeta\left(  2\right)  =\frac{\pi^{2}}{6}=1.64493$, it holds that
$0<\gamma<\frac{\pi^{2}}{6}$.

The bounds can be sharpened by the inequality $\frac{1}{n+1}<\int_{n}%
^{n+1}\frac{dx}{x}=\log\frac{n+1}{n}<\frac{1}{n}$ for any $n>0$. Indeed, since
$\frac{1}{n}-\log\frac{n+1}{n}>0$, the identity (\ref{gamma_finite_k_identity}%
) provides the lower bound%
\[
\sum_{n=1}^{k}\frac{1}{n}-\log k=1-\log2+\sum_{n=2}^{k}\left(  \frac{1}%
{n}-\log\frac{n+1}{n}\right)  +\log\frac{k+1}{k}>1-\log2=0.30683
\]
Similarly, rewriting the identity (\ref{gamma_finite_k_identity}) and using
$\log\frac{n+1}{n}-\frac{1}{n+1}>0$ gives us the upper bound
\begin{align*}
\sum_{n=1}^{k}\frac{1}{n}-\log k  &  =1+\sum_{n=2}^{k}\frac{1}{n}-\sum
_{n=1}^{k-1}\log\frac{n+1}{n}=1-\sum_{n=1}^{k-1}\left(  \log\frac{n+1}%
{n}-\frac{1}{n+1}\right) \\
&  <1-\left(  \log2-\frac{1}{2}\right)  =0.80683
\end{align*}
In summary, we find $1-\log2<\gamma<\frac{3}{2}-\log2$. Sharper bounds follows
from Poisson's integral (\ref{Digamma_Binets_integral}) in \textbf{art}.
\ref{art_Stirling_formula}.

\medskip\refstepcounter{article}{\noindent\textbf{\thearticle. }%
}\ignorespaces\label{art_reflection_Gamma} \emph{Reflection formula}. Gauss
\cite[p. 148]{Gauss1813} derives the reflection formula for the Gamma function
via his contiguous relations of the hypergeometric function, that, for
specific parameters, reduce to the sinus function. Gauss then finds the
infinite product of the sinus function \cite[4.3.89]{Abramowitz}%
\begin{equation}
\sin\left(  \pi x\right)  =\pi x%
{\displaystyle\prod\limits_{k=1}^{\infty}}
\left(  1-\frac{x^{2}}{k^{2}}\right)  \label{infinite_product_sinz}%
\end{equation}
Reversely, if we consider the infinite product (\ref{infinite_product_sinz}%
)\ as known\footnote{After integration of the Taylor series of $\pi\cot\left(
\pi x\right)  =\frac{1}{x}-2\sum_{n=1}^{\infty}\zeta\left(  2n\right)
\,x^{2n-1}$, valid for $\left\vert x\right\vert <1$, leading to $\log\frac{\pi
x}{\sin\left(  \pi x\right)  }=2\sum_{n=1}^{\infty}\frac{\zeta\left(
2n\right)  }{2n}\,x^{2n}$, in which the Zeta functions $\zeta\left(
2n\right)  =\sum_{k=1}^{\infty}\frac{1}{k^{2n}}$, one deduces that%
\[
\log\frac{\pi x}{\sin\left(  \pi x\right)  }=\sum_{k=1}^{\infty}\sum
_{n=1}^{\infty}\frac{1}{n}\left(  \frac{x}{k}\right)  ^{2n}=-\sum
_{k=1}^{\infty}\log\left(  1-\frac{x^{2}}{k^{2}}\right)
\]
from which (\ref{infinite_product_sinz}) follows. Although derived under the
restriction $\left\vert x\right\vert <1$, (\ref{infinite_product_sinz}) can be
shown to hold for any complex $x$.}, then it follows from Gauss's infinite
product (\ref{Gauss_product_Gamma_function}) for $\Gamma\left(  z\right)  $
that%
\[
\frac{1}{\Gamma\left(  z\right)  \Gamma\left(  -z\right)  }=-z^{2}\prod
_{n=1}^{\infty}\left(  1-\frac{z^{2}}{n^{2}}\right)  =-\frac{\sin\pi z}{\pi z}%
\]
which establishes the reflection formula \cite[6.1.17]{Abramowitz} of the
Gamma function, valid for all $z$,%
\begin{equation}
\Gamma\left(  z\right)  \Gamma\left(  1-z\right)  =\frac{\pi}{\sin\pi z}
\label{Gamma_reflection_formula}%
\end{equation}
For example, if $z=\frac{1}{2}$, then the reflection formula
(\ref{Gamma_reflection_formula}) shows that $\Gamma\left(  \frac{1}{2}\right)
=\sqrt{\pi}$.

\medskip\refstepcounter{article}{\noindent\textbf{\thearticle. }%
}\ignorespaces\label{art_multiplication_formula_Gamma} \emph{Multiplication
formula}. Gauss \cite[p. 149]{Gauss1813} derives his elegant\footnote{Gauss
writes \textquotedblleft Unde habemus theorema elegans\textquotedblright.}
multiplication formula \cite[6.1.20]{Abramowitz}
\begin{equation}
\Gamma\left(  nz\right)  =\left(  2\pi\right)  ^{\frac{1}{2}\left(
1-n\right)  }n^{nz-\frac{1}{2}}%
{\displaystyle\prod\limits_{k=0}^{n-1}}
\Gamma\left(  z+\frac{k}{n}\right)  \label{Gamma_multiplication_formula}%
\end{equation}
as follows. Gauss observes that%
\begin{equation}
\frac{n^{nz}\prod_{j=0}^{n-1}\Pi\left(  k,z-\frac{j}{n}\right)  }{\Pi\left(
nk,nz\right)  }=\frac{\left(  \Gamma\left(  k+1\right)  \right)  ^{n}}%
{\Gamma\left(  nk+1\right)  k^{\frac{\left(  n-1\right)  }{2}}}
\label{Gauss_product_Gauss_factors}%
\end{equation}
does not depend on $z$. However, he does not give the derivation of
(\ref{Gauss_product_Gauss_factors}), but we do. Multiplying the Gauss factors
(\ref{Gauss_finite_factoring_Gamma}),%
\begin{align*}
\prod_{j=0}^{n-1}\Pi\left(  k,z-\frac{j}{n}\right)   &  =\prod_{j=0}%
^{n-1}\frac{\Pi\left(  k\right)  \Pi\left(  z-\frac{j}{n}\right)  }{\Pi\left(
k+z-\frac{j}{n}\right)  }k^{z-\frac{j}{n}}\\
&  =\left(  \Pi\left(  k\right)  \right)  ^{n}k^{nz}k^{-\frac{1}{n}\sum
_{j=1}^{n-1}j}\prod_{j=0}^{n-1}\frac{\Pi\left(  z-\frac{j}{n}\right)  }%
{\Pi\left(  z+k-\frac{j}{n}\right)  }%
\end{align*}
and, changing the notation $\Pi\left(  z\right)  =\Gamma\left(  z+1\right)  $
and with $\sum_{j=1}^{n-1}j=\frac{n\left(  n-1\right)  }{2}$, we have
\[
\prod_{j=0}^{n-1}\Pi\left(  k,z-\frac{j}{n}\right)  =\left(  \Gamma\left(
k+1\right)  \right)  ^{n}k^{nz}k^{-\frac{\left(  n-1\right)  }{2}}\prod
_{j=0}^{n-1}\frac{\Gamma\left(  z+1-\frac{j}{n}\right)  }{\Gamma\left(
z+k+1-\frac{j}{n}\right)  }%
\]
Introducing Gauss's infinite product (\ref{Gauss_product_Gamma_function})
yields
\begin{align*}
\prod_{j=0}^{n-1}\frac{\Gamma\left(  z+1-\frac{j}{n}\right)  }{\Gamma\left(
z+k+1-\frac{j}{n}\right)  }  &  =\prod_{j=0}^{n-1}\prod_{m=1}^{\infty}%
\frac{m^{k}}{\left(  m+1\right)  ^{k}}\prod_{m=0}^{\infty}\frac{\left(
nm+nz+nk+n-j\right)  }{\left(  nm+nz+n-j\right)  }\\
&  =\prod_{m=1}^{\infty}\frac{m^{nk}}{\left(  m+1\right)  ^{nk}}\prod
_{m=0}^{\infty}\prod_{l=1}^{n}\frac{\left(  nm+nz+nk+l\right)  }{\left(
nm+nz+l\right)  }\\
&  =\prod_{m=1}^{\infty}\frac{m^{nk}}{\left(  m+1\right)  ^{nk}}\prod
_{m=0}^{\infty}\prod_{l=1}^{n}\left(  1+\frac{nk}{nm+nz+l}\right)
\end{align*}
Because $r=nm+l$ runs over all integers, we observe that%
\[
\prod_{m=0}^{\infty}\prod_{l=1}^{n}\left(  1+\frac{nk}{nm+l+nz}\right)
=\prod_{r=1}^{\infty}\left(  1+\frac{nk}{r+nz}\right)
\]
and we arrive at%
\[
\prod_{j=0}^{n-1}\Pi\left(  k,z-\frac{j}{n}\right)  =\left(  \Gamma\left(
k+1\right)  \right)  ^{n}k^{nz}k^{-\frac{\left(  n-1\right)  }{2}}\prod
_{m=1}^{\infty}\frac{m^{nk}}{\left(  m+1\right)  ^{nk}}\prod_{r=1}^{\infty
}\left(  1+\frac{nk}{r+nz}\right)
\]
Gauss divides $\prod_{j=0}^{n-1}\Pi\left(  k,z-\frac{j}{n}\right)  $ by
$\Pi\left(  nk,nz\right)  $, which equals
\begin{align*}
\Pi\left(  nk,nz\right)   &  =\frac{\Pi\left(  nk\right)  \Pi\left(
nz\right)  }{\Pi\left(  nk+nz\right)  }n^{nz}k^{nz}=\frac{\Gamma\left(
nk+1\right)  \Gamma\left(  nz+1\right)  }{\Gamma\left(  nk+nz+1\right)
}n^{nz}k^{nz}\\
&  =n^{nz}k^{nz}\Gamma\left(  nk+1\right)  \prod_{m=1}^{\infty}\frac{m^{nk}%
}{\left(  m+1\right)  ^{nk}}\prod_{m=1}^{\infty}\frac{m+nz+nk}{m+nz}\\
&  =n^{nz}k^{nz}\Gamma\left(  nk+1\right)  \prod_{m=1}^{\infty}\frac{m^{nk}%
}{\left(  m+1\right)  ^{nk}}\prod_{r=1}^{\infty}\left(  1+\frac{nk}%
{r+nz}\right)
\end{align*}
resulting in (\ref{Gauss_product_Gauss_factors}) and demonstrating that the
left-hand side of (\ref{Gauss_product_Gauss_factors}) is independent of $z$.
As follows from the Gauss factors (\ref{Gauss_finite_factoring_Gamma}),
$\Pi\left(  k,0\right)  =\Pi\left(  nk,0\right)  =1$ and the choice of $z=0$
in the left-hand side of (\ref{Gauss_product_Gauss_factors}) is%
\[
\frac{n^{nz}\prod_{j=0}^{n-1}\Pi\left(  k,z-\frac{j}{n}\right)  }{\Pi\left(
nk,nz\right)  }=\prod_{j=0}^{n-1}\Pi\left(  k,-\frac{j}{n}\right)
\]
After taking the limit of $k\rightarrow\infty$ of both sides, the definition
$\Pi\left(  z\right)  =\Gamma\left(  z+1\right)  =\lim_{k\rightarrow\infty}%
\Pi\left(  k,z\right)  $ leads to%
\[
\frac{n^{nz}\prod_{j=0}^{n-1}\Gamma\left(  z+1-\frac{j}{n}\right)  }%
{\Gamma\left(  nz+1\right)  }=\prod_{j=1}^{n-1}\Gamma\left(  1-\frac{j}%
{n}\right)
\]
Let $l=n-j$, then $\prod_{j=1}^{n-1}\Gamma\left(  1-\frac{j}{n}\right)
=\prod_{j=1}^{n-1}\Gamma\left(  \frac{n-j}{n}\right)  =\prod_{l=1}^{n-1}%
\Gamma\left(  \frac{l}{n}\right)  $ so that%
\[
\left(  \prod_{j=1}^{n-1}\Gamma\left(  1-\frac{j}{n}\right)  \right)
^{2}=\prod_{j=1}^{n-1}\Gamma\left(  1-\frac{j}{n}\right)  \Gamma\left(
\frac{j}{n}\right)  =\prod_{j=1}^{n-1}\frac{\pi}{\sin\pi\frac{j}{n}}=\frac
{\pi^{n-1}}{\prod_{j=1}^{n-1}\sin\frac{\pi j}{n}}%
\]
where the reflection formula (\ref{Gamma_reflection_formula}) has been
invoked. Using Euler's formula\footnote{The polynomial $z^{n}-1$ has as zeros
the $n$-th roots of unity, $z^{n}-1=\left(  z-1\right)  \prod_{k=1}%
^{n-1}\left(  1-z\,e^{-\frac{2\pi ik}{n}}\right)  $. Choosing $z=e^{\frac
{2ix}{n}}$ leads, after some manipulations, to
(\ref{sinx_as_product_sin_Euler}).}
\begin{equation}
\sin x=2^{n-1}\prod_{k=0}^{n-1}\sin\frac{(\pi k+x)}{n}
\label{sinx_as_product_sin_Euler}%
\end{equation}
for $x\rightarrow0$ gives $\prod_{k=1}^{n-1}\sin\frac{\pi k}{n}=\frac
{n}{2^{n-1}}$ and finally leads to%
\[
\Gamma\left(  nz+1\right)  =\left(  2\pi\right)  ^{-\frac{n-1}{2}}%
n^{nz+\frac{1}{2}}\prod_{k=1}^{n}\Gamma\left(  z+\frac{k}{n}\right)
\]
that equals (\ref{Gamma_multiplication_formula}), because $\prod_{k=1}%
^{n}\Gamma\left(  z+\frac{k}{n}\right)  =\Gamma\left(  z+1\right)  \prod
_{k=1}^{n-1}\Gamma\left(  z+\frac{k}{n}\right)  =z\prod_{k=0}^{n-1}%
\Gamma\left(  z+\frac{k}{n}\right)  $.

\medskip\refstepcounter{article}{\noindent\textbf{\thearticle. }%
}\ignorespaces\label{art_Gauss_integral_digamma_function} \emph{Gauss's
integral for the digamma function}. We will demonstrate Gauss's integral
\cite[{p. 160, formula [78]}]{Gauss1813}%
\begin{equation}
\psi\left(  z\right)  =\frac{d}{dz}\log\Gamma\left(  z\right)  =\int
_{0}^{\infty}\left(  \frac{e^{-t}}{t}-\frac{e^{-zt}}{1-e^{-t}}\right)  dt
\label{Digamma_Gauss_integral}%
\end{equation}
Instead of following Gauss's deduction, a more elegant derivation \cite[p.
247]{Whittaker_Watson} is obtained from Weierstrass's infinite product
(\ref{Weierstrass_product_Gamma_function}). The logarithm of Weierstrass's
infinite product (\ref{Weierstrass_product_Gamma_function}) is
\[
-\log\Gamma(z+1)=\gamma\,z+\sum_{n=1}^{\infty}\left(  \log\left(  1+\frac
{z}{n}\right)  -\frac{z}{n}\right)
\]
and differentiation yields%
\[
\psi\left(  z+1\right)  =\frac{d}{dz}\log\Gamma\left(  z+1\right)
=-\gamma-\sum_{n=1}^{\infty}\left(  \frac{1}{z+n}-\frac{1}{n}\right)
=-\gamma+z\sum_{n=1}^{\infty}\frac{1}{\left(  z+n\right)  n}%
\]
from which $\psi\left(  1\right)  =-\gamma$. The functional equation
(\ref{func_eq_digamma}) of $\psi\left(  z\right)  $ then shows that
\begin{equation}
\psi\left(  z\right)  =-\gamma-\frac{1}{z}-\lim_{k\rightarrow\infty}\sum
_{n=1}^{k}\left(  \frac{1}{z+n}-\frac{1}{n}\right)  \label{digamma_series}%
\end{equation}
The polygamma functions, defined as $\psi^{(n)}(z)=$ $\frac{d^{n+1}\ln
\Gamma(z)}{dz^{n+1}}$ with $\psi^{(0)}(z)=\psi(z)$, follow immediately from
(\ref{digamma_series}) as \cite[6.4.10]{Abramowitz}
\begin{equation}
\psi^{(n)}(z)=(-1)^{n+1}n!\sum_{k=0}^{\infty}\frac{1}{(z+k)^{n+1}}
\label{polygamma_series}%
\end{equation}
\qquad

Substituting $\int_{0}^{\infty}e^{-t\left(  z+n\right)  }dt=\frac{1}{z+n}$,
valid for $\operatorname{Re}\left(  t\right)  >0$, into (\ref{digamma_series}%
)\ yields%
\[
\psi\left(  z\right)  =-\gamma-\int_{0}^{\infty}e^{-tz}dt-\lim_{k\rightarrow
\infty}\int_{0}^{\infty}\left(  e^{-tz}-1\right)  \sum_{n=1}^{k}e^{-tn}dt
\]
With $\sum_{n=1}^{k}e^{-tn}=e^{-t}\frac{1-e^{-kt}}{1-e^{-t}}$, we have%
\begin{align*}
\psi\left(  z\right)   &  =-\gamma-\int_{0}^{\infty}e^{-tz}dt+\lim
_{k\rightarrow\infty}\int_{0}^{\infty}\frac{\left(  e^{-t}-e^{-t\left(
z+1\right)  }\right)  \left(  1-e^{-kt}\right)  }{1-e^{-t}}dt\\
&  =-\gamma-\int_{0}^{\infty}e^{-tz}dt+\int_{0}^{\infty}\frac{e^{-t}%
-e^{-t\left(  z+1\right)  }}{1-e^{-t}}dt-\lim_{k\rightarrow\infty}\int
_{0}^{\infty}\frac{\left(  1-e^{-tz}\right)  }{1-e^{-t}}e^{-\left(
k+1\right)  t}dt
\end{align*}
and%
\[
\psi\left(  z\right)  =-\gamma+\int_{0}^{\infty}\frac{e^{-t}-e^{-tz}}%
{1-e^{-t}}dt+\lim_{k\rightarrow\infty}\int_{0}^{\infty}\frac{\left(
1-e^{-tz}\right)  }{1-e^{-t}}e^{-\left(  k+1\right)  t}dt
\]
Since $\left\vert \int_{0}^{\infty}\frac{\left(  1-e^{-tz}\right)  }{1-e^{-t}%
}e^{-\left(  k+1\right)  t}dt\right\vert \leq\int_{0}^{\infty}\left\vert
\frac{\left(  1-e^{-tz}\right)  }{1-e^{-t}}\right\vert e^{-\left(  k+1\right)
t}dt\leq\max_{t\geq0}\left\vert \frac{\left(  1-e^{-tz}\right)  }{1-e^{-t}%
}\right\vert \frac{1}{k+1}\rightarrow0$ for large $k$, because the function
$\frac{\left(  1-e^{-tz}\right)  }{1-e^{-t}}$ is bounded for $t\geq0$, we
arrive at%
\[
\psi\left(  z\right)  =-\gamma+\int_{0}^{\infty}\frac{e^{-t}-e^{-tz}}%
{1-e^{-t}}dt
\]

Similarly as above, we rewrite Euler's constant $\gamma$ in
(\ref{def_Euler_constant_gamma}) by invoking again $\frac{1}{s}=\int
_{0}^{\infty}e^{-st}dt$ for $\operatorname{Re}\left(  t\right)  >0$ and by
using $\log x=\int_{1}^{x}\frac{1}{s}ds=\int_{1}^{x}\int_{0}^{\infty}%
e^{-st}dtds=\int_{0}^{\infty}\frac{e^{-t}-e^{-xt}}{t}dt$,%
\begin{align*}
\sum_{n=1}^{k}\frac{1}{n}-\log k  &  =\int_{0}^{\infty}e^{-t}\frac{1-e^{-kt}%
}{1-e^{-t}}dt-\int_{0}^{\infty}\frac{e^{-t}-e^{-kt}}{t}dt\\
&  =\int_{0}^{\infty}\left(  \frac{1}{1-e^{-t}}-\frac{1}{t}\right)
e^{-t}dt-\int_{0}^{\infty}\left(  \frac{e^{-t}}{1-e^{-t}}-\frac{1}{t}\right)
e^{-kt}dt
\end{align*}
The generating function of the Bernoulli numbers, written as
\begin{equation}
\frac{1\,}{e^{t}-1}=\frac{1}{t}-\frac{1}{2}+\sum_{n=1}^{\infty}B_{2n}%
\,\frac{t^{2n-1}}{(2n)!}\hspace{1cm}\text{for }\left\vert t\right\vert
\leq2\pi\label{gf_Bernoullinumbers}%
\end{equation}
shows that $\frac{e^{-t}}{1-e^{-t}}-\frac{1}{t}$ is continuous at $t=0$ and
bounded for $t\geq0$, leading for $k\rightarrow\infty$ to%
\begin{equation}
\gamma=\int_{0}^{\infty}\left(  \frac{1}{1-e^{-t}}-\frac{1}{t}\right)
e^{-t}dt \label{Euler_constant_as_integral}%
\end{equation}
and to Gauss's integral (\ref{Digamma_Gauss_integral}).

\medskip\refstepcounter{article}{\noindent\textbf{\thearticle. }%
}\ignorespaces\label{art_Stirling_formula} \emph{Stirling's asymptotic
formula}. Inspired by the Bernoulli generating function
(\ref{gf_Bernoullinumbers}), we rewrite Gauss's integral
(\ref{Digamma_Gauss_integral}) for $z\rightarrow z+1$,
\begin{align*}
\frac{d}{dz}\log\Gamma\left(  z+1\right)   &  =\int_{0}^{\infty}\left(
\frac{e^{-t}}{t}-\frac{e^{-zt}}{e^{t}-1}\right)  dt\\
&  =\int_{0}^{\infty}\left(  \frac{e^{-t}-e^{-zt}}{t}+\frac{e^{-zt}}%
{2}-e^{-zt}\left(  \frac{1}{e^{t}-1}-\frac{1}{t}+\frac{1}{2}\right)  \right)
dt\\
&  =\log z+\frac{1}{2z}-\int_{0}^{\infty}e^{-zt}\left(  \frac{1}{e^{t}%
-1}-\frac{1}{t}+\frac{1}{2}\right)  dt
\end{align*}

Since $0\leq\frac{1}{e^{t}-1}-\frac{1}{t}+\frac{1}{2}\leq\frac{1}{2}$ is
bounded for $t\geq0$ and continuous at $t=0$, the last integral is uniformly
convergent for $\operatorname{Re}\left(  z\right)  >0$ and, hence, can be
integrated from $1$ to $z$,%
\[
\log\Gamma\left(  z+1\right)  =z\log z-z+1+\frac{1}{2}\log z+\int_{0}^{\infty
}\frac{e^{-zt}-e^{-t}}{t}\left(  \frac{1}{e^{t}-1}-\frac{1}{t}+\frac{1}%
{2}\right)  dt
\]
With the functional equation (\ref{func_eq_Gamma}), we find%
\[
\log\Gamma\left(  z\right)  =\left(  z-\frac{1}{2}\right)  \log z-z+1+\int
_{0}^{\infty}\frac{e^{-zt}}{t}\left(  \frac{1}{e^{t}-1}-\frac{1}{t}+\frac
{1}{2}\right)  dt-\int_{0}^{\infty}\frac{e^{-t}}{t}\left(  \frac{1}{e^{t}%
-1}-\frac{1}{t}+\frac{1}{2}\right)  dt
\]
The last integral can be evaluated \cite[p. 249]{Whittaker_Watson} as
$1-\frac{1}{2}\log\left(  2\pi\right)  $ and we arrive at Binet's first
integral for $\operatorname{Re}\left(  z\right)  >0$,%
\begin{equation}
\log\Gamma\left(  z\right)  =\left(  z-\frac{1}{2}\right)  \log z-z+\frac
{1}{2}\log\left(  2\pi\right)  +\int_{0}^{\infty}\frac{e^{-zt}}{t}\left(
\frac{1}{e^{t}-1}-\frac{1}{t}+\frac{1}{2}\right)  dt
\label{Binet_logGamma_integral_1ste}%
\end{equation}
We observe from (\ref{Binet_logGamma_integral_1ste}) that
\[
\left\vert \int_{0}^{\infty}\frac{e^{-zt}}{t}\left(  \frac{1}{e^{t}-1}%
-\frac{1}{t}+\frac{1}{2}\right)  dt\right\vert \leq\max_{t\geq0}\left\vert
\frac{1}{t}\left(  \frac{1}{e^{t}-1}-\frac{1}{t}+\frac{1}{2}\right)
\right\vert \frac{1}{z}=\frac{B_{2}}{2z}=\frac{1}{12z}%
\]
so that, for any $\operatorname{Re}\left(  z\right)  >0$, a sharp upper bound
is found,%
\begin{equation}
\log\Gamma\left(  z\right)  \leq\left(  z-\frac{1}{2}\right)  \log
z-z+\frac{1}{2}\log\left(  2\pi\right)  +\frac{1}{12z}
\label{LogGamma(z)_bound_Re(z)>0}%
\end{equation}
In particular, the integral in (\ref{Binet_logGamma_integral_1ste}) is
positive for positive real $x>0$, leading to the bounds%
\[
\left(  x-\frac{1}{2}\right)  \log x-x+\frac{1}{2}\log\left(  2\pi\right)
\leq\log\Gamma\left(  x\right)  \leq\left(  x-\frac{1}{2}\right)  \log
x-x+\frac{1}{2}\log\left(  2\pi\right)  +\frac{1}{12x}%
\]

Invoking the generating function of the Bernoulli numbers
(\ref{gf_Bernoullinumbers}) in (\ref{Binet_logGamma_integral_1ste}), ignoring
the restriction $\left\vert t\right\vert <2\pi$, yields Stirling's
approximation \cite[6.1.40]{Abramowitz},%
\begin{align}
\log\Gamma\left(  z\right)   &  =\left(  z-\frac{1}{2}\right)  \log
z-z+\frac{1}{2}\log\left(  2\pi\right)  +\sum_{n=1}^{\infty}\,\frac{B_{2n}%
}{2n\left(  2n-1\right)  z^{2n-1}}\label{Stirling_logGamma_series}\\
&  =\left(  z-\frac{1}{2}\right)  \log z-z+\frac{1}{2}\log\left(  2\pi\right)
+\frac{1}{12z}-\frac{1}{360z^{3}}+O\left(  \frac{1}{z^{5}}\right) \nonumber
\end{align}
By substituting the partial fraction expansion $\frac{1}{t}\left(  \frac
{1}{e^{t}-1}-\frac{1}{t}+\frac{1}{2}\right)  =\sum_{n=1}^{\infty}\frac
{2}{t^{2}+4\pi^{2}n^{2}}$ in (\ref{Binet_logGamma_integral_1ste}),%
\[
\log\Gamma\left(  z\right)  =\left(  z-\frac{1}{2}\right)  \log z-z+\frac
{1}{2}\log\left(  2\pi\right)  +2\int_{0}^{\infty}e^{-zt}\sum_{n=1}^{\infty
}\frac{1}{t^{2}+4\pi^{2}n^{2}}dt
\]
and changing the integration variable,%
\begin{align*}
\int_{0}^{\infty}e^{-zt}\sum_{n=1}^{\infty}\frac{1}{t^{2}+4\pi^{2}n^{2}}dt  &
=\int_{0}^{\infty}\sum_{n=1}^{\infty}\frac{e^{-zt}}{\left(  1+\left(  \frac
{t}{2\pi n}\right)  ^{2}\right)  \left(  2\pi n\right)  ^{2}}dt=\frac{1}{2\pi
}\int_{0}^{\infty}\frac{du}{1+u^{2}}\sum_{n=1}^{\infty}\frac{e^{-2\pi nzu}}%
{n}\\
&  =-\frac{1}{2\pi}\int_{0}^{\infty}\frac{du}{1+u^{2}}\log\left(  1-e^{-2\pi
zu}\right)
\end{align*}
we find, after partial integration, Binet's second form \cite[p.
251]{Whittaker_Watson},\cite[p. 217]{Sansone},%
\begin{equation}
\log\Gamma\left(  z\right)  =\left(  z-\frac{1}{2}\right)  \log z-z+\frac
{1}{2}\log\left(  2\pi\right)  +2z\int_{0}^{\infty}\frac{\arctan t}{e^{2\pi
tz}-1}dt \label{Binet_logGamma_integral_2de}%
\end{equation}
Substituting in Binet's second form (\ref{Binet_logGamma_integral_2de}), the
Taylor series $\arctan t=\sum_{n=0}^{\infty}\frac{\left(  -1\right)  ^{n}%
}{2n+1}t^{2n+1}$ around $t=0$, valid for $\left\vert t\right\vert <1$,
reversing the integral and summation while ignoring the restriction
$\left\vert t\right\vert <1$, leads, with $\,B_{2n}=(-1)^{n-1}4n\int
_{0}^{\infty}\frac{t^{2n-1}}{e^{2\pi t}-1}dt$ for $n\geq1$, again to
Stirling's approximation (\ref{Stirling_logGamma_series}).

Titchmarsh \cite[p. 151]{Titchmarshfunctions} gives a third form%
\[
\log\Gamma\left(  z\right)  =\left(  z-\frac{1}{2}\right)  \log z-z+\frac
{1}{2}\log\left(  2\pi\right)  +\int_{0}^{\infty}\frac{\left[  t\right]
-t+\frac{1}{2}}{t+z}dt
\]
where $\left[  t\right]  $ is the nearest integer smaller than or equal to
$t$. Blagouchine \cite{Blagouchine_2016} gives seven series for $\log
\Gamma\left(  z\right)  $. Edwards \cite[p. 106-114]{Edwards} derives the
Stirling approximation (\ref{Stirling_logGamma_series}) for $\log\Pi\left(
z\right)  =\log\Gamma\left(  z+1\right)  $ by Euler-Maclaurin summation
(\ref{sum_EulerMaclaurin}).

Differentiating Binet's second form (\ref{Binet_logGamma_integral_2de})
results in a companion of Gauss's integral (\ref{Digamma_Gauss_integral}) for
the digamma function,
\begin{equation}
\psi\left(  z\right)  =\log z-\frac{1}{2z}-2\int_{0}^{\infty}\frac{t}{e^{2\pi
tz}-1}\frac{dt}{1+t^{2}}=\log z-\frac{1}{2z}-2\int_{0}^{\infty}\frac
{t}{e^{2\pi t}-1}\frac{dt}{z^{2}+t^{2}} \label{Digamma_Binets_integral}%
\end{equation}
from which we find that Euler's constant $\gamma=-\psi\left(  1\right)
=\frac{1}{2}+2\int_{0}^{\infty}\frac{t}{e^{2\pi t}-1}\frac{dt}{1+t^{2}}%
>\frac{1}{2}$. Since $\int_{0}^{\infty}\frac{t}{e^{2\pi t}-1}\frac{dt}%
{1+t^{2}}<\int_{0}^{\infty}\frac{tdt}{e^{2\pi t}-1}$ and $\Gamma\left(
s\right)  \zeta\left(  s\right)  =\int_{0}^{\infty}\frac{t^{s-1}dt}{e^{t}-1}$,
an upper bound follows as $\gamma<\frac{1}{2}+\frac{2}{\left(  2\pi\right)
^{2}}\zeta\left(  2\right)  =0.5833$.

\medskip\refstepcounter{article}{\noindent\textbf{\thearticle. }%
}\ignorespaces\label{art_asymptotic_behavior_LogGamma} \emph{Asymptotic
behavior of }$\log\Gamma\left(  z\right)  $. For large $r$ and $\theta\neq\pi
$, Stirling's formula (\ref{Stirling_logGamma_series}) shows that
\[
\log\Gamma\left(  b+re^{i\theta}\right)  =\left(  b+re^{i\theta}-\frac{1}%
{2}\right)  \log\left(  b+re^{i\theta}\right)  -b-re^{i\theta}+\frac{1}{2}%
\log\left(  2\pi\right)  +O\left(  \frac{1}{r}\right)
\]
With
\begin{align*}
\log\left(  b+re^{i\theta}\right)   &  =\log\left(  re^{i\theta}\left(
1+r^{-1}e^{-i\theta}b\right)  \right)  =\log re^{i\theta}+\log\left(
1+r^{-1}e^{-i\theta}b\right) \\
&  =\log r+i\theta+r^{-1}e^{-i\theta}b+O\left(  \frac{1}{r^{2}}\right)
\end{align*}
we have
\begin{align*}
\log\Gamma\left(  b+re^{i\theta}\right)   &  =\left(  b+re^{i\theta}-\frac
{1}{2}\right)  \log r+i\theta\left(  b+re^{i\theta}-\frac{1}{2}\right)
+\left(  br^{-1}+e^{i\theta}-\frac{1}{2}r^{-1}\right)  e^{-i\theta}b\\
&  \hspace{1cm}-b-re^{i\theta}+\frac{1}{2}\log\left(  2\pi\right)  +O\left(
\frac{1}{r}\right) \\
&  =\left(  b-\frac{1}{2}+r\cos\theta\right)  \log r-\theta\left(  r\sin
\theta\right)  -r\cos\theta+\frac{1}{2}\log\left(  2\pi\right) \\
&  \hspace{1cm}+i\left\{  r\ln r\sin\theta+\left(  b-\frac{1}{2}+r\cos
\theta\right)  \theta-r\sin\theta\right\}  +O\left(  \frac{1}{r}\right)
\end{align*}
Thus,%
\[
\Gamma\left(  b+re^{i\theta}\right)  =\sqrt{2\pi}e^{r\left(  \log\left(
r\right)  -1\right)  \cos\theta-\theta r\sin\theta+\left(  b-\frac{1}%
{2}\right)  \log\left(  r\right)  }e^{i\left(  r\left(  \log\left(  r\right)
-1\right)  \sin\theta+r\theta\cos\theta+\left(  b-\frac{1}{2}\right)
\theta\right)  }\left(  1+O\left(  \frac{1}{r}\right)  \right)
\]
from which%
\[
\left\vert \Gamma\left(  b+re^{i\theta}\right)  \right\vert =r^{b-\frac{1}{2}%
}e^{r\left(  \log\left(  r\right)  -1\right)  \cos\theta-\theta r\sin\theta
}\left(  1+O\left(  \frac{1}{r}\right)  \right)
\]
Hence, it holds that%
\begin{equation}
\frac{1}{\left\vert \Gamma\left(  b+re^{i\theta}\right)  \right\vert }%
=\frac{r^{\frac{1}{2}-b}}{\sqrt{2\pi}}e^{-r\left(  \ln\left(  r\right)
-1\right)  \cos\theta+\theta r\sin\theta}\left(  1+O\left(  \frac{1}%
{r}\right)  \right)  \label{asymptotic_Gamma(b+rexp(it))}%
\end{equation}
illustrating that $\lim_{r\rightarrow\infty}\frac{1}{\left\vert \Gamma\left(
b+re^{i\theta}\right)  \right\vert }=0$ for $-\frac{\pi}{2}<\theta<\frac{\pi
}{2}$. However, if $\theta=\frac{\pi}{2}$, then in contrast
\[
\frac{1}{\left\vert \Gamma\left(  b+ir\right)  \right\vert }=\frac{r^{\frac
{1}{2}-b}}{\sqrt{2\pi}}e^{\frac{\pi}{2}r}\left(  1+O\left(  \frac{1}%
{r}\right)  \right)  .
\]

\subsection{Complex integrals for the Gamma function}

\label{sec_integrals_Gamma_function}

\medskip\refstepcounter{article}{\noindent\textbf{\thearticle. }%
}\ignorespaces\label{art_Hankel_integral} \emph{Hankel's integral}. We derive
Hankel's integral%
\begin{equation}
\frac{1}{\Gamma\left(  z\right)  }=\frac{1}{2\pi i}\int_{C}w^{-z}e^{w}dw
\label{Hankel_contour_integral_1opGamma}%
\end{equation}
where the contour $C$ starts at $-\infty$ below the real axis, encircles the
origin at $z=0$ and returns above the negative real axis again to $-\infty$.
Such a contour around the branch-cut (here the negative real axis) is
\textquotedblleft classical\textquotedblright\ in integrals containing
$w^{\alpha}$ such as Mellin transforms \cite{Titchmarshfourier}. If
$\varepsilon$ is the radius of a circle at the origin, then evaluation of the
integral along the contour $C$ yields%
\begin{align*}
\int_{C}w^{-z}e^{w}dw  &  =-\int_{\infty}^{\varepsilon}\left(  xe^{-i\pi
}\right)  ^{-z}e^{-x}dx+i\int_{-\pi}^{\pi}\left(  \varepsilon e^{i\theta
}\right)  ^{-z}e^{\varepsilon e^{i\theta}}\varepsilon e^{i\theta}d\theta
-\int_{\varepsilon}^{\infty}\left(  xe^{i\pi}\right)  ^{-z}e^{-x}dx\\
&  =\left(  e^{i\pi z}-e^{-i\pi z}\right)  \int_{\varepsilon}^{\infty}%
x^{-z}e^{-x}dx+i\varepsilon^{1-z}\int_{-\pi}^{\pi}e^{i\theta\left(
1-z\right)  }e^{\varepsilon e^{i\theta}}d\theta
\end{align*}
If $\operatorname{Re}\left(  1-z\right)  >0$, then we take the limit
$\varepsilon\rightarrow0$ and find with Euler's integral
(\ref{Euler_integral_Gamma_functie})%
\[
\int_{C}w^{-z}e^{w}dw=2i\sin\pi z\int_{0}^{\infty}x^{1-z-1}e^{-x}dx=2i\sin\pi
z\Gamma\left(  1-z\right)
\]
After replacing $z$ by $1-z$, we obtain a contour integral for the Gamma
function,
\begin{equation}
\Gamma\left(  z\right)  =\frac{\pi}{\sin\pi z}\frac{1}{2\pi i}\int_{C}%
w^{z-1}e^{w}dw \label{Contour_Gamma_function}%
\end{equation}
The reflection formula (\ref{Gamma_reflection_formula}), $\Gamma\left(
z\right)  \Gamma\left(  1-z\right)  =\frac{\pi}{\sin\pi z}$, of the Gamma
function leads to Hankel's contour integral
(\ref{Hankel_contour_integral_1opGamma}). Although the derivation was
restricted to $\operatorname{Re}\left(  z\right)  >0$ in
(\ref{Contour_Gamma_function}), by analytic continuation (see e.g.
\cite[Chapter IV]{Titchmarshfunctions},\cite[Chapter III]{Evgrafov_1965}),
(\ref{Hankel_contour_integral_1opGamma}) as well as
(\ref{Contour_Gamma_function}) hold for all $z\in\mathbb{C}$. The contour
integral (\ref{Contour_Gamma_function}) for the Gamma function demonstrates
that $\Gamma\left(  z\right)  $ has simple poles at $z=-n$ for each integer
$n$ (due to $\sin\pi z$), in agreement with Gauss's product
(\ref{Gauss_product_Gamma_function}).

We can deform\footnote{The same deformation holds for the contour integral
(\ref{Contour_Gamma_function}) as well.} the contour $C$ into $C_{\phi}$ by
tilting the straight line above the negative real axis over an angle $\phi$
and the straight line below the negative real axis over an angle $-\phi$.
Indeed, consider the contour $L$ consisting of the contour $C$, the circle
segment at infinity from the angle $\pi$ to the angle $\phi$, followed by the
line $w=re^{i\phi}$, where $r$ ranges from infinity towards $\rho$, the circle
centered at $w=0$ with radius $\rho$ turning from angle $\phi$ towards $-\phi$
and complemented by the line $w=re^{-i\phi}$ and infinite circle segment
towards the begin of the contour $C$. The integral $\int_{L}w^{-z}e^{w}dw=0$,
because the contour $L$ encloses an analytic region of the function
$w^{-z}e^{w}$. If $\phi>\frac{\pi}{2}$, then the part of $L$ along the circle
segment at infinity, $\lim_{r\rightarrow\infty}i\int_{\pi}^{\phi}%
r^{1-z}e^{\left(  1-z\right)  i\theta}e^{re^{i\theta}}d\theta=0$. Hence,
combining all contributions results in
\begin{equation}
\frac{1}{\Gamma\left(  z\right)  }=\frac{1}{2\pi i}\int_{C_{\phi}}w^{-z}%
e^{w}dw \label{Hankel_contour_phi_integral_1opGamma}%
\end{equation}
where the contour $C_{\phi}$ starts at infinity on the straight line at the
angle $-\frac{\pi}{2}>\phi>-\pi$ below the real axis until the circle at the
origin with radius $\rho$ that turns over the angle $-\phi$ to $\phi$ until
hitting the line $w=re^{i\phi}$ along which it passes towards infinity again.
The contour in (\ref{Hankel_contour_integral_1opGamma}) is the particular case
where $C=C_{\pi}$ and $\rho=\varepsilon$. Hankel's integral in
(\ref{Hankel_contour_phi_integral_1opGamma}) can be written as%
\begin{equation}
\frac{1}{\Gamma\left(  z\right)  }=\frac{1}{2\pi i}\int_{c-\infty e^{-i\phi}%
}^{c+\infty e^{i\phi}}w^{-z}e^{w}dw\hspace{1cm}\text{ with }c>0\text{ and
}\frac{\pi}{2}<\phi<\pi\text{ } \label{Hankel_contour_long_line}%
\end{equation}
Since multiplying $c^{\prime}=ac>0$ provided $a>0$, the map $w\rightarrow aw$
for any real, positive number $a$ yields%
\begin{equation}
\frac{a^{z-1}}{\Gamma\left(  z\right)  }=\frac{1}{2\pi i}\int_{c^{\prime
}-\infty e^{-i\phi}}^{c^{\prime}+\infty e^{i\phi}}w^{-z}e^{aw}dw\hspace
{1cm}\text{with }c^{\prime}>0\text{ and }\frac{\pi}{2}<\phi<\pi
\label{Hankel_e^(aw)_a^(1-z)_op_Gamma(z)}%
\end{equation}

\section{Complex function theorem}

\label{sec_complex_function_theorema}We evaluate integrals of a general kind.

\medskip\refstepcounter{article}{\noindent\textbf{\thearticle. }%
}\ignorespaces\label{art_Cauchy_type_integral} \emph{A Cauchy-type integral}.

\begin{theorem}
\label{theorem_entire_function_Cauchy} Let $f\left(  z\right)  $ be an entire
function, that is real on the real axis and $\lim_{r\rightarrow\infty}$
$\frac{f\left(  re^{i\theta}\right)  }{r^{2}}=0$ for $\theta=0$ and
$\theta=\frac{\pi}{2}$. If $\sigma$ is a positive real number, then it holds
that
\begin{equation}
\int_{-\infty}^{\infty}\frac{f(\sigma+it)}{|\sigma+it|^{2}}dt=\frac{\pi
}{\sigma}\;f(2\,\sigma) \label{Entire_function_Cauchy_integral}%
\end{equation}

\end{theorem}

\textbf{Proof}: Consider the integral $I=\frac{1}{2\pi i}\int_{L}\frac
{f(z)}{z(2\sigma-z)}dz$ where the contour $L$ is taken counter-clockwise round
the rectangle formed by the lines $\operatorname{Im}(z)=-T$,
$\operatorname{Re}(z)=a>2\sigma$, $\operatorname{Im}(z)=T$ and
$\operatorname{Re}(z)=\sigma>0$. The contour $L$ encloses the pole on the real
axis at $x=2\sigma$. Since the entire function $f(z)$ is analytic inside and
on $L$, the integral equals $I=-\frac{f(2\sigma)}{2\sigma}$. On the other
hand, evaluating the integral $I$ along the contour $L$ yields,
\begin{align*}
I  &  =\frac{1}{2\pi i}\int_{\sigma}^{a}\frac{f(x-iT)}{(x-iT)(2\sigma
-x+iT)}\;dx+\frac{1}{2\pi}\;\int_{-T}^{T}\frac{f(a+it)}{(a+it)(2\sigma
-a-it)}dt\\
&  -\frac{1}{2\pi i}\int_{\sigma}^{a}\frac{f(x+iT)}{(x+iT)(2\sigma
-x-iT)}\;dx-\frac{1}{2\pi}\;\int_{-T}^{T}\frac{f(\sigma+it)}{|\sigma+it|^{2}%
}\;dt
\end{align*}
Combined and rewritten leads to%
\begin{align*}
f(2\sigma)  &  =\frac{\sigma}{\pi}\;\int_{-T}^{T}\frac{f(\sigma+it)}%
{\sigma^{2}+t^{2}}dt-\frac{\sigma}{\pi}\;\int_{-T}^{T}\frac{f(a+it)}%
{(a+it)(2\sigma-a-it)}dt\\
&  +\frac{\sigma}{\pi i}\int_{\sigma}^{a}\operatorname{Im}\left[
\frac{f(x+iT)}{(x+iT)(2\sigma-x-iT)}\right]  dx
\end{align*}
Since $\lim_{z\rightarrow\infty}\frac{f\left(  z\right)  }{z^{2}}=0$ for
$z=re^{i\theta}$ with $\theta=0$ and $\theta=\frac{\pi}{2}$, the second and
third integral vanish, demonstrating the Theorem.\hfill$\square\medskip$

\medskip\refstepcounter{article}{\noindent\textbf{\thearticle. }%
}\ignorespaces\label{art_Mellin_transform_product_Gamma} \emph{Mellin
transform of a product of Gamma functions}. Let $p_{1},p_{2},\ldots,p_{n}$ be
different real, positive numbers, then the Mellin transform%

\begin{equation}%
{\displaystyle\prod\limits_{j=1}^{n}}
\Gamma\left(  s+p_{j}\right)  =\int_{0}^{\infty}u^{s-1}g\left(  u;\left\{
p_{j}\right\}  _{1\leq j\leq n}\right)  du \label{Melling_prod_Gamma}%
\end{equation}
has an inverse
\begin{equation}
g\left(  u;\left\{  p_{j}\right\}  _{1\leq j\leq n}\right)  =\frac{1}{2\pi
i}\int_{c-i\infty}^{c+i\infty}%
{\displaystyle\prod\limits_{j=1}^{n}}
\Gamma\left(  s+p_{j}\right)  u^{-s}ds\hspace{1cm}\text{with \thinspace}c>0
\label{Mellin_inverse_prod_Gamma}%
\end{equation}
which can be evaluated. The contour in (\ref{Mellin_inverse_prod_Gamma}) can
be closed over the negative $\operatorname{Re}\left(  s\right)  $-plane and
encloses the simple poles at $s=-q_{j}-p_{j}$ for integers $q_{j}\geq0$ for
$1\leq j\leq n$ on the negative real axis. Cauchy's residue theorem
\cite{Titchmarshfunctions} leads to%

\[
g\left(  u;\left\{  p_{j}\right\}  _{1\leq j\leq n}\right)  =\sum_{j=1}%
^{n}\sum_{q_{j}=0}^{\infty}\lim_{s\rightarrow-q_{j}-p_{j}}\left(
s+q_{j}+p_{j}\right)  \Gamma\left(  s+p_{j}\right)
{\displaystyle\prod\limits_{m=1;m\neq j}^{n}}
\Gamma\left(  s+p_{m}\right)  u^{-s}%
\]
Iterating the functional equation (\ref{func_eq_Gamma}) of the Gamma function
$\left(  q_{j}+1\right)  $-times gives
\begin{align*}
\lim_{s\rightarrow-q_{j}-p_{j}}\left(  s+q_{j}+p_{j}\right)  \Gamma\left(
s+p_{j}\right)   &  =\lim_{s\rightarrow-q_{j}-p_{j}}\frac{\Gamma\left(
s+p_{j}+q_{j}+1\right)  \left(  s+p_{j}+q_{j}\right)  }{\left(  s+p_{j}%
\right)  \left(  s+p_{j}+1\right)  \ldots\left(  s+p_{j}+q_{j}-1\right)
\left(  s+p_{j}+q_{j}\right)  }\\
&  =\frac{\Gamma\left(  1\right)  }{\left(  -q_{j}\right)  \left(
-q_{j}+1\right)  \ldots\left(  -1\right)  }=\frac{\left(  -1\right)  ^{q_{j}}%
}{\left(  q_{j}\right)  !}%
\end{align*}
Hence, we obtain%
\[
g\left(  u;\left\{  p_{j}\right\}  _{1\leq j\leq n}\right)  =\sum_{j=1}%
^{n}\sum_{q_{j}=0}^{\infty}\frac{\left(  -1\right)  ^{q_{j}}u^{q_{j}+p_{j}}%
}{\left(  q_{j}\right)  !}%
{\displaystyle\prod\limits_{m=1;m\neq j}^{n}}
\Gamma\left(  -q_{j}+p_{m}-p_{j}\right)
\]
With the reflection formula (\ref{Gamma_reflection_formula})%
\[
\Gamma\left(  -q_{j}+p_{m}-p_{j}\right)  =\frac{\pi\left(  -1\right)  ^{q_{j}%
}}{\sin\pi\left(  p_{m}-p_{j}\right)  \Gamma\left(  1+q_{j}+p_{j}%
-p_{m}\right)  }%
\]
we have%
\begin{align*}%
{\displaystyle\prod\limits_{m=1;m\neq j}^{n}}
\Gamma\left(  -q_{j}+p_{m}-p_{j}\right)   &  =%
{\displaystyle\prod\limits_{m=1;m\neq j}^{n}}
\frac{\pi\left(  -1\right)  ^{q_{j}}}{\sin\pi\left(  p_{m}-p_{j}\right)
\Gamma\left(  1+q_{j}+p_{j}-p_{m}\right)  }\\
&  =\frac{\pi^{n-1}\left(  -1\right)  ^{\left(  n-1\right)  q_{j}}}{%
{\displaystyle\prod\limits_{m=1}^{j-1}}
\sin\pi\left(  p_{m}-p_{j}\right)
{\displaystyle\prod\limits_{m=j+1}^{n}}
\sin\pi\left(  p_{m}-p_{j}\right)  }\\
&  \times\frac{1}{%
{\displaystyle\prod\limits_{m=1}^{j-1}}
\Gamma\left(  1+q_{j}+p_{j}-p_{m}\right)
{\displaystyle\prod\limits_{m=j+1}^{n}}
\Gamma\left(  1+q_{j}+p_{j}-p_{m}\right)  }%
\end{align*}
we arrive at the series%
\begin{equation}
g\left(  u;\left\{  p_{j}\right\}  _{1\leq j\leq n}\right)  =\sum_{j=1}%
^{n}\frac{\left(  -1\right)  ^{j-1}\pi^{n-1}u^{p_{j}}}{%
{\displaystyle\prod\limits_{m=1}^{j-1}}
\sin\pi\left(  p_{j}-p_{m}\right)
{\displaystyle\prod\limits_{m=j+1}^{n}}
\sin\pi\left(  p_{m}-p_{j}\right)  }\sum_{q_{j}=0}^{\infty}\frac{\left(
-1\right)  ^{nq_{j}}u^{q_{j}}}{\left(  q_{j}\right)  !%
{\displaystyle\prod\limits_{m=1;m\neq j}^{n}}
\Gamma\left(  1+q_{j}+p_{j}-p_{m}\right)  }
\label{Mellin_inverse_prod_Gamma_Taylor_series}%
\end{equation}
with the convention that $%
{\displaystyle\prod\limits_{m=a}^{b}}
f_{m}=1$ if $b<a$.

\textbf{Examples} If $n=1$, then we retrieve the classical Mellin transform of
the pair $e^{-u}$ and $\Gamma\left(  s\right)  $,%
\[
g\left(  u;p_{1}\right)  =u^{p_{1}}\sum_{q_{1}=0}^{\infty}\frac{\left(
-1\right)  ^{q_{1}}u^{q_{1}}}{\left(  q_{1}\right)  !}=u^{p_{1}}e^{-u}%
\]

If $p_{j}=\frac{j}{n}$ for $1\leq j\leq n-1$ as in Gauss's multiplication
formula (\ref{Gamma_multiplication_formula}) without $j=0$ factor, then
(\ref{Mellin_inverse_prod_Gamma_Taylor_series}), denoted as $g\left(
u;\left\{  \frac{j}{n}\right\}  _{1\leq j\leq n-1}\right)  =h_{n}\left(
u\right)  $ reduces to
\begin{equation}
h_{n}\left(  u\right)  =\pi^{n-2}\sum_{j=1}^{n-1}\frac{\left(  -1\right)
^{j-1}u^{\frac{j}{n}}}{%
{\displaystyle\prod\limits_{k=1}^{j-1}}
\sin\left(  \frac{\pi k}{n}\right)
{\displaystyle\prod\limits_{k=1}^{n-j-1}}
\sin\left(  \frac{\pi k}{n}\right)  }\sum_{q_{j}=0}^{\infty}\frac{\left(
\left(  -1\right)  ^{n-1}\right)  ^{q_{j}}u^{q_{j}}}{\left(  q_{j}\right)  !%
{\displaystyle\prod\limits_{k=1}^{j-1}}
\Gamma\left(  1+q_{j}+\frac{k}{n}\right)
{\displaystyle\prod\limits_{m=1}^{n-j-1}}
\Gamma\left(  1+q_{j}-\frac{m}{n}\right)  }
\label{Mellin_inverse_prod_Gamma_multiplication_Taylor_series}%
\end{equation}
In particular for $n=3$, then%
\[
h_{3}\left(  u\right)  =\pi\frac{u^{\frac{1}{3}}}{\sin\left(  \frac{\pi}%
{3}\right)  }\sum_{q_{1}=0}^{\infty}\frac{u^{q_{1}}}{\left(  q_{j}\right)
!\Gamma\left(  1+q_{1}-\frac{1}{3}\right)  }-\pi\frac{u^{\frac{2}{3}}}%
{\sin\left(  \frac{\pi}{3}\right)  }\sum_{q_{2}=0}^{\infty}\frac{u^{q_{2}}%
}{\left(  q_{j}\right)  !\Gamma\left(  1+q_{2}+\frac{1}{3}\right)  }%
\]

We briefly summarize the theory of the modified Bessel function $K_{\nu}(z)$,
but strongly advise to consult the monumental treatise of Watson
\cite{Watson}. The modified Bessel functions $I_{\nu}\left(  z\right)  $ and
$K_{\nu}(z)$ are defined \cite[p. 77-78]{Watson} as the two independent
solutions of the modified Bessel differential equation%
\[
z^{2}\frac{d^{2}y}{dz^{2}}+z\frac{dy}{dz}-\left(  z^{2}+\nu^{2}\right)  y=0
\]
Both $I_{\nu}\left(  z\right)  $ and $K_{\nu}(z)$ are entire functions in $v$
for $z\neq0$ and analytic in $z$, except for a cut along the negative real
axis. The function $I_{\nu}\left(  z\right)  $ in $z$ possesses the Taylor
series%
\begin{equation}
I_{\nu}\left(  z\right)  =\sum_{k=0}^{\infty}\frac{\left(  \frac{1}%
{2}z\right)  ^{\nu+2k}}{k!\Gamma(\nu+k+1)} \label{def_Taylor_series_BesselI}%
\end{equation}
The modified Bessel functions $K_{\nu}(z)$ is defined as%
\begin{equation}
K_{\nu}(z)=\frac{\pi}{2\sin\pi\nu}\left(  I_{-\nu}\left(  z\right)  -I_{\nu
}\left(  z\right)  \right)  \label{def_BesselK}%
\end{equation}
clearly even in the order $\nu$, $K_{\nu}(z)=K_{-\nu}(z)$ for all $z\neq0$,
with Taylor series%
\begin{equation}
K_{\nu}(z)=\frac{\pi}{2\sin\pi\nu}\left(  \left(  \frac{1}{2}z\right)  ^{-\nu
}\sum_{k=0}^{\infty}\frac{\left(  \frac{1}{4}z^{2}\right)  ^{k}}%
{k!\Gamma(1-\nu+k)}-\left(  \frac{1}{2}z\right)  ^{\nu}\sum_{k=0}^{\infty
}\frac{\left(  \frac{1}{4}z^{2}\right)  ^{k}}{k!\Gamma(1+\nu+k)}\right)
\label{Taylor_series_BesselK}%
\end{equation}
The Taylor series (\ref{Taylor_series_BesselK}) shows that%
\begin{equation}
h_{3}\left(  u\right)  =2u^{\frac{1}{2}}K_{\frac{1}{3}}(2\sqrt{u})
\label{h_3_BesselK}%
\end{equation}
On the other hand, the integral \cite[eq. (28.73), p. 55]{Rademacher}, valid
for $\operatorname{Re}\left(  s\right)  >0$ and $\operatorname{Re}\left(
p\right)  >0$,%
\[
\Gamma\left(  s\right)  \Gamma\left(  s+p\right)  =2\int_{0}^{\infty}%
K_{p}\left(  2\sqrt{x}\right)  x^{\frac{p}{2}+s-1}dx
\]
is a special case of (\ref{Melling_prod_Gamma}) and corresponds to
$h_{3}\left(  u\right)  $ in (\ref{h_3_BesselK}) with $p=\frac{1}{3}$ and
$s\rightarrow s+\frac{1}{3}$.

\section{Inverse Laplace transform}

\label{sec_inverse_Laplace_transform}The Laplace transform for complex $z$ is
defined (see e.g. \cite{Titchmarshfourier}, \cite[Chapter VII]{Evgrafov_1965},
\cite{Widder})\ as
\begin{equation}
\varphi(z)=\int_{0}^{\infty}e^{-zt}f(t)dt \label{def_pgf_continuousrv}%
\end{equation}
with the inverse transform,
\begin{equation}
f(t)=\frac{1}{2\pi i}\int_{c-i\infty}^{c+i\infty}\varphi(z)e^{zt}dz
\label{inverse_Laplace}%
\end{equation}
where $c$ is the smallest real value of $\operatorname{Re}(z)$ for which the
integral in (\ref{def_pgf_continuousrv}) converges.

We evaluate the integral in (\ref{inverse_Laplace}) along the line $z=c+iw$,
\[
f(t)=\frac{e^{tc}}{2\pi}\int_{-\infty}^{\infty}\varphi(c+iw)e^{itw}dw
\]
After writing the integrand in a real and an imaginary part,%
\begin{align*}
f(t)  &  =\frac{e^{tc}}{2\pi}\int_{-\infty}^{\infty}\left\{  \operatorname{Re}%
\varphi(c+iw)\cos tw-\operatorname{Im}\left(  \varphi(c+iw)\right)  \sin
tw\right\}  dw\\
&  +i\frac{e^{tc}}{2\pi}\int_{-\infty}^{\infty}\left\{  \operatorname{Re}%
\left(  \varphi(c+iw)\right)  \sin tw+\operatorname{Im}\left(  \varphi
(c+iw)\right)  \cos tw\right\}  dw
\end{align*}
we find, after separating the real and imaginary part, that%
\[
\left\{
\begin{array}
[c]{c}%
f(t)=\frac{e^{tc}}{2\pi}\int_{-\infty}^{\infty}\left\{  \operatorname{Re}%
\varphi(c+iw)\cos tw-\operatorname{Im}\varphi(c+iw)\sin tw\right\}  dw\\
0=\int_{-\infty}^{\infty}\left\{  \operatorname{Re}\varphi(c+iw)\sin
tw+\operatorname{Im}\varphi(c+iw)\cos tw\right\}  dw
\end{array}
\right.
\]
On the other hand, it follows from (\ref{def_pgf_continuousrv}) that%
\[
\left\{
\begin{array}
[c]{c}%
\operatorname{Re}\varphi(c+iw)=\int_{0}^{\infty}e^{-ct}f(t)\cos wtdt\\
\operatorname{Im}\varphi(c+iw)=-\int_{0}^{\infty}e^{-ct}f(t)\sin wtdt
\end{array}
\right.
\]
and that $\operatorname{Re}\varphi(c+iw)\cos tw$ and $\operatorname{Im}%
\varphi(c+iw)\sin tw$ are even in $w$. Likewise, $\operatorname{Re}%
\varphi(c+iw)\sin tw$ and $\operatorname{Im}\varphi(c+iw)\cos tw$ are odd in
$w$. Hence, we arrive at%
\begin{equation}
\left\{
\begin{array}
[c]{c}%
f(t)=\frac{e^{tc}}{\pi}\int_{0}^{\infty}\left\{  \operatorname{Re}%
\varphi(c+iw)\cos tw-\operatorname{Im}\varphi(c+iw)\sin tw\right\}  dw\\
\int_{-\infty}^{\infty}\operatorname{Re}\varphi(c+iw)\sin twdw=\int_{-\infty
}^{\infty}\operatorname{Im}\varphi(c+iw)\cos twdw=0
\end{array}
\right.  \label{inverse_Laplace_Berberan-Santos_1}%
\end{equation}
The derivation is a corrected version of \cite{Berberan-SantosI_2005}.

Berberan-Santos\footnote{Private communication.} suggests to proceed a step
further by defining $f\left(  t\right)  =0$ for $t<0$. In that case, it
follows from (\ref{inverse_Laplace_Berberan-Santos_1}) that%
\[
0=\int_{0}^{\infty}\left\{  \operatorname{Re}\varphi(c+iw)\cos
tw-\operatorname{Im}\varphi(c+iw)\sin tw\right\}  dw\hspace{1cm}\text{for
}t<0
\]
which leads, after replacing $t\rightarrow-t$, to%
\[
\int_{0}^{\infty}\operatorname{Re}\varphi(c+iw)\cos twdw=-\int_{0}^{\infty
}\operatorname{Im}\varphi(c+iw)\sin twdw\hspace{1cm}\text{for }t>0
\]
The final resulting set of integral equations, subject to \textquotedblleft%
$f\left(  t\right)  =0$ for $t<0$\textquotedblright, simplifies to%
\begin{equation}
\left\{
\begin{array}
[c]{cc}%
f(t)=\frac{2e^{tc}}{\pi}\int_{0}^{\infty}\operatorname{Re}\varphi(c+iw)\cos
twdw & \text{for }t>0\\
f(t)=-\frac{2e^{tc}}{\pi}\int_{0}^{\infty}\operatorname{Im}\varphi(c+iw)\sin
twdw & \text{for }t>0\\
\int_{-\infty}^{\infty}\operatorname{Re}\varphi(c+iw)\sin twdw=\int_{-\infty
}^{\infty}\operatorname{Im}\varphi(c+iw)\cos twdw=0 &
\end{array}
\right.  \label{inverse_Laplace_transform_Berberan_Santos}%
\end{equation}
The general form (\ref{inverse_Laplace}) does not impose the restriction
\textquotedblleft$f\left(  t\right)  =0$ for $t<0$\textquotedblright\ and is
continuous\footnote{Indeed,%
\begin{align*}
\left\vert f(\varepsilon)-f(-\varepsilon)\right\vert  &  \leq\frac{1}{2\pi
}\int_{-\infty}^{\infty}\left\vert \varphi(c+iw)\right\vert \left\vert
e^{\varepsilon\left(  c+iw\right)  }-e^{-\varepsilon\left(  c+iw\right)
}\right\vert dw\\
&  =\frac{1}{\pi}\lim_{T\rightarrow\infty}\int_{-T}^{T}\left\vert
\varphi(c+iw)\right\vert \left\vert \sinh\left(  \varepsilon\left(
c+iw\right)  \right)  \right\vert dw
\end{align*}
Since $\left\vert \sinh\left(  \varepsilon\left(  c+iw\right)  \right)
\right\vert =\left\vert \sinh\varepsilon c\cos\varepsilon w+\sin\varepsilon
w\cosh\varepsilon c\right\vert <\left\vert \sinh\varepsilon c\right\vert
+\left\vert \sin\varepsilon w\right\vert +O\left(  \varepsilon^{2}\right)  $
and choosing $\varepsilon=T^{-2-\delta}$,%
\begin{align*}
\int_{-T}^{T}\left\vert \varphi(c+iw)\right\vert \left\vert \sinh\left(
\varepsilon\left(  c+iw\right)  \right)  \right\vert dw  &  \leq
cT^{-2-\delta}\int_{-T}^{T}\left\vert \varphi(c+iw)\right\vert dw+T^{-2-\delta
}\int_{-T}^{T}\left\vert w\varphi(c+iw)\right\vert dw\\
&  \leq2cT^{-1-\delta}\max\left\vert \varphi(c+iw)\right\vert +T^{-\delta}%
\max\left\vert \varphi(c+iw)\right\vert \leq AT^{-\delta}%
\end{align*}
which can be made arbitrarily small for large $T$ so that $\left\vert
f(\varepsilon)-f(-\varepsilon)\right\vert \rightarrow0$ when $\varepsilon
\rightarrow0$.} at $t=0$. The restriction \textquotedblleft$f\left(  t\right)
=0$ for $t<0$\textquotedblright\ is only continuous at $t=0$ if $f\left(
0\right)  =0$. In spite of this concern, Berberan-Santos
\cite{Berberan-SantosI_2005} evaluates various Laplace inverses via
$f(t)=\frac{2e^{tc}}{\pi}\int_{0}^{\infty}\operatorname{Re}\varphi(c+iw)\cos
twdw$.

In Fourier transforms, the inverse has a similar form as the transform itself
and tables of Fourier transforms can thus be used in two directions. Gross
\cite{Gross_1950} has written\footnote{Professor Apelblat has informed me
about this note.} a note on the question when the inverse Laplace transform
(\ref{inverse_Laplace}) is of the same form as the Laplace transform
(\ref{def_pgf_continuousrv}) itself. In particular, if $\varphi(z)=\int
_{0}^{\infty}e^{-zt}f(t)dt=\mathcal{L}\left[  f(t)\right]  $, then Gross
\cite{Gross_1950} asks when the inverse is of the form
\[
f\left(  z\right)  =\int_{0}^{\infty}e^{-zt}g(t)dt=\mathcal{L}\left[
g(t)\right]
\]
for real $f$ and $g$. Formal subsitution of the latter into the former and
reversing the integrals yields a Stieltjes transform \cite[Chapter
VIII]{Widder},\cite[9.15, p. 269]{Titchmarshfourier}%
\[
\varphi(z)=\int_{0}^{\infty}\left(  \int_{0}^{\infty}e^{-\left(  z+s\right)
t}dt\right)  g(s)ds=\int_{0}^{\infty}\frac{g(s)}{z+s}ds
\]
which is inverted as $g\left(  s\right)  =\frac{1}{\pi}\operatorname{Im}%
\left(  \varphi(se^{-i\pi})\right)  =-\frac{1}{\pi}\operatorname{Im}\left(
\varphi(se^{i\pi})\right)  $ for real $s$. Hence, if the Laplace transform and
its inverse are of the same form, then it holds that%
\begin{equation}%
\begin{array}
[c]{ccc}%
\varphi(z)=\mathcal{L}\left[  f(t)\right]  & \Longleftrightarrow & f\left(
z\right)  =\mathcal{L}\left[  \frac{1}{\pi}\operatorname{Im}\left(
\varphi(te^{-i\pi})\right)  \right]
\end{array}
\label{Gross_Laplace_inverse_pair}%
\end{equation}
and Gross \cite{Gross_1950} briefly states conditions on the validity of
(\ref{Gross_Laplace_inverse_pair}), which essentially relate to Theorem
\ref{theorem_entire_function_Cauchy}. Berberan-Santos' set
(\ref{inverse_Laplace_transform_Berberan_Santos}) of equations is essentially
a further development of Gross's Laplace pair
(\ref{Gross_Laplace_inverse_pair}).

\end{document}